\documentclass[12pt]{article}
\usepackage{amsfonts}
\usepackage{amssymb}
\usepackage{amsbsy}
\usepackage{amsthm}

\newtheorem{thm}{Theorem}[section]
\newtheorem{cor}[thm]{Corollary}
\newtheorem{lem}[thm]{Lemma}
\newtheorem{prop}[thm]{Proposition}
\theoremstyle{definition}
\newtheorem{defn}[thm]{Definition}
\theoremstyle{remark}
\newtheorem{rem}[thm]{Remark}

\newcommand{\bt}{\begin{thm}}
\newcommand{\et}{\end{thm}}
\newcommand{\bc}{\begin{cor}}
\newcommand{\ec}{\end{cor}}
\newcommand{\bl}{\begin{lem}}
\newcommand{\el}{\end{lem}}
\newcommand{\bp}{\begin{prop}}
\newcommand{\ep}{\end{prop}}
\newcommand{\bd}{\begin{defn}}
\newcommand{\ed}{\end{defn}}
\newcommand{\br}{\begin{rem}}
\newcommand{\er}{\end{rem}}
\newcommand{\bpr}{\begin{proof}}
\newcommand{\epr}{\end{proof}}

\newcommand{\bi}{\begin{itemize}}
\newcommand{\ei}{\end{itemize}}
\newcommand{\be}{\begin{enumerate}}
\newcommand{\ee}{\end{enumerate}}
\newcommand{\ds}{\displaystyle}
\newcommand{\ba}{\begin{array}}
\newcommand{\ea}{\end{array}}
\newcommand{\beq}{\begin{equation}}
\newcommand{\eeq}{\end{equation}}
\newcommand{\beqa}{\begin{eqnarray}}
\newcommand{\eeqa}{\end{eqnarray}}

\newcommand{\N}{{\mathbb N}}
\newcommand{\Z}{{\mathbb Z}}
\newcommand{\R}{{\mathbb R}}
\newcommand{\C}{{\mathbb C}}
\newcommand{\T}{{\mathbb T}}
\newcommand{\D}{{\mathbb D}}
\newcommand{\E}{{\mathbb E}}
\newcommand{\U}{{\mathbb U}}

\newcommand{\LL}{{\mathbb L}}
\newcommand{\B}{{\mathbb B}}

\newcommand{\cA}{{\mathcal  A}}
\newcommand{\cB}{{\mathcal  B}}
\newcommand{\cC}{{\mathcal  C}}
\newcommand{\cD}{{\mathcal  D}}
\newcommand{\cE}{{\mathcal  E}}
\newcommand{\cF}{{\mathcal  F}}
\newcommand{\cH}{{\mathcal  H}}
\newcommand{\cI}{{\mathcal  I}}
\newcommand{\cL}{{\mathcal  L}}
\newcommand{\cM}{{\mathcal  M}}

\newcommand{\cP}{{\mathcal  P}}

\newcommand{\cS}{{\mathcal  S}}
\newcommand{\cU}{{\mathcal  U}}
\newcommand{\cV}{{\mathcal  V}}

\newcommand{\cX}{{\mathcal  X}}

\newcommand{\frI}{{\mathfrak I}}

\newcommand{\frP}{{\mathfrak P}}

\newcommand{\bsa}{{\boldsymbol a}}
\newcommand{\bsb}{{\boldsymbol b}}
\newcommand{\bsalpha}{{\boldsymbol \alpha}}
\newcommand{\bsbeta}{{\boldsymbol \beta}}

\newcommand{\spn}{\mathrm{span}}
\newcommand{\supp}{\mathrm{supp}}

\newcommand{\re}{\mathrm{Re}}
\newcommand{\im}{\mathrm{Im}}

\begin{document}

\title{\bf Spectral methods for \break orthogonal rational functions}
\author{L. Vel\'azquez
\thanks{This work was partially realized during a stay of the author at the
Norwegian University of Science and Technology financed by Secretar\'{\i}a de Estado de
Universidades e Investigaci\'{o}n from the Ministry of Education and Science of Spain. The
work of the author was also partly supported by a research grant from the Ministry of
Education and Science of Spain, project code MTM2005-08648-C02-01, and by Project E-64
of Diputaci\'on General de Arag\'on (Spain).} \\
\small{Departamento de Matem\'atica Aplicada, Universidad de Zaragoza, Spain} \\
\small{\texttt{velazque@unizar.es}}
}
\maketitle

\kern-20pt

\begin{abstract}

An operator theoretic approach to orthogonal rational functions on the unit circle with
poles in its exterior is presented in this paper. This approach is based on the
identification of a suitable matrix representation of the multiplication operator
associated with the corresponding orthogonality measure. Two different alternatives are
discussed, depending whether we use for the matrix representation the standard basis of
orthogonal rational functions, or a new one with poles alternatively located in the
exterior and the interior of the unit circle. The corresponding representations are
linear fractional transformations with matrix coefficients acting respectively on
Hessenberg and five-diagonal unitary matrices.

In consequence, the orthogonality measure can be recovered from the spectral measure of
an infinite unitary matrix depending uniquely on the poles and the parameters of the
recurrence relation for the orthogonal rational functions. Besides, the zeros of the
orthogonal and para-orthogonal rational functions are identified as the eigenvalues of
matrix linear fractional transformations of finite Hessenberg and five-diagonal matrices.

As an application of this operator approach, we obtain new relations between the support
of the orthogonality measure and the location of the poles and parameters of the
recurrence relation, generalizing to the rational case known results for orthogonal
polynomials on the unit circle.

Finally, we extend these results to orthogonal polynomials on the real line with poles in
the lower half plane.

\end{abstract}


\medskip

\noindent{\it (2000) AMS Mathematics Subject Classification}: 42C05, 47B36.

\section{Introduction} \label{INT}

The connection with Jacobi matrices has led to numerous applications of spectral
techniques for self-adjoint operators in the theory of orthogonal polynomials on the real
line. The direct extension of these ideas to the orthogonal polynomials on the unit
circle yields a connection with unitary Hessenberg matrices (see
\cite{Al77,Ge44,Gr93,Si105,Te92}) which has provided some applications (see for instance
\cite{Go00,Go00b,Go00c,GoNeAs95,Te92}). Nevertheless, the authentic analogue of the
Jacobi matrices for the unit circle is a class of unitary five-diagonal matrices which
has been only recently discovered (see \cite{Wa93,FIVE}). This discovery has caused an
explosion of applications of spectral methods for unitary operators in the theory of
orthogonal polynomials on the unit circle, among which the numerous applications
appearing in the monograph \cite{Si105,Si205} have been only the starting point.

The orthogonal polynomials are a particular case of a more general kind of orthogonal
functions with interest in many pure and applied sciences: the orthogonal rational
functions with prescribed poles (see \cite{BGHN99} and the references therein). The
natural generalization of the orthogonal polynomials on the real line and the unit circle
requires the poles to be in the extended real line and in the exterior of the closed unit
disk respectively. The first situation presents special complications, an indication of
this being the fact that the poles can lie on the support of the orthogonality measure.
Indeed, considered as orthogonal rational functions, the main difference between the
orthogonal polynomials on the real line and the unit circle is not the location of the
support of the measure, but the relative location of the poles with respect to this
support. Actually, the Cayley transform maps the orthogonal rational functions on the
unit circle with poles in the exterior of the closed unit disk onto the orthogonal
rational functions on the real line with poles in the lower half plane, so both of them
can be thought as generalizations of the orthogonal polynomials on the unit circle. The
purpose of the paper is to generalize to this kind of orthogonal rational functions the
above referred spectral techniques for the orthogonal polynomials on the unit circle.

An important ingredient in the theory of orthogonal rational functions are the linear
fractional transformations $z \to (a_1z+a_2)(a_3z+a_4)^{-1}$ on the complex plane, where
$a_i$ are complex numbers. It is natural to expect the related spectral methods to have a
close relationship with the operator version of such transformations, i.e., the maps $T
\to (A_1T+A_2)(A_3T+A_4)^{-1}$ in the space of linear operators on a Hilbert space, where
the coefficients $A_i$ are now operators on the same Hilbert space. The theory of linear
fractional transformations with operator coefficients goes back to the work \cite{KrSm67}
of M. G. Krein and Yu. L. \v Smuljan, motivated by the study of operators in spaces with
an indefinite metric initiated by M. G. Krein in \cite{Kr50, Kr64}. As we will see, the
matrices related to the rational analogue of the orthogonal polynomials on the unit
circle are the result of applying a linear fractional transformation with matrix
coefficients to the Hessenberg and five-diagonal unitary matrices associated with the
polynomial case.

This reason, and also a better understanding of the subsequent rational generalizations,
motivates Section \ref{OP-ORF}, which summarizes the basics on spectral methods for
orthogonal polynomials on the unit circle and describes the main results needed about
orthogonal rational functions on the unit circle with poles in the exterior of the closed
unit disk. Section \ref{ORF-MT} introduces the operator linear fractional transformations
of interest for such orthogonal rational functions. The corresponding spectral theory is
developed in Sections \ref{ORF-HM} and \ref{ORF-5M}, which are devoted to the approaches
based on Hessenberg and five-diagonal matrices respectively. Section \ref{APPL} presents
some applications of the above spectral theory to the study of the relation between the
support of the orthogonality measure and the poles and parameters of the recurrence
relation for the orthogonal rational functions. Finally, the Appendix remarks the main
analogies and differences with the spectral theory for orthogonal rational functions on
the real line with poles lying on the lower half plane.

\section{OP and ORF on the unit circle} \label{OP-ORF}

In what follows a measure on the unit circle $\T=\{z\in\C:|z|=1\}$ will be a probability
Borel measure $\mu$ supported on a subset $\supp\mu$ of $\T$. Let $\mu$ be one of such
measures and consider the Hilbert space $L^2_\mu$ of $\mu$-square-integrable functions
with inner product
$$
\langle f,g \rangle_\mu = \int \overline{f(z)} g(z) \,d\mu(z), \qquad f,g \in L^2_\mu.
$$
Unless we say the opposite we will suppose that $\supp\mu$ is an infinite set. Then,
$(z^n)_{n\geq0}$ is a linearly independent subset of $L^2_\mu$ whose orthonormalization
gives the orthogonal polynomials (OP) $(\varphi_n)_{n\geq0}$ with respect to $\mu$. If we
choose these polynomials with positive leading coefficient, they satisfy the recurrence
relation \vskip-5pt
$$
\kern-205pt \varphi_0=1,
$$
\beq \label{rr}
\rho_n \pmatrix{\varphi_n \cr \varphi_n^*} =
\pmatrix{1 & a_n \cr \overline a_n & 1}
\pmatrix{z\varphi_{n-1} \cr \varphi_{n-1}^*},
\qquad n\geq1,
\eeq
$$
a_n={\varphi_n(0)\over\varphi_n^*(0)}, \qquad
\rho_n=\sqrt{1-|a_n|^2},
$$
where $\varphi_n^*(z)=z^n\overline \varphi_n(1/z)$ and $|a_n|<1$. This establishes a
bijection between measures $\mu$ on $\T$ and sequences $(a_n)_{n\geq1}$ in the unit disk
$\D=\{z\in\C:|z|<1\}$.

A central problem in the theory of OP on the unit circle is to find relations between the
orthogonality measure $\mu$ and the sequence $(a_n)_{n\geq1}$ appearing in the recurrence
relation for the OP. There are several approaches to this problem but these last years
have seen a rapid and impressive development of new operator theory techniques (see
\cite{Si105,Si205,Si} and references therein) based on the recent discovery of the
analogue for the unit circle of the Jacobi matrix related to OP on the real line (see
\cite{FIVE,Wa93}).

The main tool for the operator theoretic approach to the OP on $\T$ is
the unitary multiplication operator
$$
T_\mu \colon \mathop{L^2_\mu \to L^2_\mu} \limits_{f(z) \; \to \; zf(z)}
$$
It is known that the spectrum of $T_\mu$ coincides with $\supp\mu$ and the eigenvalues of
$T_\mu$, which have geometric multiplicity 1, are the mass points of $\mu$. The
eigenvectors of a given eigenvalue $\lambda$ are spanned by the characteristic function
$\cX_{\{\lambda\}}$ of the set $\{\lambda\}$. Even more, if $E$ is the spectral measure
of $T_\mu$ then $\mu(\Delta)=\langle 1,E(\Delta)1 \rangle_\mu$ for any Borel subset
$\Delta$ of $\T$. All these properties are true no matter whether $\supp\mu$ is finite or
infinite.

If $(f_n)_{n\geq0}$ is a basis of $L^2_\mu$, the matrix of $T_\mu$ with respect to
$(f_n)_{n\geq0}$ is the matrix $M$ whose $(i,j)$-th element is $M_{ij}=\langle f_i,T_\mu
f_j \rangle_\mu$. In other words,
\beq \label{M}
\pmatrix{zf_0(z) & zf_1(z) & \cdots} = \pmatrix{f_0(z) & f_1(z) & \cdots} M.
\eeq
Any matrix representation $M$ of $T_\mu$ can be identified with the unitary operator
$$
\mathop{\ell^2 \to \ell^2} \limits_{x \; \to \; Mx}
$$
on the space $\ell^2$ of square-sumable complex sequences. This operator is unitarily
equivalent to $T_\mu$. Therefore, once we know the dependence of $M$ on the parameters
$(a_n)_{n\geq1}$, this matrix permits us to recover the orthogonality measure $\mu$
starting from the recurrence relation of the OP. Regarding this problem, the utility of
the matrix representation $M$ depends on its simplicity as a function of the parameters
$(a_n)_{n\geq1}$.

For instance, when the polynomials are dense in $L^2_\mu$, the representation of $T_\mu$
with respect to the OP $(\varphi_n)_{n\geq0}$ is the irreducible Hessenberg matrix (see
\cite{Al77,Ge44,Gr93,Si105,Te92})
\beq \label{HESS}
\cH = \pmatrix{
-a_1 & -\rho_1 a_2 & - \rho_1 \rho_2 a_3 & - \rho_1 \rho_2 \rho_3 a_4 & \cdots
\cr
\rho_1 & -\overline a_1 a_2 & - \overline a_1 \rho_2 a_3 & - \overline a_1 \rho_2 \rho_3 a_4 & \cdots
\cr
0 & \rho_2 & -\overline a_2 a_3 & - \overline a_2 \rho_3 a_4 & \cdots
\cr
0 & 0 & \rho_3 & -\overline a_3 a_4 & \cdots
\cr
\cdots & \cdots & \cdots & \cdots & \cdots}.
\eeq
$\cH=(h_{i,j})$ is called a Hessenberg matrix because $h_{i,j}=0$ for $i>j+1$, and the
irreducibility means that $h_{j+1,j}\neq0$ for any $j$. Using the $2\times2$ symmetric
unitary matrices
\beq \label{THETA}
\Theta_n = \pmatrix{-a_n & \rho_n \cr \rho_n & \overline a_n}, \qquad n\geq1,
\eeq
the Hessenberg representation can be factorized as
\beq \label{FAC1} \kern-7pt
\cH = \lim_n \pmatrix{\Theta_1 & \cr & I} \kern-3pt
\pmatrix{I_1 && \cr & \Theta_2 & \cr && I} \kern-3.5pt
\pmatrix{I_2 && \cr & \Theta_3 & \cr && I} \kern-1pt \cdots
\pmatrix{I_{n-1} && \cr & \Theta_n & \cr && I},
\eeq
where $I$ and $I_n$ mean the identity matrix of order infinite and $n$ respectively and
the limit has to be understood in the strong sense.

Apart from its complexity, the Hessenberg representation has the inconvenience of being
valid only when the polynomials are dense in $L^2_\mu$. In the general case $\cH$ is a
matrix representation of the restriction $T_\mu \upharpoonright \cP \colon \cP \to \cP$
of $T_\mu$ to the $T_\mu$-invariant subspace given by the closure $\cP$ of the
polynomials in $L^2_\mu$. As a restriction of a unitary operator, $T_\mu \upharpoonright
\cP$ is isometric but not necessarily unitary. $\cH$ is a representation of $T_\mu$ iff
any of the following equivalent conditions hold (see \cite{Ge44,Si105}):
$$
\cP=L^2_\mu \Leftrightarrow \log\mu' \notin L^1_m \Leftrightarrow
(a_n)_{n\geq1}\notin\ell^2 \Leftrightarrow \cH \hbox{ is unitary}.
$$
We denote by $m$ the Lebesgue measure on $\T$.

A way to avoid the problems of the Hessenberg representation is to use as a basis of
$L^2_\mu$ the Laurent OP $(\chi_n)_{n\geq0}$ that arise from the orthonormalization of
$(1,z,z^{-1},z^2,z^{-2},\dots)$, which are given by (see \cite{FIVE,Si105,Th88,Wa93})
$$
\chi_{2n}(z)=z^{-n}\varphi_{2n}^*(z), \qquad \chi_{2n+1}(z)=z^{-n}\varphi_{2n+1}(z),
\qquad n\geq0.
$$
The corresponding representation of $T_\mu$ is the five-diagonal matrix
(see \cite{FIVE,Si105,Wa93})
\beq \label{CMV}
\kern-5pt
\cC = \pmatrix{
-a_1 & -\rho_1 a_2 & \rho_1 \rho_2 & 0 & 0 & 0 & 0 & \cdots
\cr
\rho_1 & -\overline a_1 a_2 & \overline a_1 \rho_2 & 0 & 0 & 0 & 0 & \cdots
\cr
0 & -\rho_2 a_3 & -\overline a_2 a_3 & -\rho_3 a_4 & \rho_3 \rho_4 & 0 & 0 & \cdots
\cr
0 & \rho_2 \rho_3 & \overline a_2 \rho_3 & -\overline a_3 a_4 & \overline a_3 \rho_4 & 0 & 0 & \cdots
\cr
0 & 0 & 0 & -\rho_4 a_5 & -\overline a_4 a_5 & -\rho_5 a_6 & \rho_5 \rho_6 & \cdots
\cr
0 & 0 & 0 & \rho_4 \rho_5 & \overline a_4 \rho_5 & -\overline a_5 a_6 & \overline a_5 \rho_6 & \cdots
\cr
0 & 0 & 0 & 0 & 0 & -\rho_6 a_7 & -\overline a_6 a_7 & \cdots
\cr
\cdots & \cdots & \cdots & \cdots & \cdots & \cdots & \cdots & \cdots},
\eeq
which, apart from being valid for any measure $\mu$ on $\T$, is a band instead of a
Hessenberg matrix. Also, it has a much simpler dependence on the parameters
$(a_n)_{n\geq1}$ than in the Hessenberg case. Moreover, this five-diagonal representation has a much
better factorization than the Hessenberg one since $\cC=\cC_o\cC_e$, where $\cC_o$ and
$\cC_e$ are the $2 \times 2$-block-diagonal symmetric unitary matrices
\beq \label{FAC2}
\cC_o=\pmatrix{\Theta_1 &&& \cr & \Theta_3 && \cr && \Theta_5 & \cr &&& \ddots}, \qquad
\cC_e=\pmatrix{I_1 &&& \cr & \Theta_2 && \cr && \Theta_4 & \cr &&& \ddots}.
\eeq

Alternatively, it is possible to orthonormalize $(1,z^{-1},z,z^{-2},z^2,\dots)$. This
leads to the Laurent OP $(\chi_{n*})_{n\geq0}$ where $\chi_{n*}(z)=\overline\chi_n(1/z)$,
i.e.,
$$
\chi_{2n*}(z)=z^{-n}\varphi_{2n}(z), \qquad \chi_{2n+1*}(z)=z^{-n-1}\varphi_{2n+1}^*(z),
\qquad n\geq0.
$$
The related representation of $T_\mu$ is simply the transposed matrix $\cC^T=\cC_e\cC_o$
of $\cC$.

The Hessenberg and five-diagonal matrices given in (\ref{HESS}) and (\ref{CMV}) represent
a multiplication operator (on $\cP$ or $L^2_\mu$) only when $(a_n)_{n\geq1}$ lies on
$\D$. Nevertheless, they are well defined matrices for any sequence $(a_n)_{n\geq1}$ in
the closed unit disk $\overline\D$. Indeed, factorizations (\ref{FAC1}) and (\ref{FAC2})
show that, even in this case, the Hessenberg representation is isometric while the
five-diagonal one is unitary. Furthermore, when some $a_n\in\T$, these Hessenberg and
five-diagonal matrices decompose as a direct sum of an $n \times n$ and an infinite
matrix. This decomposition property is of interest when trying to make perturbative
spectral analysis of such matrix representations.

The Hessenberg and five-diagonal representations of $T_\mu$ also give a spectral
interpretation for the zeros of the OP in terms of the parameters of the recurrence. This
result comes from the relation between the OP and certain orthogonal truncations of the
operator $T_\mu$. The restriction $T_\mu \upharpoonright \cP_{n,l}$ of the multiplication
operator $T_\mu$ to the subspace $\cP_{n,l}=\spn\{z^l,z^{l+1},\dots,z^{l+n-1}\}$ has no
sense since $\cP_{n,l}$ is not invariant under $T_\mu$. To give sense to this kind of
restriction we must multiply $T_\mu$ on the left by a projection on $\cP_{n,l}$. In
particular, if $P_{n,l} \colon L^2_\mu \to L^2_\mu$ is the orthogonal projection on
$\cP_{n,l}$, the operator $T_\mu^{(\cP_{n,l})}=P_{n,l} T_\mu \upharpoonright \cP_{n,l}$ is
called the orthogonal truncation of $T_\mu$ on $\cP_{n,l}$. The key point is that, for
any $l\in\Z$, the characteristic polynomial of $T_\mu^{(\cP_{n,l})}$ is, up to factors, the
$n$-th OP $\varphi_n$ (see \cite{Si105}).

The first $n$ OP $(\varphi_k)_{k=0}^{n-1}$ are a basis of $\cP_n=\cP_{n,0}$ and the
related matrix of $T_\mu^{(\cP_n)}$ is the principal submatrix $\cH_n$ of $\cH$ of order
$n$. So, $\varphi_n$ is proportional to the characteristic polynomial of $\cH_n$, whose
eigenvalues are therefore the zeros of $\varphi_n$. Furthermore, for $l=-[(n-1)/2]$, the
first $n$ Laurent OP $(\chi_k)_{k=0}^{n-1}$ are a basis of $\cP_{n,l}$ and the
corresponding matrix of $T_\mu^{(\cP_{n,l})}$ is the principal submatrix $\cC_n$ of $\cC$
of order $n$. Hence, $\varphi_n$ is proportional to the characteristic polynomial of $\cC_n$
and, thus, the zeros of $\varphi_n$ are the eigenvalues of $\cC_n$.

Contrary to the full infinite matrix, $\cH_n$ and $\cC_n$ are not unitary and depend only
on the first $n$ parameters $a_1,\dots,a_n$. However, factorizations (\ref{FAC1}) and
(\ref{FAC2}) show that if we change in these principal submatrices the last parameter
$a_n\in\D$ by a complex number $u\in\T$, then we obtain a unitary matrix. The
corresponding characteristic polynomial is the result of performing $n$ steps of
recurrence (\ref{rr}), but substituting in the last one $a_n\in\D$ by $u\in\T$, i.e., it
is a multiple of
$$
z\varphi_{n-1}(z)+u\varphi_{n-1}^*(z).
$$
Using (\ref{rr}), this polynomial can be alternatively written up to factors as
$$
\varphi_n(z)+v\varphi_n^*(z), \qquad v = \frac{u-a}{1-\overline a_nu},
$$
and the arbitrariness of $u\in\T$ translates into a similar arbitrariness for $v\in\T$.
These polynomials, called para-orthogonal polynomials (POP), were introduced for the
first time in \cite{JoNjTh89}. There it was proved that such POP have simple zeros lying
on $\T$, which play the role of nodes in the Szeg\H{o} quadrature formulas on $\T$ (the
analogue of the Gaussian quadrature formulas on $\R$), thus, providing finitely supported
measures on $\T$ that $*$-weakly converge to the measure $\mu$. Therefore, the nodes of
the Szeg\H o quadrature formulas can be obtained as eigenvalues of Hessenberg or
five-diagonal unitary matrices.

\medskip

Our aim is to generalize the above results to the orthogonal rational functions with
poles outside of the support of the orthogonality measure. Two archetypical situations
will be considered: measures on the unit circle $\T$ and measures on the extended real
line $\overline\R=\R\cup\{\infty\}$. For convenience, the analysis will be done in a
detailed way for measures on the unit circle, the discussion of the special features in
the case of the real line being relegated to the Appendix. So, for the moment we will
consider a measure $\mu$ on $\T$ and the corresponding orthogonal rational functions with
poles arbitrary located in the exterior of the unit circle
$\E=\overline\C\setminus\overline\D$. We consider the extended complex plane
$\overline\C=\C\cup\{\infty\}$ to include for the poles the possibility of being located
at $\infty$. Indeed, the OP with respect to $\mu$ correspond to the special case of the
orthogonal rational functions with all the poles at $\infty$.

An important transformation in $\overline\C$ is $\hat z = 1/\overline z$, which leaves
invariant any element of $\T$ and establishes a bijection between $\D$ and $\E$. This
transformation induces the $*$-involution $f_*(z)=\overline{f(\hat z)}$ in the set of
complex functions, which defines an anti-unitary operator on $L^2_\mu$ for any measure
$\mu$ on $\T$. As a consequence, a sequence $(f_n)_{n\geq0}$ of functions is a basis of
$L^2_\mu$ iff $(f_{n*})_{n\geq0}$ is a basis too. Moreover, the $*$-involution on
(\ref{M}) gives
$$
\pmatrix{zf_{0*}(z) & zf_{1*}(z) & \cdots} = \pmatrix{f_{0*}(z) & f_{1*}(z) & \cdots}
\overline M {}^{-1},
$$
which, taking into account that $M$ is unitary, shows that the matrix of $T_\mu$ with
respect to $(f_{n*})_{n\geq0}$ is the transposed $M^T$ of the matrix $M$ associated with
$(f_n)_{n\geq0}$. This relation holds when $\mu$ is finitely supported too, with the only
difference that the basis of $L^2_\mu$ are finite.

Another essential ingredient in the theory of orthogonal rational functions on $\T$ are
the M\"{o}bius transformations $\zeta_\alpha$ defined for any $\alpha\in\D$ by
$$
\zeta_\alpha(z) = \frac{\varpi_\alpha^*(z)}{\varpi_\alpha(z)}, \qquad
\cases{
\varpi_\alpha(z) =
1 - \overline\alpha z,
\medskip \cr
\varpi_\alpha^*(z) = z\varpi_{\alpha*}(z) = z - \alpha.
}
$$
Up to factors in $\T$, they are all the automorphisms of $\D$. Indeed, $\zeta_\alpha$ is a
bijection of $\overline\C$ onto $\overline\C$ that leaves invariant $\T$, $\D$ and $\E$.
The inverse transformation of $\zeta_\alpha$ is $\tilde\zeta_\alpha=\zeta_{-\alpha}$. It
is also remarkable that $\zeta_{\alpha*}=1/\zeta_\alpha$. We distinguish the value
$\alpha_0=0$ that gives $\zeta_{\alpha_0}(z)=z$.

To get rational functions with fixed poles in $\E$ we introduce a sequence
$(\alpha_n)_{n\geq1}$ in $\D$. This sequence defines the finite Blaschke
products $(B_n)_{n\geq0}$ given by
\beq \label{BLAS}
\ba{l}
B_0=1,
\smallskip \cr
B_n=\zeta_{\alpha_1}\cdots\zeta_{\alpha_n}, \qquad n\geq1. \ea
\eeq
Notice that $B_{n*}=1/B_n$. The subspace
$$
\cL_n=\spn\{B_0,B_1,\dots,B_{n-1}\}=
\frac{\cP_n}{\varpi_{\alpha_1}\cdots\,\varpi_{\alpha_{n-1}}}
$$
consists of those rational functions whose poles, counted with multiplicity, lie on
$(\hat\alpha_k)_{k=1}^{n-1}$. We use the notation
$\cL_\infty=\spn\{B_n\}_{n\geq0}=\cup_{n\geq1}\cL_n$ for the set of rational functions
with poles lying on $(\hat\alpha_n)_{n\geq1}$, counted with multiplicity, and $\cL$ for
the closure of $\cL_\infty$ in $L^2_\mu$.

If $\mu$ is a measure on $\T$ we can consider the rational functions $(\Phi_n)_{n\geq0}$
that arise from the orthonormalization of $(B_n)_{n\geq0}$ in $L^2_\mu$.
$(\Phi_n)_{n\geq0}$ are called orthogonal rational functions (ORF) with respect to $\mu$
associated with $(\alpha_n)_{n\geq1}$. When referring to $(\Phi_n)_{n\geq0}$ we will call
it in short a sequence of ORF on the unit circle. These functions satisfy a recurrence
relation which, with an appropriate normalization of $(\Phi_n)_{n\geq0}$, has the form
(see \cite[Theorem 4.1.3]{BGHN99})
\beq \label{RR0}
\ba{l}
\Phi_0=1,
\smallskip \cr \ds
\pmatrix{\Phi_n \cr \Phi_n^*} = e_n {\varpi_{n-1}\over\varpi_n}
\pmatrix{1 & b_n \cr \overline b_n & 1}
\pmatrix{z_n \zeta_{n-1} \Phi_{n-1} \cr \Phi_{n-1}^*},
\qquad n\geq1,
\ea
\eeq
where
$$
b_n={\Phi_n(\alpha_{n-1})\over\Phi_n^*(\alpha_{n-1})}, \kern12pt
z_n=\cases{-\frac{|\alpha_n|}{\alpha_n} & if $\alpha_n\neq0$, \cr 1 & if $\alpha_n=0$,} \kern10pt
e_n=\sqrt{{\varpi_n(\alpha_n)\over\varpi_{n-1}(\alpha_{n-1})}{1 \over 1-|b_n|^2}},
$$
and we use the notation
$$
\zeta_n=\zeta_{\alpha_n}, \quad
\varpi_n=\varpi_{\alpha_n}, \quad
\varpi_n^*=\varpi_{\alpha_n}^*, \quad
\Phi_n^* = z_1z_2 \cdots z_nB_n\Phi_{n*}.
$$
Notice that we do not follow the standard notation $\zeta_n=z_n\zeta_{\alpha_n}$ and $B_n
= z_1\zeta_{\alpha_1} \cdots z_n\zeta_{\alpha_n}$ (see for instance \cite{BGHN99}). In
fact, concerning the matrix representations of the multiplication operator, it is more
convenient to avoid the presence of the factors $z_n$ in recurrence (\ref{RR0}),
something that we can get using the ORF $(\phi_n)_{n\geq0}$ given by
$$
\ba{l}
\phi_0=1,
\smallskip \cr
\phi_n = \overline z_1 \overline z_2 \cdots \overline z_n \Phi_n, \qquad n\geq1,
\ea
$$
and defining the superstar operation omitting the factors $z_n$, that is,
$$
\phi_n^* = B_n \phi_{n*}.
$$
Then, (\ref{RR0}) is equivalent to
\beq \label{RR}
\ba{l}
\phi_0=1,
\smallskip \cr \ds
\pmatrix{\phi_n \cr \phi_n^*} = e_n {\varpi_{n-1}\over\varpi_n}
\pmatrix{1 & a_n \cr \overline a_n & 1}
\pmatrix{\zeta_{n-1} \phi_{n-1} \cr \phi_{n-1}^*},
\qquad n\geq1,
\ea
\eeq
with
$$
a_n={\phi_n(\alpha_{n-1})\over\phi_n^*(\alpha_{n-1})}=
\overline z_1 \overline z_2 \cdots \overline z_n b_n, \qquad
e_n=\sqrt{{\varpi_n(\alpha_n)\over\varpi_{n-1}(\alpha_{n-1})}{1 \over 1-|a_n|^2}}.
$$
In the polynomial case, corresponding to $\alpha_n=0$ for all $n$, (\ref{RR}) gives
exactly (\ref{rr}). As in the polynomial situation, the parameters $(a_n)_{n\geq1}$ of
(\ref{RR}) lie on $\D$. A Favard-type theorem also holds (see \cite[Theorem
8.1.4]{BGHN99}): given a sequence $(a_n)_{n\geq1}$ in $\D$, the functions
$(\phi_n)_{n\geq0}$ defined by recurrence (\ref{RR}) are orthonormal with respect to some
measure on $\T$. This measure is unique when the infinite Blaschke product
$B(z)=\prod_{n=1}^\infty\zeta_n(z)$ diverges to zero for $z\in\D$, i.e., when
$\sum_{n=1}^\infty (1-|\alpha_n|) = \infty$. This condition means that the sequence
$(\alpha_n)_{n\geq1}$ can not approach to $\T$ very quickly.

Notice that, given a measure $\mu$ on $\T$ and a sequence $(\alpha_n)_{n\geq1}$ in $\D$,
the parameters $(a_n)_{n\geq1}$ are uniquely defined. To see this, suppose that
$(\hat\phi_n)_{n\geq0}$ is another sequence of ORF satisfying a recurrence like
(\ref{RR}), but with parameters $(\hat a_n)_{n\geq1}$ instead of $(a_n)_{n\geq1}$. Then,
$\hat\phi_n=\epsilon_n\phi_n$ with $\epsilon_n\in\T$ and $\epsilon_0=1$. Hence, comparing
the recurrences for $(\phi_n)_{n\geq0}$ and $(\hat\phi_n)_{n\geq0}$ gives
$$
\frac{1}{\sqrt{1-|a_n|^2}}
\pmatrix{\epsilon_n & 0 \cr 0 & \overline\epsilon_n}
\pmatrix{1 & a_n \cr \overline a_n & 1}
\pmatrix{\overline\epsilon_{n-1} & 0 \cr 0 & \epsilon_{n-1}} =
\frac{1}{\sqrt{1-|\hat a_n|^2}}
\pmatrix{1 & \hat a_n \cr \overline{\hat a}_n & 1}.
$$
Taking determinants in both sides of the above equality we obtain $|\hat a_n|=|a_n|$.
In consequence $\epsilon_n=\epsilon_{n-1}$ for $n\geq1$, which yields
$\epsilon_n=\epsilon_0=1$. Therefore, $\hat a_n = a_n$ and $\hat\phi_n=\phi_n$.

The above results show that any sequence $\bsalpha=(\alpha_n)_{n\geq1}$ in $\D$ defines a
surjective application
$$
\cS_\bsalpha \colon
\mathop{\kern-40pt \frP \longrightarrow \D^\infty}
\limits_{\kern5pt \ds \mu \longrightarrow \bsa=(a_n)_{n\geq1}}
$$
between the set $\frP$ of probability measures on $\T$ and the set $\D^\infty$ of
sequences in $\D$. Furthermore, $\cS_\bsalpha$ is a bijection when
$\sum_{n=1}^\infty(1-|\alpha_n|) = \infty$. The study of the application $\cS_\bsalpha$
is one of the main interests to find a matrix representation of the multiplication
operator $T_\mu$ with a simple dependence on the parameters
$\bsalpha=(\alpha_n)_{n\geq1}$ and $\bsa=(a_n)_{n\geq1}$. Indeed, in the polynomial case,
corresponding to $\bsalpha=0$, the five-diagonal representation $\cC=\cC(\bsa)$ of
$T_\mu$ given in (\ref{CMV}) has revealed to be a powerful tool in the study of $\cS_0$.

To find such a matrix representation, it is convenient to write recurrence (\ref{RR}) in
a different way. For any $\alpha\in\D$ we can define the positive number
$$
\eta_\alpha = \varpi_\alpha(\alpha)^{1/2} = \sqrt{1-|\alpha|^2}.
$$
Denoting $\eta_n=\eta_{\alpha_n}$ and introducing the parameters
\beq \label{RHO}
\rho_n=\sqrt{1-|a_n|^2}, \qquad
\rho_n^+=\frac{\eta_{n-1}}{\eta_n}\,\rho_n, \qquad
\rho_n^-=\frac{\eta_n}{\eta_{n-1}}\,\rho_n,
\eeq
(\ref{RR}) yields
\beq \label{RR1}
\cases{
\varpi_{n-1}^*\phi_{n-1} = \rho_n^+\varpi_n\phi_n - a_n\varpi_{n-1}\phi_{n-1}^*,
\medskip \cr
\varpi_n\phi_n^* = \overline a_n\varpi_n\phi_n + \rho_n^-\varpi_{n-1}\phi_{n-1}^*,
}
\qquad n\geq1.
\eeq
This way of writing (\ref{RR}) will be useful later.

\section{Operator M\"{o}bius transformations} \label{ORF-MT}

As we will see, the operator version of the scalar M\"{o}bius transformations $\zeta_\alpha$
appears in a natural way in the spectral theory of ORF on the unit circle. Analogously to
the scalar case, such operator M\"{o}bius transformations are a particular case of the linear
fractional transformations with operator coefficients introduced by M. G. Krein in
\cite{Kr50, Kr64} for the study of spaces with an indefinite metric. A detailed study of
these operator M\"{o}bius transformations in the general context of linear fractional
transformations can be found, for instance, in the original paper of M. G. Krein and Yu
L. \v Smuljan \cite{KrSm67} or in the most recent survey of T. Ya Azizov and I. S.
Iokhvidov \cite{AzIo89} and the references therein. We will introduce the operator M\"{o}bius
transformations summarizing the main properties of interest for us.

Before doing this, we will fix some notations and conventions for linear operators. In
what follows $(H,(\cdot,\cdot))$ means a separable Hilbert space. Given a linear operator
$T$ on $H$, $T^\dag$ denotes its adjoint, $\sigma(T)$ its spectrum and $\sigma_p(T)$ its
point spectrum. As usual, we omit the identity operator $\boldsymbol 1$ on $H$ so we use
the same symbol $z$ for the complex number $z\in\C$ and for the operator $z \boldsymbol
1$, the meaning being clear from the context in any case. In general, we will deal with
the Banach space $(\B_H,\|\cdot\|)$ of everywhere defined bounded linear operators on
$H$.

In particular, $\B_{\C^n}$ and $\B_{\ell^2}$ can be identified with the sets of $n \times
n$ complex matrices and infinite bounded complex matrices respectively. In this
identification we associate any bounded square matrix $M$ with the operator $x \to Mx$,
where $x$ is a column vector of $\C^n$ or $\ell^2$. However, we could also consider the
operator $x \to xM$, where $x$ is a row vector of $\C^n$ or $\ell^2$. Both operators have
the same spectrum, although their eigenvalues can be different in the case of $\ell^2$.
Nevertheless, we will normally work with normal or finite-dimensional matrices, for which
the eigenvalues are the same in both situations. However, even in these cases, the
eigenvectors are in general different. So, we will distinguish between right eigenvectors
(or just eigenvectors) for $x \to Mx$ and left eigenvectors for $x \to xM$. That is,
right eigenvectors are the standard ones while left eigenvectors are the transposed of
the eigenvectors of $M^T$ (in particular, when $M$ is normal, right eigenvectors are the
adjoints of left eigenvectors). In the subsequent discussions, this convention often
permits us to avoid the ${}^T$ superindex, something convenient because many indices
appear later.

The operator M\"{o}bius transformations on $H$ are linear fractional transformations with
operator coefficients that transform bijectively the unit ball
$\D_H=\{T\in\B_H:\|T\|<1\}$ of $\B_H$ onto itself. The role of the complex parameter
$\alpha\in\D$ of $\zeta_\alpha$ is played by an operator $A\in\D_H$, so that
$$
\eta_A = \sqrt{1 - A\kern1ptA^\dag}
$$
defines a positive operator with bounded inverse. Therefore, for any operator $T$ in the
closed unit ball $\overline\D_H=\{T\in\B_H:\|T\|\leq1\}$ we can define the operators
$\zeta_A(T),\tilde\zeta_A(T)\in\B_H$ by
$$
\ba{l}
\zeta_A(T) = \eta_A \,\varpi_A(T)^{-1} \varpi_A^*(T) \,\eta_{A^\dag}^{-1}, \qquad
\cases{
\varpi_A(T) = 1 - TA^\dag,
\medskip \cr
\varpi_A^*(T) = T - A,
}
\medskip \cr
\tilde\zeta_A(T) = \eta_A^{-1} \tilde\varpi_A^*(T) \,\tilde\varpi_A(T)^{-1} \eta_{A^\dag},
\kern22pt
\cases{
\tilde\varpi_A(T) = 1 + A^\dag T,
\medskip \cr
\tilde\varpi_A^*(T) = T + A.
}
\ea
$$
As in the scalar case, $\eta_A=\varpi_A(A)^{1/2}$. As we will see, the spectral theory of
ORF is related to transformations $\zeta_A,\tilde\zeta_A$ with $A$ normal, so that
$\eta_{A^\dag}=\eta_A$ in such a case.

The transformations $\zeta_A$ and $\tilde\zeta_A$ are the operator analogs of the scalar
M\"{o}bius transformations $\zeta_\alpha $ and $\tilde\zeta_\alpha$ respectively. The factors
$\eta_A,\eta_{A^\dag}$ disappear in the scalar case due to the commutativity.
Nevertheless, these factors are necessary for these operator transformations to keep
similar properties to the scalar ones. Actually, $\zeta_A$ and $\tilde\zeta_A$ map
$\overline\D_H$ on $\overline\D_H$, as follows from the identities
\beq \label{INV}
\ba{l}
\varpi_A(T)\,\eta_A^{-1}(1-\zeta_A(T)\,\zeta_A(T)^\dag)\,\eta_A^{-1}\varpi_A(T)^\dag =
1-TT^\dag,
\medskip \cr
\tilde\varpi_A(T)^\dag\,\eta_{A^\dag}^{-1}(1-\tilde\zeta_A(T)^\dag\,\tilde\zeta_A(T))
\,\eta_{A^\dag}^{-1}\tilde\varpi_A(T) = 1-T^\dag T. \ea \eeq Besides, for any
$S,T\in\overline\D_H$, a direct calculation shows that $S=\zeta_A(T)$ iff
$T=\tilde\zeta_A(S)$, so $\zeta_A$ and $\tilde\zeta_A$ are mutually inverse
transformations that map $\overline\D_H$ onto itself. Furthermore, (\ref{INV}) also
proves that $\zeta_A$ and $\tilde\zeta_A$ leave invariant $\D_H$ and
$\T_H=\{T\in\B_H:\|T\|=1\}$, mapping onto itself the set of isometries as well as the set
of unitary operators on $H$. Indeed, as it was proven in \cite{KrSm67}, up to unitary
left and right factors, these operator M\"{o}bius transformations are the only linear
fractional transformations with operator coefficients mapping bijectively $\D_H$ onto
itself.

Using the relation $\eta_A^2A=A\eta_{A^\dag}^2$ it is straightforward to verify the
identities
\beq \label{DAG}
\zeta_A(T)^\dag=\zeta_{A^\dag}(T^\dag), \qquad
\tilde\zeta_A(T)^\dag=\tilde\zeta_{A^\dag}(T^\dag),
\eeq
which imply that $\tilde\zeta_A(T)=\tilde\zeta_{A^\dag}(T^\dag)^\dag=\zeta_{-A}(T)$ as in
the scalar case. Notice that the equalities $\zeta_A=\tilde\zeta_{-A}$ and
$\tilde\zeta_A=\zeta_{-A}$ provide alternative expressions for $\zeta_A$ and
$\tilde\zeta_A$.

Some formulas for the operator M\"{o}bius transformations will be of interest. From the
relations $(\eta_A^2)^n A = A (\eta_{A^\dag}^2)^n$ for $n=0,1,2,\dots$, and using the
functional calculus for self-adjoint operators, we find that
$$
\eta_A \, A = A \, \eta_{A^\dag}.
$$
Thus, if we define
$$
T_A = \eta_A^{-1} \,T \,\eta_{A^\dag}
$$
for any linear operator $T$ on $H$, then, for all $T\in\overline\D_H$,
\beq \label{MOB2}
\zeta_A(T_A) = \varpi_A(T)^{-1} \varpi_A^*(T),
\qquad
\tilde\zeta_A(T)= \tilde\varpi_A^*(T_A) \,\tilde\varpi_A(T_A)^{-1}.
\eeq
This, together with the immediate identity
\beq \label{MOB1}
\varpi_A^*(T)-\varpi_A(T)\,S = T\,\tilde\varpi_A(S)-\tilde\varpi_A^*(S),
\eeq
yields
\beq \label{MOB}
\varpi_A(T) \,(\zeta_A(T_A)-S_A) = (T-\tilde\zeta_A(S)) \,\tilde\varpi_A(S_A)
\eeq
for all $T,S\in\overline\D_H$. Substituting $S$ by $\zeta_A(S)$ in (\ref{MOB})
gives
\beq \label{MOB'}
T-S = \varpi_A(T)\,\eta_A^{-1}(\zeta_A(T)-\zeta_A(S))\,\eta_{A^\dag}^{-1}\,\tilde\varpi_{-A}(S),
\eeq
where we have used that
$$
\tilde\varpi_A(\zeta_A(S_A)) =
\tilde\varpi_A(\tilde\zeta_{-A}(S_A)) = \tilde\varpi_{-A}(S)^{-1} \eta_{A^\dag}^2.
$$
If we take $A=\alpha$ and $S=z$ with $\alpha\in\D$ and $z\in\overline\D$, (\ref{MOB'}) becomes
\beq \label{MOB3}
z-T = \frac{\varpi_\alpha(z)}{\varpi_\alpha(\alpha)}
\,(\zeta_\alpha(z)-\zeta_\alpha(T)) \,\varpi_\alpha(T).
\eeq
In particular, choosing $T=\lambda$ with $\lambda\in\overline\D$,
\beq \label{MOB4}
\zeta_\alpha(z)-\zeta_\alpha(\lambda) =
\frac{\varpi_\alpha(\alpha)}{\varpi_\alpha(z)\,\varpi_\alpha(\lambda)} \,(z-\lambda).
\eeq
Notice that (\ref{MOB3}) and (\ref{MOB4}) actually hold for any
$z,\lambda\in\C\setminus\{\hat\alpha\}$.

\section{ORF and Hessenberg matrices} \label{ORF-HM}

In this section we will prove that the orthogonality measure of a sequence of ORF, as
well as the zeros of the ORF, have a spectral interpretation in terms of Hessenberg
matrices. Our first aim is to find the matrix representation of a unitary multiplication
operator with respect to a basis of ORF. Before stating the result, let us see which kind
of matrix representation we can expect. Let $\nu$ be a measure on $\T$ and
$(\varphi_n)_{n\geq0}$ the corresponding OP with positive leading coefficient. Given
$\alpha\in\D$, the functions $\phi_n(z)=\varphi_n(\zeta_\alpha(z))$ define a sequence
$(\phi_n)_{n\geq0}$ of ORF with fixed poles at $\hat\alpha$. The corresponding
orthogonality measure is $\mu=\nu_\alpha$, where
$\nu_\alpha(\Delta)=\nu(\zeta_\alpha(\Delta))$ for any Borel subset $\Delta$ of $\T$. It
is straightforward to see that recurrence (\ref{rr}) for $(\varphi_n)_{n\geq0}$ is
rewritten in terms of $(\phi_n)_{n\geq0}$ as recurrence (\ref{RR}) with the same
parameters $\bsa=(a_n)_{n\geq1}$, i.e., $\bsa=\cS_0(\nu)=\cS_\bsalpha(\mu)$. The matrix
representation of the isometric operator $T_\nu \upharpoonright \cP$ with respect to
$(\varphi_n)_{n\geq0}$ is a Hessenberg matrix $\cH=\cH(\bsa)$ with the form (\ref{HESS}).
Therefore, the matrix of $T_\mu \upharpoonright \cL$ with respect to the ORF
$(\phi_n)_{n\geq0}$ is
$$
( \langle \phi_i(z),z\phi_j(z) \rangle_\mu )_{i,j=0}^\infty =
( \langle \varphi_i(z),\tilde\zeta_\alpha(z)\varphi_j(z) \rangle_\nu )_{i,j=0}^\infty =
\tilde\zeta_\alpha(\cH).
$$
The following theorem is a natural generalization of this particular situation.

\bt \label{HESSORF}

Let $\bsalpha=(\alpha_n)_{n\geq1}$ be compactly included in $\D$, $\mu$ a measure on $\T$
and $\bsa=(a_n)_{n\geq1}=\cS_\bsalpha(\mu)$. Then, $\cL$ is $T_\mu$-invariant and the
matrix of the isometric operator $T_\mu \upharpoonright \cL$ with respect to the
corresponding ORF $(\phi_n)_{n\geq0}$ is $\cV=\tilde\zeta_\cA(\cH)$, where
$\cH=\cH(\bsa)$ is given in (\ref{HESS}) and
$$
\cA = \cA(\bsalpha) = \pmatrix{\alpha_0 &&& \cr & \alpha_1 && \cr && \alpha_2 & \cr &&& \ddots}.
$$
The isometric matrix $\cV$ represents the full operator $T_\mu$ iff any of the following equivalent
conditions is fulfilled:
$$
\cL=L^2_\mu \Leftrightarrow \cP=L^2_\mu \Leftrightarrow \log\mu' \notin L^1_m
\Leftrightarrow \bsa\notin\ell^2 \Leftrightarrow \cV \hbox{ is unitary}.
$$

\et

\bpr

$\|\cA\|<1$ because $\bsalpha$ is compactly included in $\D$, thus $\tilde\zeta_\cA$ maps
onto theirselves the sets of infinite isometric and unitary matrices. Therefore, taking
into account that $\cH$ is isometric, $\tilde\zeta_\cA(\cH)$ is a well defined isometric
matrix too.

The starting point to prove the theorem is recurrence (\ref{RR}) written as (\ref{RR1}).
The second relation in (\ref{RR1}) yields
\beq \label{PHI*}
\varpi_n\phi_n^* = \overline a_n\varpi_n\phi_n +
\sum_{k=0}^{n-1}\rho_n^-\rho_{n-1}^-\cdots\rho_{k+1}^-\overline a_k\varpi_k\phi_k,
\qquad n\geq1,
\eeq
where we set $a_0=1$. This identity, together with the first relation in (\ref{RR1}),
gives
\beq \label{HESS-ORF}
\ba{l}
\ds \varpi_n^*\phi_n = \sum_{k=0}^\infty\hat h_{k,n}\varpi_k\phi_k,
\medskip \cr
\hat h_{k,n} = \cases{
-a_{n+1}\rho_n^-\rho_{n-1}^-\cdots\rho_{k+1}^-\overline a_k & if $k<n$, \cr
-a_{n+1}\overline a_n & if $k=n$, \cr
\rho_{n+1}^+ & if $k=n+1$, \cr
0 & if $k>n+1$. \cr
}
\ea
\eeq
If we define the matrix $\hat\cH=(\hat h_{i,j})$, equality (\ref{HESS-ORF}) can be written as
\beq \label{HESS-ORF1}
\pmatrix{\phi_0(z) & \phi_1(z) & \cdots}
\left( \varpi^*_\cA(z) - \varpi_\cA(z) \,\hat\cH \right) = 0.
\eeq
Using (\ref{RHO}) we find that the Hessenberg matrix
\beq \label{hatH}
\hat\cH =
\pmatrix{
-a_1 & -\rho_1^- a_2 & - \rho_1^- \rho_2^- a_3 &
- \rho_1^- \rho_2^- \rho_3^- a_4 & \cdots
\cr
\rho_1^+ & -\overline a_1 a_2 & - \overline a_1 \rho_2^- a_3 &
- \overline a_1 \rho_2^- \rho_3^- a_4 & \cdots
\cr
0 & \rho_2^+ & -\overline a_2 a_3 & - \overline a_2 \rho_3^- a_4 & \cdots
\cr
0 & 0 & \rho_3^+ & -\overline a_3 a_4 & \cdots
\cr
\cdots & \cdots & \cdots & \cdots & \cdots}
\eeq
can be related to the isometric Hessenberg matrix $\cH$ given in (\ref{HESS}) by
\beq \label{hat-H}
\hat\cH = \eta_\cA^{-1} \cH \,\eta_\cA = \cH_\cA,
\eeq
where we have used that $\eta_{\cA^\dag}=\eta_\cA$ since $\cA$ is diagonal, so normal.
From this relation, (\ref{MOB2}) and (\ref{MOB1}) we see that (\ref{HESS-ORF1}) is
equivalent to
\beq \label{HESS-ORF2}
\pmatrix{\phi_0(z) & \phi_1(z) & \cdots} \left(z-\tilde\zeta_\cA(\cH) \right) = 0.
\eeq
This equality implies that $\cL$ is invariant under $T_\mu$, so the restriction $T_\mu
\upharpoonright \cL$ is well defined and $\tilde\zeta_\cA(\cH)$ is its matrix
representation with respect to $(\phi_n)_{n\geq0}$.

$T_\mu \upharpoonright \cL$ is an isometry because it is the restriction of a unitary
operator, which agrees with the fact that $\tilde\zeta_\cA(\cH)$ is isometric. Also,
$\tilde\zeta_\cA(\cH)$ and $\cH$ are unitary at the same time, that is, when
$\bsa\notin\ell^2$. Besides, $T_\mu \upharpoonright \cL$ is unitary iff $T_\mu\cL=\cL$.
This implies that $T_\mu^n\cL=\cL$ for any $n\in\Z$, so $\{z^n\}_{n\in\Z}\subset\cL$.
Hence $\cL=L^2_\mu$ because $\spn\{z^n\}_{n\in\Z}$ is dense in $L^2_\mu$. Conversely, if
$\cL=L^2_\mu$, then $T_\mu \upharpoonright \cL = T_\mu$ is unitary. Therefore,
$\tilde\zeta_\cA(\cH)$ is unitary iff the ORF $(\phi_n)_{n\geq0}$ are a basis of
$L^2_\mu$, i.e., iff $\tilde\zeta_\cA(\cH)$ represents the full operator $T_\mu$.
Finally, it is known that the condition $\sum_{n=1}^\infty(1-|\alpha_n|)=\infty$, which
is satisfied for $\bsalpha$ compactly included in $\D$, ensures that $\cL=\cP$ (see
\cite[Theorem 7.2.2]{BGHN99}) and so it implies the equivalence between $\cL=L^2_\mu$,
$\cP=L^2_\mu$ and $\log\mu' \notin L^1_m$ (see \cite[Corollary 7.2.4]{BGHN99}).

\epr

Given a measure $\mu$ on $\T$, the parameters $\bsa=\cS_\bsalpha(\mu)$ corresponding to
the ORF $(\phi_n)_{n\geq0}$ associated with $\bsalpha$ are in general different from the
parameters $\bsa^{(0)}=\cS_0(\mu)$ related to the OP $(\varphi_n)_{n\geq0}$. For
instance, if $\alpha_n=\alpha$ for all $n$, the comments at the beginning of this section
show that $\cS_\bsalpha(\mu_\alpha)=\cS_0(\mu)$. Taking into account that $\cS_\bsalpha$
is a bijection for a constant sequence $\bsalpha$, we conclude that
$\cS_\bsalpha(\mu)\neq\cS_0(\mu)$. Therefore, the equivalence $\cP=L^2_\mu
\Leftrightarrow \bsa\notin\ell^2$ is not trivial in the general case since the known
result in the polynomial situation is $\cP=L^2_\mu \Leftrightarrow
\bsa^{(0)}\notin\ell^2$.

Contrary to the polynomial case, the unitary matrix $\cV$ of the multiplication operator
with respect the ORF basis is not a Hessenberg matrix in general, but $\zeta_\cA(\cV)$ is
a Hessenberg matrix, where $\zeta_\cA$ is an operator M\"{o}bius transformation constructed
using all the poles of the related ORF. In the polynomial case $\cA=0$, thus
$\zeta_\cA(\cV)=\cV$ and $\cV$ becomes a Hessenberg matrix.

As a consequence of Theorem \ref{HESSORF} and the spectral properties of the unitary
multiplication operator, we have the following spectral interpretation of the support
of the orthogonality measure for ORF.

\bt \label{HESSORF1}

Let $\bsalpha$ be a sequence compactly included in $\D$, $\mu$ a measure on $\T$ such
that $\log\mu' \notin L^1_m$ and $(\phi_n)_{n\geq0}$ the corresponding ORF. If
$$
\cV=\tilde\zeta_{\cA}(\cH), \qquad \cA=\cA(\bsalpha),
\qquad \cH=\cH(\bsa), \qquad \bsa=\cS_\bsalpha(\mu),
$$
and $\cE$ is the spectral measure of $\cV$, then $\mu=\cE_{1,1}$. Besides,
$\supp\mu=\sigma(\cV)$ and the mass points of $\mu$ are the eigenvalues of $\cV$, which
have geometric multiplicity 1. $\lambda$ is a mass point iff
$(\phi_n(\lambda))_{n\geq0}\in\ell^2$. Given a mass point $\lambda$, the corresponding
eigenvectors of $\cV$ are spanned by $\pmatrix{\phi_0(\lambda) & \phi_1(\lambda) &
\cdots}^\dag$ and
$\mu(\{\lambda\})=\left(\sum_{n=0}^\infty|\phi_n(\lambda)|^2\right)^{-1}$.

\et

\bpr

Under the hypothesis of the theorem, $\cV$ is the matrix representation of the full
operator $T_\mu$ with respect to $(\phi_n)_{n\geq0}$. Hence, if $E$ is the spectral
measure of $T_\mu$, then $\mu(\cdot) = \langle \phi_0,E(\cdot)\phi_0 \rangle_\mu =
\cE_{1,1}(\cdot)$. Also, $\supp\mu=\sigma(T_\mu)=\sigma(\cV)$ and the mass points of
$\mu$ are the eigenvalues of $T_\mu$, that is, the eigenvalues of $\cV$, which have
therefore geometric multiplicity 1. If $\lambda$ is a mass point, we know that the
characteristic function $\cX_{\{\lambda\}}$ of $\{\lambda\}$ is a related eigenvector of
$T_\mu$, so,
$\langle\phi_n,\cX_{\{\lambda\}}\rangle_\mu=\mu(\{\lambda\})\,\overline{\phi_n(\lambda)}$
is the $(n+1)$-th component of a corresponding eigenvector of $\cV$. This implies that
$(\phi_n(\lambda))_{n\geq0}\in\ell^2$. Conversely, if $\lambda$ is an arbitrary complex
number such that $(\phi_n(\lambda))_{n\geq0}\in\ell^2$, relation (\ref{HESS-ORF2}) shows
that $\pmatrix{\phi_0(\lambda) & \phi_1(\lambda) & \cdots}$ is a left eigenvector of
$\cV$ with eigenvalue $\lambda$. Due to the unitarity of $\cV$, $\lambda\in\T$ and the
above statement is equivalent to saying that $\pmatrix{\phi_0(\lambda) & \phi_1(\lambda)
& \cdots}^\dag$ is a (right) eigenvector of $\cV$ with eigenvalue $\lambda$. Therefore,
$\lambda$ is a mass point of $\mu$. Also, the identity
$$
\mu(\{\lambda\}) =
\langle\cX_{\{\lambda\}},\cX_{\{\lambda\}}\rangle_\mu =
\sum_{n=0}^\infty
\langle\cX_{\{\lambda\}},\phi_n\rangle_\mu \langle\phi_n,\cX_{\{\lambda\}}\rangle_\mu =
\sum_{n=0}^\infty \mu(\{\lambda\})^2 |\phi_n(\lambda)|^2
$$
proves that $\mu(\{\lambda\})=\left(\sum_{n=0}^\infty|\phi_n(\lambda)|^2\right)^{-1}$.

\epr

The fact that the representation $\cV$ is not a Hessenberg matrix, but a M\"{o}bius
transformation of a Hessenberg matrix, makes the rational case more complicated than the
polynomial one. However, the Hessenberg structure can be kept if we formulate the
spectral results in terms of pairs of operators.

Remember that, given a Hilbert space $H$ and two operators $T,S\in\B_H$, the spectrum
and point spectrum of the pair $(T,S)$ are respectively the sets
$$
\ba{l} \sigma(T,S) =
\{ \lambda\in\overline\C : T - \lambda S \hbox{ has no inverse in } \B_H \},
\medskip \cr
\sigma_p(T,S) =
\{ \lambda\in\overline\C : T - \lambda S \hbox{ is not injective} \}.
\ea
$$
In the finite-dimensional case both sets coincide. The elements of $\sigma_p(T,S)$ are
called eigenvalues of the pair, and the eigenvectors of $(T,S)$ corresponding to an
eigenvalue $\lambda$ are the elements $x \in H\setminus\{0\}$ such that $(T -\lambda S)x
= 0$. In these definitions it is assumed that, if $\lambda=\infty$, $T - \lambda S$ must
be substituted by $S$.

With the above terminology, the isometric matrix $\cV$ and the Hessenberg pair
$(\tilde\varpi_\cA^*(\cH_\cA),\tilde\varpi_\cA(\cH_\cA))$ have the same spectrum and
eigenvalues because $\tilde\varpi_\cA(\cH_\cA)^{\pm1}\in\B_{\ell^2}$ when
$\bsalpha$ is compactly included in $\D$. So, Theorem \ref{HESSORF1} can be
obviously rewritten substituting $\cV$ by the pair
$(\tilde\varpi_\cA^*(\cH_\cA),\tilde\varpi_\cA(\cH_\cA))$. Notice that, given an
eigenvalue $\lambda$, $\pmatrix{\phi_0(\lambda) & \phi_1(\lambda) & \cdots}$ is a left
eigenvector of the pair, i.e.,
$$
\pmatrix{\phi_0(\lambda) & \phi_1(\lambda) & \cdots}
(\tilde\varpi_\cA^*(\cH_\cA) - \lambda \,\tilde\varpi_\cA(\cH_\cA)) = 0.
$$
Moreover, $\cH_\cA=\eta_\cA^{-1}\cH\,\eta_\cA$ with $\eta_\cA^{\pm1}\in\B_{\ell^2}$ due
to the restrictions on $\bsalpha$. Therefore, Theorem \ref{HESSORF1} also
holds substituting $\cV$ by the Hessenberg pair
$(\tilde\varpi_\cA^*(\cH),\tilde\varpi_\cA(\cH))$, but the left eigenvectors with
eigenvalue $\lambda$ are spanned by $\pmatrix{\phi_0(\lambda) & \phi_1(\lambda) & \cdots}
\eta_\cA^{-1}$.

\subsection{Zeros of ORF and Hessenberg matrices} \label{Z-ORF-HM}

Let $\mu$ be a measure on $\T$ and $\bsa=\cS_0(\mu)$. As we pointed out in Section
\ref{OP-ORF}, the characteristic polynomial of the orthogonal truncation of $T_\mu$ on
$\cP_n=\spn\{1,z,\dots,z^{n-1}\}$ is a multiple of the $n$-th OP related to $\mu$. From
this result, the relation between the zeros of the $n$-th OP and the eigenvalues of the
principal submatrices of $\cH(\bsa)$ follows.

To obtain a similar result for the ORF associated with a sequence $\bsalpha$ we have to
consider the operator multiplication by $\zeta_n$ in $L^2_\mu$, i.e.,
$$
\zeta_n(T_\mu) \colon \mathop{L^2_\mu \to L^2_\mu} \limits_{f \; \longrightarrow \; \zeta_nf}
$$
and the orthogonal truncation of $\zeta_n(T_\mu)$ on $\cL_n$. This orthogonal truncation
is defined by $\zeta_n(T_\mu)^{(\cL_n)} = L_n \zeta_n(T_\mu) \upharpoonright \cL_n$, where
the operator $L_n \colon L^2_\mu \to L^2_\mu$ is the orthogonal projection on $\cL_n$.
The following theorem is the starting point to identify the zeros of the ORF as the
eigenvalues of some finite matrices related to $\cA(\bsalpha)$ and $\cH(\bsa)$.

We remind that the $n$-th ORF $\phi_n$ has the form
$$
\phi_n = \frac{p_n}{\pi_n},
\qquad p_n\in\cP_{n+1}\setminus\cP_n,
\qquad \pi_n=\varpi_1\cdots\varpi_n,
$$
and the zeros of $\phi_n$, which are the zeros of the polynomial $p_n$, lie on $\D$ (see
\cite[Corollary 3.2.2]{BGHN99}).

\bt \label{ORF-ZEROS}

Let $\bsalpha$ be an arbitrary sequence in $\D$, $\mu$ a measure on $\T$ and
$\phi_n=p_n/\pi_n$ the related $n$-th ORF. Then:
\be
\item If $Z_n$ is the set of zeros of $\phi_n$, $\zeta_n(Z_n)$ is the set of eigenvalues
      of $\zeta_n(T_\mu)^{(\cL_n)}$ and these eigenvalues have geometric multiplicity 1.
\item If $p_n(z)\propto\prod_{k=1}^n(z-\lambda_k)$, the characteristic polynomial of
      $\zeta_n(T_\mu)^{(\cL_n)}$ is
      $$
      \prod_{k=1}^n(z-\zeta_n(\lambda_k)).
      $$
\ee

\et

\bpr

$f\in\cL_n\setminus\{0\}$ is an eigenvector of $\zeta_n(T_\mu)^{(\cL_n)}$ with eigenvalue
$w$ iff $(L_n\zeta_n-w)f=0$, that is, $L_n(\zeta_n-w)f=0$. This is equivalent to state
that $(\zeta_n-w)f\in\cL_n^{\bot\cL_{n+1}}=\spn\{\phi_n\}$, or, in other words,
$f\propto\phi_n(\zeta_n-w)^{-1}$. Writing $w=\zeta_n(\lambda)$ and using (\ref{MOB4}) we
find that this condition can be expressed as $f(z) \propto
p_n(z)(z-\lambda)^{-1}/\pi_{n-1}(z)$ with $\lambda \in Z_n$. This proves item 1.

Item 2 is equivalent to assert that the algebraic multiplicity $m_w$ of any eigenvalue
$w=\zeta_n(\lambda)$ of $\zeta_n(T_\mu)^{(\cL_n)}$ is equal to the multiplicity of
$\lambda$ as a root of $p_n$. Since the geometric multiplicity of $w$ is 1, $m_w \geq k$
iff there exists $f\in\cL_n$ such that $(L_n\zeta_n-w)^kf=0$ and
$(L_n\zeta_n-w)^{k-1}f\neq0$. Analogously to the previous discussion, we find that these
two conditions are equivalent to $f \in \spn\{\phi_n(\zeta_n-w)^{-j}\}_{j=1}^k \setminus
\spn\{\phi_n(\zeta_n-w)^{-j}\}_{j=1}^{k-1}$, i.e., to $f=p/\pi_{n-1}$ with $p(z) \in
\spn\{\varpi_n^{j-1}(z)p_n(z)(z-\lambda)^{-j}\}_{j=1}^k \setminus
\spn\{\varpi_n^{j-1}(z)p_n(z)(z-\lambda)^{-j}\}_{j=1}^{k-1}$, as can be seen using
(\ref{MOB4}) again. Hence, by induction on $k$ we find that $m_w \geq k$ implies that the
multiplicity of $\lambda$ as a root of $p_n$ is not less than $k$. Conversely, if the
multiplicity of $\lambda$ as a root of $p_n$ is greater than or equal to $k$,
$f(z)=\phi_n(z)(\zeta_n(z)-w)^{-k}\propto\varpi_n^{k-1}(z)p_n(z)(z-\lambda)^{-k}/\pi_{n-1}(z)\in\cL_n$
and the above results ensure that $(L_n\zeta_n-w)^kf=0$ and $(L_n\zeta_n-w)^{k-1}f\neq0$,
so $m_w \geq k$.

\epr

The next step is to obtain a matrix representation of the orthogonal truncation
$\zeta_n(T_\mu)^{(\cL_n)}$, so that we can give a matrix version of the above theorem. In
the following results the subscript ${}_n$ on a matrix means the principal submatrix of
order $n$ of such a matrix. This notation will be used throughout the rest of the paper.

\bt \label{HESSORFmatrix}

Let $\bsalpha$ be an arbitrary sequence in $\D$, $\mu$ a measure on $\T$ and
$(\phi_n)_{n\geq0}$ the related ORF. If $\cA=\cA(\bsalpha)$ and $\cH=\cH(\bsa)$ with
$\bsa=\cS_\bsalpha(\mu)$, the matrix of $\zeta_n(T_\mu)^{(\cL_n)}$ with respect to
$(\phi_k)_{k=0}^{n-1}$ is $\zeta_n(\cV^{(n)})$, where
$$
\cV^{(n)}=\tilde\zeta_{\cA_n}(\cH_n), \qquad \|\cV^{(n)}\|=1.
$$

\et

\bpr

$\|\cA_n\|<1$, thus $\tilde\zeta_{\cA_n}$ maps $\T_{\C^n}$ onto itself. Also, from
(\ref{FAC1}) we obtain the factorization
$$
\cH_n = \pmatrix{\Theta_1 & \cr & I_{n-2}} \kern-3.5pt \pmatrix{I_1 && \cr & \Theta_2 &
\cr && I_{n-3}} \kern-1pt\cdots\pmatrix{I_{n-2} & \cr & \Theta_{n-1}}
\kern-3pt\pmatrix{I_{n-1} & \cr & -a_n},
$$
all the factors being unitary except the last one which has norm 1, thus $\|\cH_n\|=1$.
Hence, $\cV^{(n)}=\tilde\zeta_{\cA_n}(\cH_n)$ is well defined and $\|\cV^{(n)}\|=1$. A
similar reason shows that $\zeta_n(\cV^{(n)})$ is a well defined matrix with norm 1.

To prove the theorem, let us write the first $n$ equations of (\ref{HESS-ORF1}) as
\beq \label{HESS-ORF-zeros}
\ba{l}
\pmatrix{\phi_0(z) & \kern-5pt \cdots \kern-5pt & \phi_{n-1}(z)}
\left( \varpi^*_{\cA_n}(z) - \varpi_{\cA_n}(z) \,\hat\cH_n \right) =
b_n \varpi_n(z) \phi_n(z),
\medskip \cr
b_n = \rho_n^+ \pmatrix{0 & 0 & \cdots & 0 & 1} \in \C^n.
\ea
\eeq
Then, identities (\ref{MOB2}), (\ref{MOB1}) and the equality
$\hat\cH_n=\eta_{\cA_n}^{-1}\cH_n\,\eta_{\cA_n}=(\cH_n)_{\cA_n}$
obtained from (\ref{hat-H}) transform (\ref{HESS-ORF-zeros}) into
\beq \label{HESS-ORF-zeros1}
\pmatrix{\phi_0(z) & \cdots & \phi_{n-1}(z)} \left( z - \cV^{(n)} \right) =
c_n \varpi_n(z) \phi_n(z),
\qquad c_n\in\C^n.
\eeq
Using (\ref{MOB3}) we get
$$
\pmatrix{\phi_0(z) & \cdots & \phi_{n-1}(z)} \left( \zeta_n(z) - \zeta_n(\cV^{(n)})
\right) = d_n \phi_n(z), \qquad d_n\in\C^n.
$$
Hence, taking into account that
$$
L_n\phi_k=\cases{\phi_k & if $k < n$, \cr 0 & if $k \geq n$,}
$$
we finally obtain
$$
\pmatrix{L_n\zeta_n\phi_0 & \cdots & L_n\zeta_n\phi_{n-1}} =
\pmatrix{\phi_0 & \cdots & \phi_{n-1}} \zeta_n(\cV^{(n)}),
$$
which proves that $\zeta_n(\cV^{(n)})$ is the matrix of $L_n \zeta_n(T_\mu)
\upharpoonright \cL_n$ with respect to $(\phi_k)_{k=0}^{n-1}$.

\epr

Theorems \ref{ORF-ZEROS} and \ref{HESSORFmatrix} provide a spectral interpretation of the
zeros of ORF in terms of operator M\"{o}bius transformations of Hessenberg matrices.

\bt \label{HESSORFzeros}

Let $\bsalpha$ be an arbitrary sequence in $\D$, $\mu$ a measure on $\T$,
$(\phi_n)_{n\geq0}$ the corresponding ORF, $\cA=\cA(\bsalpha)$ and $\cH=\cH(\bsa)$ with
$\bsa=\cS_\bsalpha(\mu)$. If $\cV^{(n)}=\tilde\zeta_{\cA_n}(\cH_n)$, then:
\be
\item The zeros of $\phi_n$ are the eigenvalues of $\cV^{(n)}$, which have geometric
      multiplicity 1. If $\lambda$ is a zero of $\phi_n$, the related left eigenvectors of
      $\cV^{(n)}$ are spanned by $\pmatrix{\phi_0(\lambda)&\cdots&\phi_{n-1}(\lambda)}$.
\item $\ds \phi_n=\frac{p_n}{\pi_n}$ with $p_n$ proportional to the characteristic
      polynomial of $\cV^{(n)}$.
\ee

\et

\bpr

From Theorems \ref{ORF-ZEROS} and \ref{HESSORFmatrix}, the eigenvalues of
$\zeta_n(\cV^{(n)})$ have geometric multiplicity 1 and
$\sigma(\zeta_n(\cV^{(n)}))=\zeta_n(Z_n)$, where $Z_n$ are the zeros of $\phi_n$. Also,
the characteristic polynomial of $\zeta_n(\cV^{(n)})$ is
$\prod_{k=1}^n(z-\zeta_n(\lambda_k))$, where $p_n(z)=\prod_{k=1}^n(z-\lambda_k)$.
$\sigma(\zeta_n(\cV^{(n)}))=\zeta_n(\sigma(\cV^{(n)}))$, so, bearing in mind that
$\zeta_n$ is bijective, $\sigma(\cV^{(n)})=Z_n$. Furthermore, given an eigenvalue
$\lambda$ of $\cV^{(n)}$, the corresponding eigenvalue $\zeta_n(\lambda)$ of
$\zeta_n(\cV^{(n)})$ has the same geometric and algebraic multiplicity. Therefore,
$\prod_{k=1}^n(z-\lambda_k)$ is the characteristic polynomial of $\cV^{(n)}$. Finally, if
$\lambda$ is a zero of $\phi_n$, (\ref{HESS-ORF-zeros1}) shows that
$\pmatrix{\phi_0(\lambda)&\cdots&\phi_{n-1}(\lambda)}$ is a left eigenvector of
$\cV^{(n)}$ with eigenvalue $\lambda$.

\epr

From (\ref{MOB2}) and (\ref{MOB1}) we know that
$$
\varpi^*_{\cA_n}(z) - \varpi_{\cA_n}(z) \,\hat\cH_n =
z \,\tilde\varpi_{\cA_n}(\hat\cH_n)-\tilde\varpi_{\cA_n}^*(\hat\cH_n) =
(z-\cV^{(n)}) \,\tilde\varpi_{\cA_n}(\hat\cH_n).
$$
This gives other alternatives to express $p_n$ as a determinant, like
$$
p_n(z) \propto \det\left(\varpi^*_{\cA_n}(z) - \varpi_{\cA_n}(z) \,\hat\cH_n\right) =
\det\left(z \,\tilde\varpi_{\cA_n}(\hat\cH_n)-\tilde\varpi_{\cA_n}^*(\hat\cH_n)\right).
$$
The interest of the above expressions is that they show that $p_n$ can be calculated as a
determinant of a Hessenberg matrix. Furthermore, the last expression provides a new
spectral interpretation of the zeros of the ORF, related to the concept of the spectrum
of a pair of operators. It shows that the zeros of $\phi_n$ are the eigenvalues of the
Hessenberg pair $(\tilde\varpi_{\cA_n}^*(\hat\cH_n),\tilde\varpi_{\cA_n}(\hat\cH_n))$.
Also, according to Theorem \ref{HESSORFzeros}, the left eigenvectors of
$(\tilde\varpi_{\cA_n}^*(\hat\cH_n),\tilde\varpi_{\cA_n}(\hat\cH_n))$ corresponding to an
eigenvalue $\lambda$ are spanned by
$\pmatrix{\phi_0(\lambda)&\cdots&\phi_{n-1}(\lambda)}$. Taking into account that
$\hat\cH_n= \eta_{\cA_n}^{-1}\cH_n\,\eta_{\cA_n}$, the zeros of $\phi_n$ can be also
understood as the eigenvalues of the Hessenberg pair
$(\tilde\varpi_{\cA_n}^*(\cH_n),\tilde\varpi_{\cA_n}(\cH_n))$, the left eigenvectors with
eigenvalue $\lambda$ being spanned by
$\pmatrix{\phi_0(\lambda)&\cdots&\phi_{n-1}(\lambda)} \eta_{\cA_n}^{-1}$. Indeed,
$$
p_n(z) \propto \det\left(\varpi^*_{\cA_n}(z) - \varpi_{\cA_n}(z) \,\cH_n\right) =
\det\left(z\,\tilde\varpi_{\cA_n}(\cH_n)-\tilde\varpi_{\cA_n}^*(\cH_n)\right).
$$

\medskip

Apart from the sequence $(\phi_n)_{n\geq0}$ of ORF, another remarkable rational functions
arise in the theory of ORF. They are the so called para-orthogonal rational functions
(PORF), given by \beq \label{PORFv} Q_n^v = \phi_n + v \phi_n^*, \qquad v\in\T. \eeq The
PORF are the generalization of the POP to the rational case. Analogously to the
POP, the interest of the PORF $Q_n^v$ relies on the fact that, contrary to the ORF
$\phi_n$, it has $n$ different zeros lying on $\T$ which, thus, play an important role in
quadrature formulas and rational moment problems (see \cite[Chapters 5 and 10]{BGHN99}).
These quadrature formulas associate with each PORF $Q_n^v$ a measure $\mu_n^v$ supported
on its zeros with a mass $(\sum_{k=0}^{n-1}|\phi_k(\lambda)|^2)^{-1}$ at each zero
$\lambda$. Such quadrature formulas are exact in $\cL_{n-1}\cL_{n-1*}$.

Due to the exactness of the quadrature formulas, the first $n$ OP $(\phi_k)_{k=0}^{n-1}$
related to $\mu$ are also an orthonormal basis of the $n$-dimensional Hilbert space
$L^2_{\mu_n^v}$. The convergence properties of the quadrature formulas imply that, for
any sequence $(v_n)_{n\geq1}$ in $\T$, the sequence $(\mu_n^{v_n})_{n\geq1}$ of measures
$*$-weakly converges to the orthogonality measure $\mu$ on $\T$ whenever
$\sum_{n=1}^\infty(1-|\alpha_n|)=\infty$.

A spectral interpretation can be also obtained for the zeros of the PORF. To understand
this, let us write the PORF in an equivalent way. Using recurrence (\ref{RR}) we find that
\beq \label{PORFu}
Q_n^v = (1+\overline a_n v) \,e_n {\varpi_{n-1}\over\varpi_n}
\,(\zeta_{n-1} \phi_{n-1} + u \phi_{n-1}^*),
\qquad u=\tilde\zeta_{a_n}(v).
\eeq
Notice that the parameter $u$ goes through the full unit circle as the parameter $v$ does
so. (\ref{PORFu}) shows that, like $\phi_n$, $Q_n^v$ is obtained from $n$ steps of
recurrence (\ref{RR}), but changing in the $n$-th step $a_n\in\D$ by
$u=\tilde\zeta_{a_n}(v)\in\T$. The analogous substitution in $\cH_n$ gives
$$
\cH_n^u =
\pmatrix{\Theta_1 & \cr & I_{n-2}}
\kern-3.5pt \pmatrix{I_1 && \cr & \Theta_2 & \cr && I_{n-3}}
\kern-1pt\cdots\pmatrix{I_{n-2} & \cr & \Theta_{n-1}}
\kern-3pt\pmatrix{I_{n-1} & \cr & -u},
$$
which is obviously a unitary matrix. Therefore, $\tilde\zeta_{\cA_n}(\cH_n^u)$ is unitary too
because $\tilde\zeta_{\cA_n}$ preserves the unitarity.

The following result provides a spectral interpretation of the zeros of the PORF $Q_n^v$
in terms of the unitary Hessenberg matrix $\cH_n^u$, as well as a connection of such a
matrix with the unitary multiplication operator $T_{\mu_n^v}$. It can be understood as a
limit case of Theorems \ref{HESSORFmatrix} and \ref{HESSORFzeros}.

\bt \label{HESSPORFzeros}

Let $\bsalpha$ be an arbitrary sequence in $\D$, $\mu$ a measure on $\T$,
$(\phi_n)_{n\geq0}$ the corresponding ORF, $\cA=\cA(\bsalpha)$ and $\cH=\cH(\bsa)$ with
$\bsa=\cS_\bsalpha(\mu)$. If $Q_n^v=\phi_n+v\phi_n^*$ is the $n$-th PORF related to
$v\in\T$ and $\mu_n^v$ is the associated measure, then:
\be
\item The matrix of $T_{\mu_n^v}$ with respect to $(\phi_k)_{k=0}^{n-1}$ is
      $$
      \cV^{(n;u)}=\tilde\zeta_{\cA_n}(\cH_n^u), \qquad u=\tilde\zeta_{a_n}(v).
      $$
\item The zeros of $Q_n^v$ are the eigenvalues of $\cV^{(n;u)}$.
      If $\lambda$ is a zero of $Q_n^v$, the related eigenvectors of $\cV^{(n;u)}$ are
      spanned by $\pmatrix{\phi_0(\lambda)&\cdots&\phi_{n-1}(\lambda)}^\dag$.
\item $\ds Q_n^v=\frac{q_n^v}{\pi_n}$ with $q_n^v$ proportional to the characteristic
      polynomial of $\cV^{(n;u)}$.
\ee

\et

\bpr

Using (\ref{PHI*}) in (\ref{PORFu}) we find that
$$
\ba{l}
\ds \varpi_{n-1}^*\phi_{n-1} =
\sum_{k=0}^{n-1}\hat h_{k,n-1}^u\varpi_k\phi_k +
\frac{\rho_n^+}{1+\overline a_nv} \, \varpi_n Q_n^v,
\medskip \cr
\hat h_{k,n-1}^u = \cases{
-u\rho_{n-1}^-\rho_{n-2}^-\cdots\rho_{k+1}^-\overline a_k & if $k<n-1$, \cr
-u\overline a_{n-1} & if $k=n-1$,
}
\qquad u=\tilde\zeta_{a_n}(v).
\ea
$$
This relation together with the first $n-1$ equations of (\ref{HESS-ORF1}) lead to the matrix
identity
$$
\ba{l}
\pmatrix{\phi_0(z) & \kern-5pt \cdots \kern-5pt & \phi_{n-1}(z)}
\left( \varpi^*_{\cA_n}(z) - \varpi_{\cA_n}(z) \,\hat\cH_n^u \right) =
b_n \varpi_n(z) Q_n^v(z),
\cr
\ds b_n = \frac{\rho_n^+}{1+\overline a_n v} \pmatrix{0 & 0 & \cdots & 0 & 1} \in \C^n.
\ea
$$
where $\hat\cH_n^u=\eta_{\cA_n}^{-1}\cH_n^u\,\eta_{\cA_n}=(\cH_n^u)_{\cA_n}$.
So, (\ref{MOB2}) and (\ref{MOB1}) give
\beq \label{HESS-PORF-zeros1}
\pmatrix{\phi_0(z) & \cdots & \phi_{n-1}(z)} \left( z - \cV^{(n;u)} \right) =
c_n \varpi_n(z) Q_n^v(z),
\qquad c_n\in\C^n.
\eeq
$Q_n^v=0$ in $L^2_{\mu_n^v}$, thus (\ref{HESS-PORF-zeros1}) implies that $\cV^{(n;u)}$ is
the matrix of $T_{\mu_n^v}$ with respect to $(\phi_k)_{k=0}^{n-1}$. The rest of the
statements are a consequence of this one and the properties of the multiplication
operators, similarly to the proof of Theorem \ref{HESSORF1}. Alternatively, they can be
obtained directly from relation (\ref{HESS-PORF-zeros1}), the unitarity of $\cV^{(n:u)}$
and the fact that $Q_n^v$ has $n$ different zeros.

\epr

Analogously to the comments after Theorem \ref{HESSORFzeros}, if $u=\tilde\zeta_{a_n}(v)$,
$$
q_n^v(z) \propto \det\left(\varpi^*_{\cA_n}(z) - \varpi_{\cA_n}(z) \,\cH_n^u\right) =
\det\left(z \,\tilde\varpi_{\cA_n}(\cH_n^u)-\tilde\varpi_{\cA_n}^*(\cH_n^u)\right),
$$
which gives $q_n^v$ as a determinant of a Hessenberg matrix too. The zeros of $Q_n^v$ are
the eigenvalues of the Hessenberg pair
$(\tilde\varpi_{\cA_n}^*(\cH_n^u),\tilde\varpi_{\cA_n}(\cH_n^u))$,
whose left eigenvectors with eigenvalue $\lambda$ are spanned by
$\pmatrix{\phi_0(\lambda)&\cdots&\phi_{n-1}(\lambda)} \eta_{\cA_n}^{-1}$.

\section{ORF and five-diagonal matrices} \label{ORF-5M}

Apart from the presence of operator M\"{o}bius transformations, there are some drawbacks in
the spectral theory of ORF previously developed: the appearance of a Hessenberg matrix
$\cH$ instead of a band one, the complicated dependence of $\cH=\cH(\bsa)$ on the
parameters $\bsa$, and the fact that it represents the full multiplication operator
$T_\mu$ only for certain measures $\mu$ on $\T$. We will not be able to avoid the
operator M\"{o}bius transformations because they are linked to the ORF, but the other
problems can be overcome by choosing a different basis of ORF in $L^2_\mu$.

The key idea is to use, instead of the ORF $(\phi_n)_{n\geq0}$ with poles in $\E$, other
ones whose poles are alternatively in $\E$ and $\D$. For this pourpose we define the
finite odd and even Blaschke products
$$
\ba{l}
B_0^o=B_0^e=1,
\smallskip \cr
B_n^o=\zeta_1\zeta_3\cdots\zeta_{2n-1},
\qquad
B_n^e=\zeta_2\zeta_4\cdots\zeta_{2n},
\qquad n\geq1.
\ea
$$
Consider the rational functions $(\chi_n)_{n\geq0}$ given by
\beq
\label{LORF}
\chi_{2n} = B_{n*}^e\phi_{2n}^*,
\qquad
\chi_{2n+1} = B_{n*}^e\phi_{2n+1},
\qquad n\geq0.
\eeq
Since $\zeta_{n*}=1/\zeta_n$, the subspaces $\cM_n=\spn\{\chi_0,\dots,\chi_{n-1}\}$ are
$$
\ba{l}
\cM_{2n}=B_{n-1*}^e\cL_{2n}=\spn\{B_{0*}^e,B_1^o,B_{1*}^e,\dots,B_{n-1*}^e,B_n^o\},
\medskip \cr
\cM_{2n+1}=B_{n*}^e\cL_{2n+1}=\spn\{B_{0*}^e,B_1^o,B_{1*}^e,\dots,B_n^o,B_{n*}^e\},
\ea
$$
i.e., $\cM_{2n}$ and $\cM_{2n+1}$ are the sets of rational functions whose poles, counted
with multiplicity, lie on $(\hat\alpha_1,\alpha_2,\dots,\alpha_{2n-2},\hat\alpha_{2n-1})$
and $(\hat\alpha_1,\alpha_2,\dots,\hat\alpha_{2n-1},\alpha_{2n})$ respectively. We will
use the notation $\cM_\infty=\spn\{\chi_n\}_{n\geq0}=\cup_{n\geq1}\cM_n$ and $\cM$ for
the closure of $\cM_\infty$ in $L^2_\mu$.

The orthonormality conditions $\phi_n \bot \cL_n$ and $\langle\phi_n,\phi_n\rangle_\mu=1$
can be rewritten using $\phi_n^*$ as $\phi_n^* \bot \zeta_n\cL_n$ and
$\langle\phi_n^*,\phi_n^*\rangle_\mu=1$. Hence, the orthonormality of $(\phi_n)_{n\geq0}$ is
equivalent to $\chi_{2n} \bot B_{n*}^e\zeta_{2n}\cL_{2n} = \cM_{2n}$, $\chi_{2n+1}
\bot B_{n*}^e\cL_{2n+1} = \cM_{2n+1}$ and $\langle\chi_n,\chi_n\rangle_\mu=1$, i.e., to the
orthonormality of $(\chi_n)_{n\geq0}$. The sequence $(\chi_n)_{n\geq0}$ is therefore the
result of orthonormalizing $(B_{0*}^e,B_1^o,B_{1*}^e,B_2^o,B_{2*}^e,\dots)$ in $L^2_\mu$.
Hence, relation (\ref{LORF}) establishes a bijection between ORF associated with the
sequences $(\alpha_n)_{n\geq1}$ and
$(\alpha_1,\hat\alpha_2,\alpha_3,\hat\alpha_4,\dots)$. We can consider also the ORF
associated with the sequence $(\hat\alpha_1,\alpha_2,\hat\alpha_3,\alpha_4,\dots)$, i.e.,
the ORF that arise from the orthonormalization of
$(B_0^e,B_{1*}^o,B_1^e,B_{2*}^o,B_2^e,\dots)$ in $L^2_\mu$. This ORF are
$(\chi_{n*})_{n\geq0}$, which are related to $(\phi_n)_{n\geq0}$ by
\beq
\label{LORF*}
\chi_{2n*} = B_{n*}^o\phi_{2n}, \qquad \chi_{2n+1*} = B_{n+1*}^o\phi_{2n+1}^*,
\qquad n\geq0.
\eeq

As a conclusion, the two possibilities to generate ORF which alternate poles in $\E$ and
$\D$ are related between them and, also, to the ORF with poles in $\E$. The last ORF have
been extensively studied, thus every known result for them can be easily translated to
the first ones. This is an interesting result because the ORF with poles arbitrarily
located in $\overline\C\backslash\T=\E\cup\D$ are not so well known than the ORF with
poles in $\E$. However, surprisingly, we will use the above connection to obtain new
results for ORF with poles in $\E$ using ORF with alternating poles in $\E$ and $\D$. The
basic idea is that the ORF $(\chi_n)_{n\geq0}$ provide new matrix tools for the analysis
of questions concerning the ORF $(\phi_n)_{n\geq0}$. The reason for this is the different
nature of the recurrence satisfied by $(\chi_n)_{n\geq0}$, which, as we will see, is a
five-term linear recurrence relation. We could think that recurrence (\ref{RR}) for
$(\phi_n)_{n\geq0}$ should be better because it is two-term, but the presence of
$\phi_n^*$ causes a non-linearity which is the origin of the difficulties to connect with
linear operator theory. Alternatively, expanding $\phi_n^*$ in the basis
$(\phi_n)_{n\geq0}$ we find the linear relation (\ref{HESS-ORF}), but this is a
recurrence with a number of terms that increases with $n$, giving rise to a
representation of $T_\mu$ related to a Hessenberg instead of a band matrix. On the
contrary, the five-term linear recurrence for the ORF $(\chi_n)_{n\geq0}$ provides a
matrix representation of $T_\mu$ in terms of five-diagonal matrices, as the following
theorem states.

\bt \label{BANDORF}

Let $\bsalpha$ be a sequence compactly included in $\D$, $\mu$ a measure on $\T$ and
$\bsa=\cS_\bsalpha(\mu)$. Then, the ORF $(\chi_n)_{n\geq0}$ associated with
$(\alpha_1,\hat\alpha_2,\alpha_3,\hat\alpha_4,\dots)$ are a basis of $L^2_\mu$ and the
matrix of $T_\mu$ with respect to $(\chi_n)_{n\geq0}$ is $\cU=\tilde\zeta_\cA(\cC)$,
where $\cA=\cA(\bsalpha)$ and $\cC=\cC(\bsa)$ is given in (\ref{CMV}).

\et

\bpr

Since $\bsalpha$ is compactly included in $\D$, $\|\cA\|<1$ and, hence, $\tilde\zeta_\cA$
maps unitarity matrices into unitary matrices. Thus $\tilde\zeta_\cA(\cC)$ is a well
defined unitary matrix because $\cC$ is unitary.

Using (\ref{RR1}) and (\ref{LORF}) we find that, for $n\geq1$,
\beq \label{RRchi}
\ba{l}
\varpi_{2n-1}^*\chi_{2n-1} = B_{n-1*}^e
(\rho_{2n}^+\varpi_{2n}\phi_{2n} - a_{2n}\varpi_{2n-1}\phi_{2n-1}^*) =
\smallskip \cr \kern62pt
= B_{n*}^e \rho_{2n}^+
(\rho_{2n+1}^+\varpi_{2n+1}\phi_{2n+1} - a_{2n+1}\varpi_{2n}\phi_{2n}^*) \, -
\smallskip \cr \kern75pt
- \, B_{n-1*}^e a_{2n}
(\overline a_{2n-1}\varpi_{2n-1}\phi_{2n-1} + \rho_{2n-1}^-\varpi_{2n-2}\phi_{2n-2}^*) =
\smallskip \cr \kern62pt
= \rho_{2n}^+\rho_{2n+1}^+ \varpi_{2n+1}\chi_{2n+1} -
\rho_{2n}^+a_{2n+1} \varpi_{2n}\chi_{2n} \, -
\smallskip \cr \kern75pt
- \, \overline a_{2n-1}a_{2n} \varpi_{2n-1}\chi_{2n-1} -
\rho_{2n-1}^-a_{2n} \varpi_{2n-2}\chi_{2n-2},
\medskip \cr
\varpi_{2n}^*\chi_{2n} = B_{n-1*}^e
(\overline a_{2n}\varpi_{2n}\phi_{2n} + \rho_{2n}^-\varpi_{2n-1}\phi_{2n-1}^*) =
\smallskip \cr \kern40pt
= B_{n*}^e \overline a_{2n}
(\rho_{2n+1}^+\varpi_{2n+1}\phi_{2n+1} - a_{2n+1}\varpi_{2n}\phi_{2n}^*) \, +
\smallskip \cr \kern53pt
+ \, B_{n-1*}^e \rho_{2n}^-
(\overline a_{2n-1}\varpi_{2n-1}\phi_{2n-1} + \rho_{2n-1}^-\varpi_{2n-2}\phi_{2n-2}^*) =
\smallskip \cr \kern40pt
= \overline a_{2n}\rho_{2n+1}^+ \varpi_{2n+1}\chi_{2n+1} -
\overline a_{2n}a_{2n+1} \varpi_{2n}\chi_{2n} \, +
\smallskip \cr \kern53pt
+ \, \overline a_{2n-1}\rho_{2n}^- \varpi_{2n-1}\chi_{2n-1} +
\rho_{2n-1}^-\rho_{2n}^-\varpi_{2n-2}\chi_{2n-2},
\ea
\eeq
while
$$
\varpi_0^*\chi_0 = \rho_1^+\varpi_1\chi_1 - a_1\varpi_0\chi_0.
$$
This is the five-term linear recurrence for $(\chi_n)_{n\geq0}$, which can be written
in the form
\beq \label{RR-ORF}
\pmatrix{\chi_0(z) & \chi_1(z) & \cdots}
\left( \varpi_\cA^*(z) - \varpi_\cA(z) \,\hat\cC \right) = 0,
\eeq
where $\hat\cC$ is the five-diagonal matrix
\beq \label{hatC}
\hat\cC = \pmatrix{
-a_1 &
-\rho_1^- a_2 & \rho_1^- \rho_2^- & 0 & 0 & 0 & \cdots
\cr
\rho_1^+ & -\overline a_1 a_2 & \overline a_1 \rho_2^- & 0 & 0 & 0 &
\cdots
\cr
0 & -\rho_2^+ a_3 & -\overline a_2 a_3 & -\rho_3^- a_4 & \rho_3^- \rho_4^- & 0 & \cdots
\cr
0 & \rho_2^+ \rho_3^+ & \overline a_2 \rho_3^+ & -\overline a_3 a_4 & \overline a_3 \rho_4^-
& 0 & \cdots
\cr
0 & 0 & 0 & -\rho_4^+ a_5 & -\overline a_4 a_5 & -\rho_5^- a_6 & \cdots
\cr
0 & 0 & 0 & \rho_4^+ \rho_5^+ & \overline a_4 \rho_5^+ & -\overline a_5 a_6 & \cdots
\cr
\cdots &\cdots & \cdots & \cdots & \cdots & \cdots & \cdots}.
\eeq
Using (\ref{RHO}) we find that $\hat\cC$ can be related to the unitary five-diagonal matrix
$\cC$ given in (\ref{CMV}) by
\beq \label{hat-C}
\hat\cC = \eta_\cA^{-1} \cC \,\eta_\cA = \cC_\cA.
\eeq
Bearing in mind (\ref{MOB2}) and (\ref{MOB1}), this relation implies that (\ref{RR-ORF}) is
equivalent to
\beq \label{BAND-ORF2}
\pmatrix{\chi_0(z) & \chi_1(z) & \cdots} \left(z-\tilde\zeta_\cA(\cC) \right) = 0,
\eeq
which shows that $\cM$ is invariant under $T_\mu$ and $\tilde\zeta_\cA(\cC)$ is
the matrix representation of $T_\mu \upharpoonright \cM$ with respect to $(\chi_n)_{n\geq0}$.

Similar arguments to those given in the proof of Theorem \ref{HESSORF} prove that $T_\mu
\upharpoonright \cM$ is unitary iff $\cM=L^2_\mu$. However, $T_\mu \upharpoonright \cM$
is unitary whenever $\bsalpha$ is compactly included in $\D$ because in this case
$\tilde\zeta_\cA(\cC)$ is unitary for any sequence $\bsa$ in $\D$. Therefore,
$\cM=L^2_\mu$, i.e., the ORF $(\chi_n)_{n\geq0}$ are a basis of $L^2_\mu$, which implies
that $\tilde\zeta_\cA(\cC)$ is a matrix of the full operator $T_\mu$.

\epr

\br \label{MATRIXchi*}

We know that $(\chi_n)_{n\geq0}$ and $(\chi_{n*})_{n\geq0}$ are basis of $L^2_\mu$ at the
same time, and the corresponding matrices of $T_\mu$ are related by transposition.
Therefore, the previous theorem can be equivalently formulated saying that
$(\chi_{n*})_{n\geq0}$ is a basis of $L^2_\mu$ whenever $\bsalpha$ is compactly included
in $\D$ and, in this case, the related matrix of $T_\mu$ is $\cU^T$. Notice that the
second equalitiy in (\ref{DAG}) implies that $\cU^T=\tilde\zeta_\cA(\cC^T)$ because $\cA$
is diagonal.

\er

Theorem \ref{BANDORF} states that, contrary to the case of the ORF $(\phi_n)_{n\geq0}$,
$(\chi_n)_{n\geq0}$ and $(\chi_{n*})_{n\geq0}$ are basis of $L^2_\mu$ for any measure
$\mu$ on $\T$ if $(\alpha_n)_{n\geq1}$ is compactly included in $\D$. Indeed, the
completeness of $(\chi_n)_{n\geq0}$ and $(\chi_{n*})_{n\geq0}$ in $L^2_\mu$ holds even
under a more general condition for $\bsalpha$, as the next theorem shows. Denoting
$\cF_*=\{f_*:f\in\cF\}$ for any set $\cF$ of complex functions, the problem is to find
sufficient conditions for the equality $\cM=L^2_\mu$ or, equivalently, $\cM_*=L^2_\mu$.

\bp \label{BASIS}

Let $\bsalpha$ be a sequence in $\D$, $\mu$ a measure on $\T$ and $(\chi_n)_{n\geq0}$ the
ORF associated with $(\alpha_1,\hat\alpha_2,\alpha_3,\hat\alpha_4,\dots)$. If
$$
\sum_{k=1}^\infty (1-|\alpha_{2k-1}|) = \sum_{k=1}^\infty (1-|\alpha_{2k}|) = \infty,
$$
then $(\chi_n)_{n\geq0}$ and $(\chi_{n*})_{n\geq0}$ are both basis of $L^2_\mu$.

\ep

\bpr

Given an arbitrary sequence $\bsbeta=(\beta_n)_{n\geq1}$ in $\overline\C\backslash\T$,
let us use the notation $\cL_\infty(\bsbeta)$ for the set of rational functions with
poles in $\hat\bsbeta=(\hat\beta_n)_{n\geq1}$, counted with multiplicity, i.e.,
\beq \label{Linfty}
\ba{c}
\ds\cL_\infty(\bsbeta) = \bigcup_{n\geq2}
\frac{\cP_n}{\varpi_{\beta_1} \cdots \, \varpi_{\beta_{n-1}}},
\ea
\eeq
where $\varpi_\infty(z)=z$. Also, let $\cL(\bsbeta)$ be the closure of
$\cL_\infty(\bsbeta)$ in $L^2_\mu$. Notice that
$\cL(\bsbeta)_*=\cL(\hat\bsbeta)$ and
$\cM=\cL(\alpha_1,\hat\alpha_2,\alpha_3,\hat\alpha_4,\dots)$.

We will show that
\bi
\item[({\it i})]
$\sum_{k=1}^\infty (1-|\alpha_{2k-1}|)=\infty
\;\Rightarrow\;
\{z^j\}_{j\in\N}\subset\cL(\alpha_1,\hat\alpha_2,\alpha_3,\hat\alpha_4,\dots)$,
\item[({\it ii})]
$\sum_{k=1}^\infty (1-|\alpha_{2k}|)=\infty
\;\Rightarrow\;
\{z^{-j}\}_{j\in\N}\subset\cL(\alpha_1,\hat\alpha_2,\alpha_3,\hat\alpha_4,\dots)$.
\ei
This demonstrates the proposition since $\spn\{z^j\}_{j\in\Z}$ is dense in $L^2_\mu$.

Indeed, we only must prove ({\it i}) since it implies ({\it ii}). To see this,
assume that ({\it i}) holds for any sequence $(\alpha_n)_{n\geq1}$ in $\D$. Then,
applying ({\it i}) to the sequence $(\alpha_2,\alpha_1,\alpha_4,\alpha_3,\dots)$ we find
that $\sum_{k=1}^\infty (1-|\alpha_{2k}|)=\infty$ ensures
$\{z^j\}_{j\in\N}\subset\cL(\alpha_2,\hat\alpha_1,\alpha_4,\hat\alpha_3,\dots) =
\cL(\hat\alpha_1,\alpha_2,\hat\alpha_3,\alpha_4,\dots)$, which, applying the
$*$-involution, becomes
$\{z^{-j}\}_{j\in\N}\subset\cL(\alpha_1,\hat\alpha_2,\alpha_3,\hat\alpha_4,\dots)$.

The conditions $\sum_{k=1}^\infty (1-|\alpha_{2k-1}|) = \infty$ and $\sum_{n=1}^\infty
(1-|\alpha_{2k}|) = \infty$ are equivalent respectively to the divergence (to zero) in
$\D$ of the Blaschke products $B^o=\prod_{k=1}^\infty\zeta_{2k-1}$ and
$B^e=\prod_{k=1}^\infty\zeta_{2k}$. Thus, all what we must prove is that the divergence
of $B^o$ implies that $z^j\in\cL(\alpha_1,\hat\alpha_2,\alpha_3,\hat\alpha_4,\dots)$ for
any $j\in\N$.

According to (\ref{Linfty}),
$\{z^j\}_{j\in\N}\subset\cL_\infty(\alpha_1,\hat\alpha_2,\alpha_3,\hat\alpha_4,\dots)=\cM_\infty$
if $\alpha_{2k-1}=0$ for infinitely many values $k\in\N$. Hence, we only need to study
the opposite case that, without loss of generality, we can suppose is
$\alpha_1=\alpha_3=\cdots=\alpha_{2s-1}=0$ and $\alpha_{2k-1}\neq0$ for $k>s$. Then,
$\{z,\dots,z^s\}\subset\cM_\infty$ and, since $\langle f,f \rangle_\mu \leq \|f\|^2_\infty$
for any $f\in\cM_\infty$, it suffices to prove that
$$
\inf_{f\in\cM_n}\|z^j-f(z)\|_\infty \stackrel{n}{\longrightarrow} 0,
\qquad \forall j>s.
$$
To measure the $L^\infty$-distance between a polynomial and a subspace like
$\cM_n$ we can use the following result (see \cite[p. 243]{Ak56} or the more recent
reference \cite[p. 150]{BGHN99}):
$$
\min_{q\in\cP_N} \left\| \frac{z^N+q(z)}{(z-w_1)\cdots(z-w_n)} \right\|_\infty =
\prod_{k=1}^n \frac{1}{\max\{|w_k|,1\}}, \quad w_k\in\C, \quad N \geq n.
$$
Therefore, if $B^o$ diverges, taking $n>s$,
$$
\ba{l}
\kern-3pt \ds
\inf_{ \parbox{35pt}{\scriptsize $f \in \cM_{2n}$ ${}\kern5pt a_k\in\C$} }
\left\| z^{s+m}+a_{m-1}z^{s+m-1}+\cdots+a_1z^{s+1} - f(z) \right\|_\infty =
\cr
\kern17pt \ds
= \inf_{ q\in\cP_{2n+m-1} }
\left\| \frac{z^{2n+m-1}+q(z)}{\prod_{k=1}^{n-1}(z-\alpha_{2k})
\prod_{k=s+1}^n(z-\hat\alpha_{2k-1})} \right\|_\infty
= \prod_{k=s+1}^n |\alpha_{2k-1}| \stackrel{n}{\longrightarrow} 0
\ea
$$
and
$$
\ba{l}
\kern-8pt \ds
\inf_{ \parbox{44pt}{\scriptsize $f \in \cM_{2n+1}$ ${}\kern9pt a_k\in\C$} }
\left\| z^{s+m}+a_{m-1}z^{s+m-1}+\cdots+a_1z^{s+1} - f(z) \right\|_\infty =
\cr
\kern12pt \ds
= \inf_{ q\in\cP_{2n+m} }
\left\| \frac{z^{2n+m}+q(z)}{\prod_{k=1}^n(z-\alpha_{2k})
\prod_{k=s+1}^n(z-\hat\alpha_{2k-1})} \right\|_\infty
= \prod_{k=s+1}^n |\alpha_{2k-1}| \stackrel{n}{\longrightarrow} 0
\ea
$$
for any $m\in\N$. This result implies by induction on $m$ that $z^{s+m}\in\cM$ for any
$m\in\N$.

\epr

In the polynomial case $\cA=0$ so $\cU=\cC$ becomes a five-diagonal matrix. However, in the
general case $\cU$ is not a band matrix but its M\"{o}bius transform $\zeta_\cA(\cU)=\cC$ is
five-diagonal. This fact makes the rational case more complicated than the polynomial one
but, as we will see later, the matrix $\cU$ can be also used in the rational case to
transcribe certain properties of the measure $\mu$ into properties of the corresponding
sequence $\bsa=\cS_\bsalpha(\mu)$.

As in the Hessenberg case, the previous theorem provides a spectral interpretation of the
support of the measure $\mu$. The arguments are similar to those given in the proof of
Theorem \ref{HESSORF1}, but now the restriction $\log\mu' \notin L^1_m$ is not necessary.
Hence, we obtain the following result.

\bt \label{BANDORF1}

Let $\bsalpha$ be a sequence compactly included in $\D$, $\mu$ a measure on $\T$,
$(\phi_n)_{n\geq0}$ the corresponding ORF and $(\chi_n)_{n\geq0}$ the ORF associated with
$(\alpha_1,\hat\alpha_2,\alpha_3,\hat\alpha_4,\dots)$. If
$$
\cU=\tilde\zeta_{\cA}(\cC), \qquad \cA=\cA(\bsalpha), \qquad \cC=\cC(\bsa), \qquad
\bsa=\cS_\bsalpha(\mu),
$$
and $\cE$ is the spectral measure of $\cU$, then $\mu=\cE_{1,1}$. Besides,
$\supp\mu=\sigma(\cU)$ and the mass points of $\mu$ are the eigenvalues of $\cU$, which
have geometric multiplicity 1. $\lambda$ is a mass point iff
$(\chi_n(\lambda))_{n\geq0}\in\ell^2$. Given a mass point $\lambda$, the corresponding
eigenvectors of $\cU$ are spanned by $\pmatrix{\chi_0(\lambda) & \chi_1(\lambda) &
\cdots}^\dag$ and
$\mu(\{\lambda\})=\left(\sum_{n=0}^\infty|\chi_n(\lambda)|^2\right)^{-1}
=\left(\sum_{n=0}^\infty|\phi_n(\lambda)|^2\right)^{-1}$.

\et

Analogously to the Hessenberg case, we can formulate the above spectral results in terms
of pairs of band operators. Theorem \ref{HESSORF1} implies that
$\supp\mu=\sigma(\tilde\varpi_\cA^*(\cC),\tilde\varpi_\cA(\cC))$ and the mass points of
$\mu$ are the eigenvalues of the five-diagonal pair
$(\tilde\varpi_\cA^*(\cC),\tilde\varpi_\cA(\cC))$. Also, given a mass point $\lambda$,
the left eigenvectors of $(\tilde\varpi_\cA^*(\cC),\tilde\varpi_\cA(\cC))$ are spanned by
$\pmatrix{\chi_0(\lambda) & \chi_1(\lambda) & \cdots} \eta_\cA^{-1/2}$. Furthermore, the
factorization $\cC=\cC_o\cC_e$ makes possible to formulate the above results using the
tridiagonal pair
$$
(\tilde\varpi_\cA^*(\cC),\tilde\varpi_\cA(\cC))\,\cC_e^\dag =
(\cC_o+\cA\cC_e^\dag,\cC_e^\dag+\cA^\dag\cC_o)
$$
instead of the five-diagonal pair.

\subsection{Zeros of ORF and five-diagonal matrices} \label{Z-ORF-5M}

The previous results suggest that it should be possible a spectral interpretation of the
zeros of ORF and PORF in terms of five-diagonal matrices. Similarly to the Hessenberg
case, an important ingredient for this is the orthogonal truncation of $\zeta_n(T_\mu)$
on $\cM_n$. Taking into account that $\cM_n=B_{l*}^e\cL_n$ with $l=[(n-1)/2]$, the
following generalization of Theorem \ref{ORF-ZEROS} is of interest to relate this
truncation to the zeros of the ORF $\phi_n$. This generalization deals with the
orthogonal truncation of $\zeta_n(T_\mu)$ on $h\cL_n$, $h \in L^2_\mu$, which is given by
$\zeta_n(T_\mu)^{(h\cL_n)} = L_n^h \zeta_n(T_\mu) \upharpoonright h\cL_n$, where $L_n^h
\colon L^2_\mu \to L^2_\mu$ is the orthogonal projection on $h\cL_n$.

\bt \label{ORF-ZEROS2}

Let $\bsalpha$ be an arbitrary sequence in $\D$, $\mu$ a measure on $\T$ and
$\phi_n=p_n/\pi_n$ the related $n$-th ORF. Then, for any Borel function $h\colon\T\to\T$:
\be
\item If $Z_n$ is the set of zeros of $\phi_n$, $\zeta_n(Z_n)$ is the set of eigenvalues
      of $\zeta_n(T_\mu)^{(h\cL_n)}$, and these eigenvalues have geometric multiplicity 1.
\item If $p_n(z)\propto\prod_{k=1}^n(z-\lambda_k)$, the characteristic polynomial of
      $\zeta_n(T_\mu)^{(h\cL_n)}$ is
      $$
      \prod_{k=1}^n(z-\zeta_n(\lambda_k)).
      $$
\ee

\et

\bpr

The operator multiplication by $h$ in $L^2_\mu$
$$
h(T_\mu) \colon \mathop{L^2_\mu \to L^2_\mu} \limits_{f \; \longrightarrow \; hf}
$$
is unitary because $h$ maps $\T$ on itself. When restricted in the following way
$$ V \colon \mathop{\cL_n \to h\cL_n} \limits_{f \kern5pt \longrightarrow \kern5pt hf}
$$
it yields a isometric isomorphism $V$ between $\cL_n$ and $h\cL_n$. The orthogonal
projection on $h\cL_n$ is $L_n^h=h(T_\mu)L_nh(T_\mu)^\dag$, where $L_n$ is the orthogonal
projection on $\cL_n$. Thus, the orthogonal truncations of $\zeta_n(T_\mu)$ on $h\cL_n$
and $\cL_n$ are related by $\zeta_n(T_\mu)^{(h\cL_n)}=V\zeta_n(T_\mu)^{(\cL_n)}V^{-1}$,
so they are unitarily equivalent. In consequence, they have the same eigenvalues and with
the same geometric and algebraic multiplicity. Hence, the result follows from Theorem
\ref{ORF-ZEROS}.

\epr

Taking $h=B_{l*}^e$, $l=[(n-1)/2]$, in the previous theorem we find that it holds for the
orthogonal truncation of $\zeta_n(T_\mu)$ on $\cM_n$, i.e.,
$\zeta_n(T_\mu)^{(\cM_n)}=M_n\zeta_n(T_\mu)\upharpoonright\cM_n$, where $M_n \colon
L^2_\mu \to L^2_\mu$ is the orthogonal projection on $\cM_n$. To give a matrix version of
this result we simply need a matrix representation of $\zeta_n(T_\mu)^{(\cM_n)}$.

\bt \label{BANDORFmatrix}

Let $\bsalpha$ be an arbitrary sequence in $\D$, $\mu$ a measure on $\T$ and
$(\chi_n)_{n\geq0}$ the ORF associated with
$(\alpha_1,\hat\alpha_2,\alpha_3,\hat\alpha_4,\dots)$. If $\cA=\cA(\bsalpha)$ and
$\cC=\cC(\bsa)$ with $\bsa=\cS_\bsalpha(\mu)$, the matrix of $\zeta_n(T_\mu)^{(\cM_n)}$
with respect to $(\chi_k)_{k=0}^{n-1}$ is $\zeta_n(\cU^{(n)})$, where
$$
\cU^{(n)}=\tilde\zeta_{\cA_n}(\cC_n), \qquad \|\cU^{(n)}\|=1.
$$

\et

\bpr

From the factorization $\cC = \cC_o\cC_e$ we find that $\cC_n = \cC_{on}\cC_{en}$. Only
one among the factors $\cC_{on}$ and $\cC_{en}$ is unitary, but the norm of the remaining
factor is 1, so $\|\cC_n\|=1$. Since $\tilde\zeta_{\cA_n}$ leaves $\T_{\C^n}$ invariant,
$\cU^{(n)}=\tilde\zeta_{\cA_n}(\cC_n)$ is a well defined matrix with $\|\cU^{(n)}\|=1$.
The same holds for $\zeta_n(\cU^{(n)})$.

To prove that $\zeta_n(\cU^{(n)})$ is the matrix representation of $\zeta_n(T_\mu)$ with
respect to $(\chi_k)_{k=0}^{n-1}$, let us consider first an odd $n$. Then, the first $n$
equations of (\ref{RR-ORF}) can be written as
$$
\pmatrix{\chi_0(z) & \kern-5pt \cdots \kern-5pt & \chi_{n-1}(z)}
\left( \varpi^*_{\cA_n}(z) - \varpi_{\cA_n}(z) \,\hat\cC_n \right) =
b_n \varpi_n(z) \chi_n(z), \quad b_n\in\C^n.
$$
(\ref{hat-C}) gives $\hat\cC_n=\eta_{\cA_n}^{-1}\cC_n\eta_{\cA_n}=(\cC_n)_{\cA_n}$. This,
together with identities (\ref{MOB2}) and (\ref{MOB1}), yields
$$
\pmatrix{\chi_0(z) & \cdots & \chi_{n-1}(z)} \left( z - \cU^{(n)} \right) = c_n
\varpi_n(z) \chi_n(z), \quad c_n\in\C^n.
$$
Using (\ref{MOB3}) we get
$$
\pmatrix{\chi_0(z) & \cdots & \chi_{n-1}(z)} \left( \zeta_n(z) - \zeta_n(\cU^{(n)})
\right) = d_n \chi_n(z), \quad d_n\in\C^n,
$$
and, taking into account that
$$
M_n\chi_k=\cases{\chi_k & if $k < n$, \cr 0 & if $k \geq n$,}
$$
we finally obtain
$$
\pmatrix{M_n\zeta_n\chi_0 & \cdots & M_n\zeta_n\chi_{n-1}} =
\pmatrix{\chi_0 & \cdots & \chi_{n-1}} \zeta_n(\cU^{(n)}).
$$
This equality proves that $\zeta_n(\cU^{(n)})$ is the matrix of $M_n \zeta_n(T_\mu)
\upharpoonright \cM_n$ with respect to $(\chi_k)_{k=0}^{n-1}$.

On the other hand, if $n$ is even, we consider the orthogonal truncation of
$\zeta_n(T_\mu)$ on $\cM_{n*}$, i.e., $\zeta_n(T_\mu)^{(\cM_{n*})} = M_{n*}
\zeta_n(T_\mu) \upharpoonright \cM_{n*}$, where $M_{n*} \colon L^2_\mu \to L^2_\mu$ is
the orthogonal projection on $\cM_{n*}$. Taking into account that $\cC$ is unitary, the
identity obtained by applying the $*$-involution on (\ref{RR-ORF}) reads
\beq \label{RR-ORF1*}
\pmatrix{\chi_{0*}(z) & \chi_{1*}(z) & \cdots}
\left( \varpi_\cA^*(z) - \varpi_\cA(z) \,(\cC^T)_\cA \right) = 0.
\eeq
A similar reasoning starting from the first $n$ equations of this equality proves that
the matrix of $\zeta_n(T_\mu)^{(\cM_{n*})}$ with respect to $(\chi_{k*})_{k=0}^{n-1}$ is
$\zeta_n(\cU^{(n)}_*)$, where $\cU^{(n)}_*=\tilde\zeta_{\cA_n}(\cC_n^T)$. Notice that
(\ref{DAG}) implies that $\cU^{(n)}_*=\cU^{(n)T}$ because $\cA_n$ is diagonal.

Given a measure $\mu$ on $\T$, for any $n\in\N$, the subspace $\cM_n$ only depends on the
parameters $\alpha_1,\dots,\alpha_{n-1}$ of the sequence $\bsalpha$, so the same holds
for the orthogonal truncations $\zeta_n(T_\mu)^{(\cM_n)}$ and
$\zeta_n(T_\mu)^{(\cM_{n*})}$. Therefore, concerning the spectral properties of these
truncations we can suppose without loss of generality that $\bsalpha$ is compactly
supported on $\D$. Then, the matrix representations of $T_\mu$ with respect to
$(\chi_k)_{k\geq0}$ and $(\chi_{k*})_{k\geq0}$ are $\cU$ and $\cU^T$ respectively. Hence,
the representations of the orthogonal truncations $\zeta_n(T_\mu)^{(\cM_n)}$ and
$\zeta_n(T_\mu)^{(\cM_{n*})}$ with respect to $(\chi_k)_{k=0}^{n-1}$ and
$(\chi_{k*})_{k=0}^{n-1}$ are the principal submatrices $(\zeta_n(\cU))_n$ and
$(\zeta_n(\cU^T))_n$ respectively. The fact that $(\zeta_n(\cU^T))_n =
\left(\zeta_n(\cU)\right)_n^T$ implies that, when $n$ is even, the matrix of
$\zeta_n(T_\mu)^{(\cM_n)}$ with respect to $(\chi_k)_{k=0}^{n-1}$ is
$\zeta_n(\cU^{(n)}_*)^T=\zeta_n(\cU^{(n)})$.

\epr

\br

The proof of the previous theorem also shows that the matrix of the orthogonal truncation
$\zeta_n(T_\mu)^{(\cM_{n*})}$ with respect to $(\chi_{k*})_{k=0}^{n-1}$ is
$\zeta_n(\cU^{(n)})^T$.

\er

As a consequence of Theorems \ref{ORF-ZEROS2} and \ref{BANDORFmatrix} we have the
following spectral interpretation of the zeros of ORF in terms of M\"{o}bius transformations
of five-diagonal matrices.

\bt \label{BANDORFzeros}

Let $\bsalpha$ be an arbitrary sequence in $\D$, $\mu$ a measure on $\T$,
$(\phi_n)_{n\geq0}$ the corresponding ORF and $(\chi_n)_{n\geq0}$ the ORF
associated with the sequence $(\alpha_1,\hat\alpha_2,\alpha_3,\hat\alpha_4,\dots)$. Let
$\cA=\cA(\bsalpha)$ and $\cC=\cC(\bsa)$ with $\bsa=\cS_\bsalpha(\mu)$. If
$\cU^{(n)}=\tilde\zeta_{\cA_n}(\cC_n)$, then:
\be
\item The zeros of $\phi_n$ are the eigenvalues of $\cU^{(n)}$, which have geometric
      multiplicity 1. If $\lambda$ is a zero of $\phi_n$, the related left and right
      eigenvectors of $\cU^{(n)}$ are spanned by $X_n(\lambda)$ and $Y_n(\lambda)^T$
      respectively, where
      $$
      X_n=B^e_{[\frac{n-1}{2}]}\pmatrix{\chi_0&\cdots&\chi_{n-1}}, \qquad
      Y_n=B^o_{[\frac{n}{2}]}\pmatrix{\chi_{0*}&\cdots&\chi_{n-1*}}.
      $$
\item $\ds \phi_n=\frac{p_n}{\pi_n}$ with $p_n$ proportional to the characteristic
      polynomial of $\cU^{(n)}$.
\ee

\et

\bpr

The vectors $X_n(z)$ and $Y_n(z)$ are rational functions with the poles lying on $\E$, so
they can be evaluated at any zero $\lambda$ of $\phi_n$ since $\lambda\in\D$. Besides,
$X_n(\lambda),Y_n(\lambda)\neq0$ because $B^e_k\chi_{2k} = \phi_{2k}^*$,
$B^o_k\chi_{2k-1*}=\phi_{2k-1}^*$ and $\phi_n^*$ has its zeros in $\E$.

The proof of the theorem is similar to the case of Theorem \ref{HESSORFzeros}, the only
difference concerning the identification of the eigenvectors. To obtain the left eigenvectors
of $\cU^{(n)}$ let us consider the first $n$ equations of (\ref{RR-ORF}) for an arbitrary
$n\in\N$. These equations can be written
as
\beq \label{equationsZEROS}
\ba{l}
\pmatrix{\chi_0(z) & \kern-5pt \cdots \kern-5pt & \chi_{n-1}(z)}
\left( \varpi^*_{\cA_n}(z) - \varpi_{\cA_n}(z) \,\hat\cC_n \right) =
\medskip \cr \kern150pt
= b_n \varpi_n(z) \chi_n(z) + d_n \varpi_{n+1}(z) \chi_{n+1}(z),
\ea
\eeq
where $b_n,d_n\in\C^n$ are
$$
\ba{l}
b_n = \cases{
\rho_n^+ \pmatrix{0 & \cdots & 0 & \rho_{n-1}^+ & \overline a_{n-1}} & odd $n$,
\smallskip \cr
- \rho_n^+ a_{n+1} \pmatrix{0 & \cdots & 0 & 1} & even $n$,
}
\medskip \cr
d_n = \cases{
0 & odd $n$,
\smallskip \cr
\rho_n^+ \rho_{n+1}^+ \pmatrix{0 & \cdots & 0 & 1} & even $n$.
}
\ea
$$
Writing $\chi_n$ and $\chi_{n*}$ in terms of $\phi_n$ and $\phi_n^*$ with the aid of
(\ref{LORF}) and (\ref{LORF*}), and using the first equation of (\ref{RR1}) in the case of
even $n$, (\ref{equationsZEROS}) reads
\beq
\label{equationsZEROS2}
\kern-11pt \ba{l}
\pmatrix{\chi_0(z) & \kern-5pt \cdots \kern-5pt & \chi_{n-1}(z)}
\kern-1pt \left( \varpi^*_{\cA_n}(z) - \varpi_{\cA_n}(z) \,\hat\cC_n \right) \kern-2pt =
\rho_n^+ \varpi_n(z) B^e_{l*}(z) \phi_n(z) \kern1pt v_n,
\medskip \cr \kern5pt
l=\left[\frac{n-1}{2}\right], \quad v_n\in\C^n, \quad
v_n = \cases{
\pmatrix{0 & \cdots & 0 & \rho_{n-1}^+ & \overline a_{n-1}} & odd $n$,
\smallskip \cr
\pmatrix{0 & \cdots & 0 & 1} & even $n$.
}
\ea
\eeq
Also, remember that (\ref{MOB2}), (\ref{MOB1}) and $\hat\cC_n=\eta_{\cA_n}^{-1}\cC_n\eta_{\cA_n}$
imply that
\beq \label{MOBIUS-PAIR2}
(z-\cU^{(n)}) \,\tilde\varpi_{\cA_n}(\hat\cC_n) =
\varpi^*_{\cA_n}(z) - \varpi_{\cA_n}(z) \,\hat\cC_n.
\eeq
Therefore, if $\lambda$ is a zero of $\phi_n$, (\ref{equationsZEROS2}) and
(\ref{MOBIUS-PAIR2}) show that $X_n(\lambda)$ is a left eigenvector of $\cU^{(n)}$ with
eigenvalue $\lambda$.

Proceeding in a similar way with the first $n$ equations of (\ref{RR-ORF1*}) we find that
$Y_n(\lambda)$ is a left eigenvector of $\cU^{(n)T}$ with eigenvalue $\lambda$ for any
zero $\lambda$ of $\phi_n$. Therefore, $Y_n(\lambda)^T$ is a right eigenvector of
$\cU^{(n)}$.

\epr

For a unitary matrix, like $\cV$ in the case $\bsa\notin\ell^2$, $\cV^{(n;u)}$ or $\cU$,
left and right eigenvectors are related by the $\dag$-operation. However, this is not the
case of the matrices $\cV^{(n)}$ or $\cU^{(n)}$. Theorem \ref{HESSORFzeros} only gives
information about the left eigenvectors of $\cV^{(n)}$, while Theorem \ref{BANDORFzeros}
provides both, the left and right eigenvectors of $\cU^{(n)}$. Apart from the simplest
form of $\cU^{(n)}$, this is another advantage of using this matrix instead of
$\cV^{(n)}$ for the spectral representation of the zeros of ORF.

Concerning the form of the eigenvectors of $\cU^{(n)}$, notice that the factors
$B^e_{[(n-1)/2]}$ and $B^o_{[n/2]}$ in $X_n$ and $Y_n$ are necessary to avoid any problem
when evaluating them on a point of $\D$. However, if a zero $\lambda$ of $\phi_n$ does
not coincide with any $\alpha_k$ for $k=1,\dots,n-1$, then we can take as left and right
eigenvectors $\pmatrix{\chi_0(\lambda),\dots,\chi_{n-1}(\lambda)}$ and
$\pmatrix{\chi_{0*}(\lambda),\dots,\chi_{n-1*}(\lambda)}^T$ respectively.

As in the Hessenberg case, there are other alternatives to express $p_n$ as a
determinant. Indeed, from (\ref{MOB2}), (\ref{MOB1}) and the identity
$\hat\cC_n=\eta_{\cA_n}^{-1}\cC_n\eta_{\cA_n}$,
$$
p_n(z) \propto \det\left(\varpi^*_{\cA_n}(z) - \varpi_{\cA_n}(z) \,\cC_n\right)
= \det\left(z
\,\tilde\varpi_{\cA_n}(\cC_n)-\tilde\varpi_{\cA_n}^*(\cC_n)\right).
$$
So $p_n$ can be calculated as a determinant of a five-diagonal matrix. Furthermore, the
last expression shows that the zeros of $\phi_n$ are the eigenvalues of the five-diagonal
pair $(\tilde\varpi_{\cA_n}^*(\cC_n),\tilde\varpi_{\cA_n}(\cC_n))$. The associated left
eigenvectors with eigenvalue $\lambda$ are spanned by
$X_n(\lambda)\,\eta_{\cA_n}^{-1/2}$.

Besides, the factorization $\cC_n=\cC_{on}\cC_{en}$ permits us to express $p_n$
as a determinant of a tridiagonal matrix. If $n$ is odd, $\cC_{en}$ is unitary, thus
$$
p_n(z) \propto \det\left(z(\cC_{en}^\dag+\cA_n^\dag\cC_{on})
-(\cC_{on}+\cA_n\cC_{en}^\dag)\right)
$$
and the zeros of $\phi_n$ are the eigenvalues of the tridiagonal pair
$$
(\cC_{on}+\cA_n\cC_{en}^\dag,\cC_{en}^\dag+\cA_n^\dag\cC_{on}),
$$
which has the same left eigenvectors as
$(\tilde\varpi_{\cA_n}^*(\cC_n),\tilde\varpi_{\cA_n}(\cC_n))$. On the contrary,
$\cC_{on}$ is unitary for an even $n$. In this situation we can use the fact that
$\cU^{(n)}$ and $\cU^{(n)T}=\tilde\zeta_{\cA_n}(\cC_n^T)$ have the same characteristic
polynomial, and the left eigenvectors of one of them are the transposed of the right
eigenvectors of the other one. Hence,
$$
p_n(z) \propto
\det\left(\varpi^*_{\cA_n}(z) - \varpi_{\cA_n}(z) \,\cC_n^T\right) =
\det\left(z \,\tilde\varpi_{\cA_n}(\cC_n^T)-\tilde\varpi_{\cA_n}^*(\cC_n^T)\right)
$$
and, bearing in mind that $\cC_n^T=\cC_{en}\cC_{on}$,
$$
p_n(z) \propto \det\left(z(\cC_{on}^\dag+\cA_n^\dag\cC_{en})
-(\cC_{en}+\cA_n\cC_{on}^\dag)\right).
$$
So, the zeros of $\phi_n$ are the eigenvalues of the tridiagonal pair
$$
(\cC_{en}+\cA_n\cC_{on}^\dag,\cC_{on}^\dag+\cA_n^\dag\cC_{en})
$$
and the left eigenvectors with eigenvalue $\lambda$ are spanned by
$Y_n(\lambda)\,\eta_{\cA_n}^{-1/2}$.

\medskip

The zeros of the PORF $Q_n^v$ have a spectral interpretation in terms of band matrices
too. Such an interpretation has to do with the matrix representation of $T_{\mu_n^v}$
with respect $(\chi_k)_{k=0}^{n-1}$, which is an orthonormal basis of $L^2_{\mu_n^v}$ due
to the exactness of the quadrature formulas associated with $\mu_n^v$. Similar arguments
to those appearing before Theorem \ref{HESSPORFzeros} show that the zeros of the PORF
should be related to the unitary matrix $\cC_n^u$ obtained from $\cC_n$ when substituting
the parameter $a_n\in\D$ by $u\in\T$. More precisely, we have the following result, which
can be understood as a limit case of Theorems \ref{BANDORFmatrix} and \ref{BANDORFzeros}.

\bt \label{BANDPORFzeros}

Let $\bsalpha$ be an arbitrary sequence in $\D$, $\mu$ a measure on $\T$,
$(\phi_n)_{n\geq0}$ the corresponding ORF, $\cA=\cA(\bsalpha)$ and $\cC=\cC(\bsa)$ with
$\bsa=\cS_\bsalpha(\mu)$. If $Q_n^v=\phi_n+v\phi_n^*$ is the $n$-th PORF related to
$v\in\T$ and $\mu_n^v$ is the associated measure, then:
\be
\item The matrix of $T_{\mu_n^v}$ with respect to $(\chi_k)_{k=0}^{n-1}$ is
      $$
      \cU^{(n;u)}=\tilde\zeta_{\cA_n}(\cC_n^u), \qquad u=\tilde\zeta_{a_n}(v).
      $$
\item The zeros of $Q_n^v$ are the eigenvalues of $\cU^{(n;u)}$.
      If $\lambda$ is a zero of $Q_n^v$, the related eigenvectors of $\cU^{(n;u)}$ are
      spanned by $\pmatrix{\chi_0(\lambda)&\cdots&\chi_{n-1}(\lambda)}^\dag$.
\item $\ds Q_n^v=\frac{q_n^v}{\pi_n}$ with $q_n^v$ proportional to the characteristic
      polynomial of $\cU^{(n;u)}$.
\ee

\et

\bpr

As in the case of Theorem \ref{HESSPORFzeros}, it suffices to prove item 1. For an odd
$n=2l+1$, using (\ref{PORFu}) in a similar computation to that of (\ref{RRchi}) gives
$$
\ba{l}
\ds \varpi_{n-2}^* \chi_{n-2} =
\rho_{n-1}^+ \rho_n^+ \varpi_n B^e_{[n/2]*} \frac{Q_n^v}{1+\overline a_nv} -
\rho_{n-1}^+ u \, \varpi_{n-1} \chi_{n-1} \, -
\smallskip \cr \kern65pt
- \, \overline a_{n-2}a_{n-1} \varpi_{n-2} \chi_{n-2} -
\rho_{n-2}^- a_{n-1} \varpi_{n-3} \chi_{n-3},
\medskip \cr
\ds \varpi_{n-1}^* \chi_{n-1} =
\overline a_{n-1} \rho_n^+ \varpi_n B^e_{[n/2]*} \frac{Q_n^v}{1+\overline a_nv} -
\overline a_{n-1} u \, \varpi_{n-1} \chi_{n-1} \, +
\smallskip \cr \kern65pt
+ \, \overline a_{n-2} \rho_{n-1}^- \varpi_{n-2} \chi_{n-2} +
\rho_{n-2}^- \rho_{n-1}^- \varpi_{n-3} \chi_{n-3},
\ea
$$
where $u=\tilde\zeta_{a_n}(v)$. These relations can be combined with the first $n-2$
equations of (\ref{RR-ORF}) in the matrix identity
$$
\ba{l}
\pmatrix{\chi_0(z) & \kern-5pt \cdots \kern-5pt & \chi_{n-1}(z)}
\left( \varpi^*_{\cA_n}(z) - \varpi_{\cA_n}(z) \,\hat\cC_n^u \right) =
\smallskip \cr \kern190pt
= b_n \varpi_n(z) B^e_{l*}(z) \, Q_n^v(z),
\quad b_n\in\C^n,
\ea
$$
with $\hat\cC_n^u=\eta_{\cA_n}^{-1}\cC_n^u\eta_{\cA_n}=(\cC_n^u)_{\cA_n}$. Thus,
using (\ref{MOB2}) and (\ref{MOB1}) we find that
$$
\pmatrix{\chi_0(z) & \kern-5pt \cdots \kern-5pt & \chi_{n-1}(z)}
\left( z - \cU^{(n;u)} \right) =
c_n \varpi_n(z) B^e_{l*}(z) \, Q_n^v(z), \quad c_n\in\C^n.
$$
Therefore, $\cU^{(n;u)}$ is the matrix of $T_{\mu_n^v}$ with respect to
$(\chi_k)_{k=0}^{n-1}$ because $Q_n^v=0$ in $L^2_{\mu_n^v}$.

On the other hand, if $n=2l$ is even, proceeding in a similar way with
(\ref{PORFu}) and (\ref{RR-ORF}) we arrive at
$$
\pmatrix{\chi_{0*}(z) & \kern-5pt \cdots \kern-5pt & \chi_{n-1*}(z)}
\left( z - {\cU^{(n;u)}}^T \right) =
c_n \varpi_n(z) B^o_{l*}(z) \, Q_n^v(z), \quad c_n\in\C^n,
$$
so ${\cU^{(n;u)}}^T$ is the matrix of $T_{\mu_n^v}$ with respect to
$(\chi_{k*})_{k=0}^{n-1}$. Consequently, the matrix of $T_{\mu_n^v}$ with respect to
$(\chi_k)_{k=0}^{n-1}$ is $\cU^{(n;u)}$.

\epr

The zeros of a PORF can be also interpreted as eigenvalues of a pair of band matrices.
If $u=\tilde\zeta_{a_n}(v)$,
$$
q_n^v(z) \propto \det\left(\varpi^*_{\cA_n}(z) - \varpi_{\cA_n}(z) \,\cC_n^u\right) =
\det\left(z \,\tilde\varpi_{\cA_n}(\cC_n^u)-\tilde\varpi_{\cA_n}^*(\cC_n^u)\right)
$$
gives $q_n^v$ as a determinant of a five-diagonal matrix. The zeros of $Q_n^v$ are the
eigenvalues of the five-diagonal pair
$(\tilde\varpi_{\cA_n}^*(\cC_n^u),\tilde\varpi_{\cA_n}(\cC_n^u))$ and, given an
eigenvalue $\lambda$, $\pmatrix{\chi_0(\lambda)&\cdots&\chi_{n-1}(\lambda)}
\eta_{\cA_n}^{-1/2}$ spans the corresponding left eigenvectors subspace.

We have also a factorization $\cC_n^u=\cC_{on}^u\cC_{en}^u$, where $\cC_{on}^u$ and
$\cC_{en}^u$ are the result of substituting $a_n$ by $u$ in $\cC_{on}$ and $\cC_{en}$
respectively (this substitution actually takes place only in $\cC_{on}$ or $\cC_{en}$,
depending whether $n$ is odd or even). $\cC_{on}^u$ and $\cC_{en}^u$ are both unitary, so
$$
p_n(z) \propto \det\left(z(\cC_{en}^{u\dag}+\cA_n^\dag\cC_{on}^u)
-(\cC_{on}^u+\cA_n\cC_{en}^{u\dag})\right)
$$
and the zeros of $\phi_n$ are the eigenvalues of the tridiagonal pair
$$
(\cC_{on}^u+\cA_n\cC_{en}^{u\dag},\cC_{en}^{u\dag}+\cA_n^\dag\cC_{on}^u),
$$
which has the same left eigenvectors as
$(\tilde\varpi_{\cA_n}^*(\cC_n^u),\tilde\varpi_{\cA_n}(\cC_n^u))$.

\section{Applications} \label{APPL}

In this section we will present some applications of the spectral theory previously
developed for the ORF on the unit circle. We will use the results involving five-diagonal
matrices due to their advantages. The corresponding spectral theory associates with each
sequence of ORF a five-diagonal unitary matrix $\cC(\bsa)$ depending on the parameters
$\bsa=(a_n)_{n\geq1}$ of the recurrence relation, and a diagonal matrix $\cA(\bsalpha)$
depending on the sequence $\bsalpha=(\alpha_n)_{n\geq1}$ which defines the poles
$\hat\alpha_n$. These band matrices keep all the information about the ORF since they
generate the full sequence of ORF through the associated recurrence. The importance of
these matrices is that they play the role of a simple short cut that connects directly
the parameters $\bsa$, $\bsalpha$ to the ORF and the related orthogonality measure.

An essential difference with the polynomial case is that the matrix directly related to
the ORF and the orthogonality measure is not the five-diagonal one, but an operator
M\"{o}bius transform of it, namely,
$\cU(\bsa,\bsalpha)=\tilde\zeta_{\cA(\bsalpha)}(\cC(\bsa))$. This introduces important
difficulties when trying to apply the spectral theory to the rational case. However, in
spite of these difficulties, the matrix tool $\cU(\bsa,\bsalpha)$ becomes powerful enough
to deal with hard problems even in the rational case. To understand the scope of the
rational spectral theory, we will use it to solve some non trivial problems about the
relation between the behavior of the sequences $\bsa$, $\bsalpha$ and the properties of
the corresponding orthogonality measure $\mu(\bsa,\bsalpha)$. The answers to these
problems are known for OP, but the generalizations to ORF are new.

The strategy will be to apply standard results of perturbation theory to the unitary
operator on $\ell^2$ defined by the matrix $\cU(\bsa,\bsalpha)$. We will apply such
perturbation results to the comparison of $\cU(\bsa,\bsalpha)$ with another normal
matrix, eventually with the form $\cU(\bsb,\bsbeta)$. A useful remark for these
comparisons is that, for $\bsbeta$ compactly supported in $\D$, $\cU(\bsb,\bsbeta)$
defines a unitary operator for any sequence $\bsb$ in $\overline\D$ since, then,
$\cC(\bsb)$ is unitary. However, $\cU(\bsb,\bsbeta)$ only represents a multiplication
operator on $\T$ when $\bsb$ lies on $\D$. When $b_n\in\T$ for some $n$ we know that
$\cC(\bsb)$ decomposes as a direct sum of an $n \times n$ and an infinite matrix. Taking
into account that $\cA(\bsbeta)$ is diagonal, a similar decomposition holds for
$\cU(\bsb,\bsbeta)$.

The results of operator theory that we will apply state that two operators $T,S$ on $H$
have some common spectral property provided that the perturbation $T-S$ belongs to
certain class of operators. We will deal with two kinds of perturbations: compact and
trace class operators. Both are subsets of $\B_H$ that are closed under sum, left and
right product by any element of $\B_H$ and also under the $\dag$-operation, that is, they
are hermitian ideals of $\B_H$. This fact is the key that permit us to use techniques of
band matrices in the spectral theory of ORF, according to the following result.

\bp \label{IDEAL}

Let $\frI$ be a hermitian ideal of $\B_H$. If $A,B\in\D_H$ are normal and $AB=BA$, the
condition $A-B\in\frI$ implies the equivalences
$$
T-S\in\frI \;\Leftrightarrow\; \zeta_A(T)-\zeta_B(S)\in\frI \;\Leftrightarrow\;
\tilde\zeta_A(T)-\tilde\zeta_B(S)\in\frI, \quad \forall\,T,S\in\overline\D_H.
$$

\ep

\bpr

It suffices to prove the first equivalence because $\tilde\zeta_A=\zeta_{-A}$. Let $\frI$
be a hermitian ideal of $\B_H$. The identities
$$
T_1T_2-S_1S_2=(T_1-S_1)\,T_2+S_1(T_2-S_2),
\quad
T^{-1}\kern-3pt-S^{-1}=-T^{-1}\kern-1pt(T-S)S^{-1}
$$
prove that
$$
\ba{l}
T_i,S_i\in\B_H, \kern7pt T_i-S_i\in\frI \kern7pt\Rightarrow\kern7pt T_1S_1-T_2S_2\in\frI,
\smallskip \cr
T^{-1},S^{-1}\in\B_H, \kern7pt T-S\in\frI \kern7pt\Rightarrow\kern7pt T^{-1}\kern-3pt-S^{-1}\in\frI.
\ea
$$
Suppose now $A,B\in\D_H$ normal such that $AB=BA$ and $A-B\in\frI$. Then
$\eta_A^2-\eta_B^2=BB^\dag-AA^\dag\in\frI$. The functional calculus for normal
operators shows that $\eta_A\eta_B=\eta_B\eta_A$, so
$\eta_A-\eta_B=(\eta_A+\eta_B)^{-1}(\eta_A^2-\eta_B^2)\in\frI$ since
$(\eta_A+\eta_B)^{-1}\in\B_H$ because $\eta_A$ and $\eta_B$ are positive with bounded
inverse. If, besides, $T,S\in\overline\D_H$ are such that $T-S\in\frI$, then
$\varpi_A(T)-\varpi_B(S)=SB^\dag-TA^\dag\in\frI$ and
$\varpi_A^*(T)-\varpi_B^*(S)=T-S+B-A\in\frI$. In consequence,
$T-S\in\frI\Rightarrow\zeta_A(T)-\zeta_B(S)\in\frI$. Substituting in this result $A,B$ by
$-A,-B$ and $T,S$ by $\zeta_A(T),\zeta_B(S)$ respectively, we also find the opposite
inclusion.

\epr

Taking into account that $\cA(\bsalpha)$ is diagonal, the above result has the following
immediate consequence.

\bc \label{IDEAL-U}

If $\frI$ is a hermitian ideal of $\B_\ell^2$ and $\bsalpha,\bsbeta$ are sequences
compactly included in $\D$, the condition $\cA(\bsalpha)-\cA(\bsbeta)\in\frI$ implies the
equivalence
$$
\cC(\bsa)-\cC(\bsb)\in\frI \;\Leftrightarrow\; \cU(\bsa,\bsalpha)-\cU(\bsb,\bsbeta)\in\frI
$$
for any sequences $\bsa,\bsb$ in $\overline\D$.

\ec

Besides, from the factorization $\cC(\bsa)=\cC_o(\bsa)\,\cC_e(\bsa)$ we find that, for
any ideal $\frI$ of $\B_{\ell^2}$,
$$
\cC_o(\bsa)-\cC_o(\bsb),\;\cC_e(\bsa)-\cC_e(\bsb)\in\frI \;\Rightarrow\;
\cC(\bsa)-\cC(\bsb)\in\frI.
$$
In fact, many perturbation results for the five diagonal matrix $\cC(\bsa)$ are known due
to the extensive use of this matrix during the last years for the spectral analysis of OP
on $\T$.

The perturbation results that we will use are the invariance of the essential spectrum
for normal operators related by a compact perturbation (Weyl's theorem: see \cite{We09}
and \cite{Be71,Si71}), and the invariance of the absolutely continuous spectrum for
unitary operators related by a trace class perturbation (Kato-Birman theorem: see
\cite{Ka57,Bi62} and \cite{BiKr62}). Given an operator $T$, its essential and absolutely
continuous spectrum will be denoted $\sigma_e(T)$ and $\sigma_{ac}(T)$ respectively. In
the case of a normal operator, $\sigma_e(T)$ is constituted by the limit points of
$\sigma(T)$ and the eigenvalues with infinite geometric multiplicity. In particular, for
any measure $\mu$ on $\T$, $\sigma_e(T_\mu)$ is the set $\{\supp\mu\}'$ of limit points
of $\supp\mu$ and $\sigma_{ac}(T_\mu)$ is the support of the absolutely continuous part
$\mu_{ac}$ of $\mu$.

There are several ways to characterize the compactness of an operator but, in the case of
an operator represented by a band matrix, a very practical characterization is available:
compactness is equivalent to stating that all the diagonals converge to zero. If the
matrix is not banded the convergence of the diagonals to zero is only a necessary
condition for the compactness. The compactness can be also used to characterize certain
properties of the essential spectrum. For instance, given a unitary operator $T$,
$\sigma_e(T)\subset\{\lambda_1,\dots,\lambda_n\}$ iff $(\lambda_1-T)\cdots(\lambda_n-T)$
is compact (Krein's theorem: see \cite{AkKr62} and \cite{Go00}).

The trace class operators, i.e., the operators $T$ such that $\sqrt{T^\dag T}$ has finite
trace, are more difficult to characterize, even if they are represented by a band matrix.
Nevertheless, any infinite matrix $(k_{i,j})$ that satisfies the condition
$\sum_{i,j}|k_{i,j}|<\infty$ represents a trace class operator on $\ell^2$.

Concerning the compactness and trace class character of
$\cU(\bsa,\bsalpha)-\cU(\bsb,\bsbeta)$, Corollary \ref{IDEAL-U} implies that it is a
consequence of the same property for $\cA(\bsalpha)-\cA(\bsbeta)$ and
$\cC(\bsa)-\cC(\bsb)$. The diagonal matrix $\cA(\bsalpha)-\cA(\bsbeta)$ represents a
compact operator iff $\lim_n(\alpha_n-\beta_n)=0$, and is trace class iff
$\bsalpha-\bsbeta\in\ell^2$. As for the compactness and trace class arguments for
$\cC(\bsa)-\cC(\bsb)$ in the applications that we will discuss, they follow the same
lines as in \cite{Si105}.

As a first group of applications in the study of the dependence $\mu(\bsa,\bsalpha)$, we
will analyze the extreme behaviors corresponding to a sequence $\bsa$ converging to zero
or (subsequently) to the unit circle. In what follows Lim$_nx_n$ means the set of limit
points of a sequence $(x_n)$ in $\C$.

\bt \label{LIMITS}

If $\bsalpha$ is compactly included in $\D$, then:
\be
\item
$\ds \lim_na_n=0 \;\Rightarrow\; \supp\mu(\bsa,\bsalpha)=\T$.
\item
$\ds \lim_n|a_n|=1 \;\Rightarrow\; \{\supp\mu(\bsa,\bsalpha)\}' =
\mathop{\hbox{\rm Lim}} \limits_n \,\tilde\zeta_n(-\overline a_n a_{n+1})$.
\item
$\ds \limsup_n|a_n|=1 \;\Rightarrow\; \mu(\bsa,\bsalpha)$ {\rm singular}.
\ee

\et

\bpr

It is straightforward to check that, for any sequence $\bsalpha$ in $\D$, the ORF
$(\phi_n)_{n\geq0}$ corresponding to the Lebesgue measure
$$
dm(e^{i\theta})=\frac{d\theta}{2\pi}=\frac{1}{2\pi i}\frac{dz}{z},
\qquad z=e^{i\theta},
\qquad \theta\in[0,2\pi),
$$
are given by
$$
\ba{l}
\ds \phi_0=1,
\smallskip \cr
\ds \phi_n=\eta_n\frac{\varpi_0^*}{\varpi_n}B_{n-1}, \qquad n\geq1,
\ea
$$
and satisfy recurrence (\ref{RR}) with parameters $a_n=0$. Therefore, when
$\bsalpha$ is compactly supported in $\D$, the unitary matrix $\cU(0,\bsalpha)$
represents the multiplication operator $T_m$. So,
$\sigma(\cU(0,\bsalpha))=\supp\,m=\T$.

Now, suppose an arbitrary sequence $\bsa$ in $\D$ such that $\lim_na_n=0$. Then
$\cC(\bsa)-\cC(0)$ is compact, thus $\cU(\bsa,\bsalpha)-\cU(0,\bsalpha)$ is compact too.
Hence, Weyl's theorem implies $\{\supp\mu(\bsa,\bsalpha)\}'=\{\supp\, m\}'=\T$, that is,
$\supp\mu(\bsa,\bsalpha)=\T$.

\medskip

If $\lim_n|a_n|=1$, then $\cC(\bsa)-\cD(\bsa)$ is compact, where $\cD(\bsa)$ is the
diagonal matrix
$$
\cD(\bsa)=\pmatrix{-a_1 \cr & -\overline a_1a_2 \cr && -\overline a_2a_3 \cr &&& \ddots}.
$$
Therefore, Proposition \ref{IDEAL} implies that
$\cU(\bsa,\bsalpha)-\tilde\zeta_{\cA(\bsalpha)}(\cD(\bsa))$ is compact too. Notice that
$$
\tilde\zeta_{\cA(\bsalpha)}(\cD(\bsa))=
\pmatrix{\tilde\zeta_0(-a_1) \cr & \tilde\zeta_1(-\overline a_1a_2) \cr
&& \tilde\zeta_2(-\overline a_2a_3) \cr &&& \ddots}
$$
is diagonal and bounded, so it is normal and Weyl's theorem states that
$\{\supp\mu(\bsa,\bsalpha)\}'=\hbox{Lim}_n\tilde\zeta_n(-\overline a_na_{n+1})$.

\medskip

Finally, assume that $\limsup_n|a_n|=1$. This means that there is a subsequence
$(a_n)_{n\in\cI}$, $\cI\subset\N$, such that $\lim_{n\in\cI}a_n=a\in\T$. Without loss of
generality we can suppose $\sum_{n\in\cI}|a_n-a|^{1/2}<\infty$, so that
$\sum_{n\in\cI}(|a_n-a|+\rho_n)<\infty$ because $\rho_n\leq\sqrt{2|a_n-a|}$. Let $\bsb$
be the sequence defined by
$$
b_n = \cases{a & if $n\in\cI$, \cr a_n & if $n\notin\cI$.}
$$
The condition $\sum_{n\in\cI}(|a-a_n|+\rho_n)<\infty$ ensures that
$\cC_o(\bsa)-\cC_o(\bsb)$ and $\cC_e(\bsa)-\cC_e(\bsb)$ are trace class, so the same
holds for $\cU(\bsa,\bsalpha)-\cU(\bsb,\bsalpha)$. The Birman-Krein theorem states that
$\supp\mu_{ac}(\bsa,\bsalpha)=\sigma_{ac}(\cU(\bsb,\bsalpha))$, but the fact that
$b_n\in\T$ for infinitely many values of $n$ implies that $\cU(\bsb,\bsalpha)$ decomposes
as a direct sum of finite matrices, so it has a pure point spectrum and, hence, it has no
absolutely continuous part. Therefore, $\mu_{ac}(\bsa,\bsalpha)=0$.

\epr

We can also obtain general conditions for the invariance of
$\{\supp\mu(\bsa,\bsalpha)\}'$ and $\supp\mu_{ac}(\bsa,\bsalpha)$.

\bt \label{GENERAL}

If $\bsalpha$ is compactly included in $\D$, then:
\be
\item
$\ds \lim_n(\alpha_n-\beta_n)=\lim_n(a_n-b_n)=0 \;\Rightarrow\;
\{\supp\mu(\bsa,\bsalpha)\}'=\{\supp\mu(\bsb,\bsbeta)\}'$.
\item
$\ds \sum_n\left(|\alpha_n-\beta_n|+|a_n-b_n|\right)<\infty \;\Rightarrow\;
\supp\mu_{ac}(\bsa,\bsalpha)=\supp\mu_{ac}(\bsb,\bsbeta)$.
\item
If $b_n=\lambda_na_n$ with $\lambda_n\in\C$, then:
$$
\kern-12pt
\ba{l}
\ds \lim_n|\lambda_n|=\lim_n\lambda_{n+1}\overline\lambda_n=1 \;\Rightarrow\;
\{\supp\mu(\bsa,\bsalpha)\}'=\{\supp\mu(\bsb,\bsalpha)\}',
\smallskip \cr
\ds \sum_n(||\lambda_n|^2-1|+|\lambda_{n+1}\overline\lambda_n-1|)<\infty \;\Rightarrow\;
\supp\mu_{ac}(\bsa,\bsalpha)=\supp\mu_{ac}(\bsb,\bsalpha).
\ea
$$
\item \vskip-10pt
$\ds \beta_n=\alpha_{n+N}, \;b_n=a_{n+N} \;\Rightarrow\;
\cases{
\{\supp\mu(\bsa,\bsalpha)\}'=\{\supp\mu(\bsb,\bsbeta)\}',
\smallskip \cr
\kern3pt
\supp\mu_{ac}(\bsa,\bsalpha)=\supp\mu_{ac}(\bsb,\bsbeta).
}$
\ee

\et

\bpr

First, notice that any of the hypothesis of the theorem ensure that $\bsbeta$ is
compactly included in $\D$ when $\bsalpha$ satisfies the same property. Thus, the
spectral theory that we have developed works for both sequences, $\bsalpha$ and
$\bsbeta$. Concerning the notation, in what follows we write $\rho_n=\sqrt{1-|a_n|^2}$,
as usually, and $\sigma_n=\sqrt{1-|b_n|^2}$.

\medskip

To prove the first item, notice that the inequality
$$
|\rho_n-\sigma_n|^2\leq|\rho_n^2-\sigma_n^2|=||a_n|^2-|b_n|^2|\leq2|a_n-b_n|
$$
implies that the conditions $\lim_n(\alpha_n-\beta_n)=\lim_n(a_n-b_n)=0$ ensure the
compactness of $\cA(\bsalpha)-\cA(\bsbeta)$, $\cC_o(\bsa)-\cC_o(\bsb)$ and
$\cC_e(\bsa)-\cC_e(\bsb)$. In consequence, $\cU(\bsa,\bsalpha)-\cU(\bsb,\bsbeta)$ is
compact too and Weyl's theorem implies the equality
$\{\supp\mu(\bsa,\bsalpha)\}'=\{\supp\mu(\bsb,\bsbeta)\}'$.

\medskip

Suppose now $\sum_n\left(|\alpha_n-\beta_n|+|a_n-b_n|\right)<\infty$. If
$\limsup_n|a_n|=1$, then $\limsup_n|b_n|=1$, so we conclude from Theorem \ref{LIMITS}.3
that $\supp\mu_{ac}(\bsa,\bsalpha)=\supp\mu_{ac}(\bsb,\bsbeta)=\emptyset$. If, on the
contrary, $\limsup_n|a_n|<1$, then $|a_n|,|b_n| \leq r$ for some $r<1$. Taking into
account the inequality
$$
|\rho_n-\sigma_n| \leq \frac{||a_n|^2-|b_n|^2|}{\rho_n+\sigma_n} \leq
\frac{r}{\sqrt{1-r^2}} |a_n-b_n|,
$$
$\sum_n\left(|\alpha_n-\beta_n|+|a_n-b_n|\right)<\infty$ implies that
$\cA(\bsalpha)-\cA(\bsbeta)$, $\cC_o(\bsa)-\cC_o(\bsb)$ and $\cC_e(\bsa)-\cC_e(\bsb)$ are
trace class. Thus, $\cU(\bsa,\bsalpha)-\cU(\bsb,\bsbeta)$ is trace class too and, from
the Birman-Krein theorem, $\supp\mu_{ac}(\bsa,\bsalpha)=\supp\mu_{ac}(\bsb,\bsbeta)$.

\medskip

Consider $b_n=\lambda_na_n$ with
$\lim_n|\lambda_n|=\lim_n\lambda_{n+1}\overline\lambda_n=1$. We can write
$\lambda_n=|\lambda_n|e^{i\theta_n}$ with
$\theta_n\in[\theta_{n-1}-\pi,\theta_{n-1}+\pi)$, so that
$\lim_n|\theta_{n+1}-\theta_n|=0$. Define
$$
U=\pmatrix{u_1 \cr & \overline u_2 \cr & & u_3 \cr & & & \overline u_4 \cr & & & &\ddots},
\qquad u_n=e^{i\theta_n/2}.
$$
The identity
$$
\ba{l}
\pmatrix{u_n \cr & \overline u_{n+1}}
\Theta_n(\bsa)
\pmatrix{u_n \cr & \overline u_{n+1}}
- \Theta_n(\bsb) =
\medskip \cr \kern150pt
= \pmatrix{a_nu_n^2(|\lambda_n|^2-1) &
\rho_n\overline u_{n+1}u_n - \sigma_n
\cr
\rho_n\overline u_{n+1}u_n - \sigma_n &
\overline a_n\overline u_n^2 (\overline u_{n+1}^2u_n^2 - |\lambda_n|^2)},
\ea
$$
together with $\lim_nu_{n+1}\overline u_n=1$ and
$|\rho_n-\sigma_n|^2\leq|\rho_n^2-\sigma_n^2|=|1-|\lambda_n|^2||a_n|^2$, shows that
$U\cC_o(\bsa)\kern1pt U-\cC_o(\bsb)$ and $U^\dag\cC_e(\bsa)\kern1pt U^\dag-\cC_e(\bsb)$
are compact. This implies the compactness of
$
U\cC(\bsa)\kern1pt U^\dag\!-\cC(\bsb)=
U\cC_o(\bsa)\kern1pt U \, U^\dag\cC_e(\bsa)\kern1pt U^\dag - \cC_o(\bsb)\,\cC_e(\bsb),
$
which, bearing in mind Proposition \ref{IDEAL}, is equivalent to the compactness of
$
U\cU(\bsa,\bsalpha)\kern1pt U^\dag-\cU(\bsb,\bsalpha)=
\zeta_{\cA(\bsalpha)}(U\cC(\bsa)\kern1pt U^\dag)-\zeta_{\cA(\bsalpha)}(\cC(\bsb)).
$
Therefore, Weyl's Theorem implies that
$\{\supp\mu(\bsa,\bsalpha)\}'=\{\supp\mu(\bsb,\bsalpha)\}'$.

When $\sum_n(||\lambda_n|^2-1|+|\lambda_{n+1}\overline\lambda_n-1|)<\infty$ we have to
consider again two possibilities. If $\limsup_n|a_n|=1$, necessarily $\limsup_n|b_n|=1$
and $\supp\mu_{ac}(\bsa,\bsalpha)=\supp\mu_{ac}(\bsb,\bsalpha)=\emptyset$ from Theorem
\ref{LIMITS}.3. If $\limsup_n|a_n|<1$, then $|a_n|,|b_n| \leq r$ for some $r<1$, so the
relations
$$
\ba{l}
\ds |\rho_n-\sigma_n| \leq \frac{||a_n|^2-|b_n|^2|}{\rho_n+\sigma_n} \leq
\frac{r^2}{2\sqrt{1-r^2}} |1-|\lambda_n|^2|,
\smallskip \cr
\ds |u_{n+1}^2\overline u_n^2 - 1| \sim
|\lambda_{n+1}\overline\lambda_n-|\lambda_{n+1}\overline\lambda_n||
\leq 2|\lambda_{n+1}\overline\lambda_n-1|,
\smallskip \cr
\ds |u_{n+1}\overline u_n - 1| = \frac{|u_{n+1}^2\overline u_n^2 - 1|}{|u_{n+1}\overline u_n + 1|}
\sim \frac{1}{2} |u_{n+1}^2\overline u_n^2 - 1| \leq |\lambda_{n+1}\overline\lambda_n-1|,
\ea
$$
ensure that $U\cC_o(\bsa)\kern1pt U-\cC_o(\bsb)$ and $U^\dag\cC_e(\bsa)\kern1pt
U^\dag-\cC_e(\bsb)$ are trace class. The Birman-Krein theorem then proves that
$\supp\mu_{ac}(\bsa,\bsalpha)=\supp\mu_{ac}(\bsb,\bsalpha)$ similarly to the previous
case.

\medskip

Finally, let $b_n=a_{n+N}$ and $\beta_n=\alpha_{n+N}$ for some $N\in\N$.
Consider the sequences $\tilde\bsa$ and $\tilde\bsalpha$ given by
$$
\tilde a_n = \cases{1 & if $n \leq N$, \cr a_n & if $n>N$,}
\qquad
\tilde \alpha_n = \cases{0 & if $n \leq N$, \cr \alpha_n & if $n>N$.}
$$
$\cA(\bsalpha)-\cA(\tilde\bsalpha)$, $\cC_o(\bsa)-\cC_o(\tilde\bsa)$ and
$\cC_e(\bsa)-\cC_e(\tilde\bsa)$ are finite rank, therefore
$\cU(\bsa,\bsalpha)-\cU(\tilde\bsa,\tilde\bsalpha)$ is compact and trace class. Besides,
we have the decomposition $\cU(\tilde\bsa,\tilde\bsalpha)=-I_N\oplus\cU(\bsb,\bsbeta)$, so
$\cU(\tilde\bsa,\tilde\bsalpha)$ and $\cU(\bsb,\bsbeta)$ have the same essential and
absolutely continuous spectrum. As a consequence of these facts, the Weyl and
Birman-Krein theorems give $\{\supp\mu(\bsa,\bsalpha)\}'=\{\supp\mu(\bsb,\bsbeta)\}'$ and
$\supp\mu_{ac}(\bsa,\bsalpha)=\supp\mu_{ac}(\bsb,\bsbeta)$.

\epr

Combining the different results of the previous theorem we can obtain a more general one.

\bt \label{GENERAL2}

For any sequence $\bsalpha$ compactly included in $\D$:
\be
\item \kern-5pt
If
$\ds\lim_n(\alpha_{n+N}-\beta_n)=\lim_n(\lambda_na_{n+N}-b_n)=0$,
$\ds\lim_n|\lambda_n|=\lim_n\lambda_{n+1}\overline\lambda_n=1$,
then
$\{\supp\mu(\bsa,\bsalpha)\}'=\{\supp\mu(\bsb,\bsbeta)\}'$.
\item \kern-5pt
If
$\ds\sum_n(|\alpha_{n+N}-\beta_n|+|\lambda_na_{n+N}-b_n|+
||\lambda_n|^2-1|+|\lambda_{n+1}\overline\lambda_n-1|)<\infty$,
then
$\supp\mu_{ac}(\bsa,\bsalpha)=\supp\mu_{ac}(\bsb,\bsbeta)$.
\ee

\et

A particular case of this theorem is worthwhile to be emphasized.

\bc \label{LOPEZ}

Let $\alpha\in\D$, $a\in[0,1]$, $\lambda\in\T$ and
$$
\Gamma_{\lambda,a}=\{\lambda e^{i\theta}:|\theta|<2\arcsin a\}.
$$
\be
\item
If $\ds\lim_n\alpha_n=\alpha$, $\ds\lim_n|a_n|=a$ and $\ds\lim_n\frac{a_{n+1}}{a_n}=\lambda$, then
$$
\{\supp\mu(\bsa,\bsalpha)\}'= \T\setminus\tilde\zeta_\alpha(\Gamma_{\lambda,a}).
$$
\item
If $\ds\sum_n\bigg(|\alpha_n-\alpha|+||a_n|-a|+\bigg|\frac{a_{n+1}}{a_n}-\lambda\bigg|\bigg)<\infty$, then
$$
\supp\mu_{ac}(\bsa,\bsalpha) = \T\setminus\tilde\zeta_\alpha(\Gamma_{\lambda,a}).
$$
\ee

\ec

\bpr

Let us write $a_n=|a_n|v_n$ with $v_n\in\T$. Notice that $\bsalpha$ is compactly included
in $\D$ because it is convergent in $\D$. Therefore, we can apply Theorem \ref{GENERAL2}
to $\mu(\bsa,\bsalpha)$ and $\mu(\bsb,\bsbeta)$ with $\beta_n=\alpha$, $b_n=\lambda^na$
and $\lambda_n=\lambda^n\overline v_n$. Taking into account the relation
$$
\ba{l}
\ds |\lambda_{n+1}\overline\lambda_n-1|=\left|\lambda-\frac{v_{n+1}}{v_n}\right|\leq
\left|\lambda-\frac{a_{n+1}}{a_n}\right|+\left|\frac{a_{n+1}}{a_n}-\frac{v_{n+1}}{v_n}\right|=
\medskip \cr \kern149pt
\ds =\left|\lambda-\frac{a_{n+1}}{a_n}\right|+\left|\frac{|a_{n+1}|}{|a_n|}-1\right|\leq
2\left|\frac{a_{n+1}}{a_n}-\lambda\right|,
\ea
$$
we find that $\{\supp\mu(\bsa,\bsalpha)\}'=\{\supp\mu(\bsb,\bsbeta)\}'$ under the
assumptions of item 1, and $\supp\mu_{ac}(\bsa,\bsalpha)=\supp\mu_{ac}(\bsb,\bsbeta)$
under the hypothesis of item 2. On the other hand, from the comments at the beginning of
Section \ref{ORF-HM}, we know that $\mu(\bsb,\bsbeta)=\nu_\alpha$, where
$\nu=\mu(\bsb,0)$ is the measure on $\T$ whose OP have parameters $\lambda^na$ and
$\nu_\alpha$ is defined by $\nu_\alpha(\Delta)=\nu(\zeta_\alpha(\Delta))$ for any Borel
subset $\Delta$ of $\T$. Therefore,
$\{\supp\nu_\alpha\}'=\tilde\zeta_\alpha(\{\supp\nu\}')$,
$\supp(\nu_\alpha)_{ac}=\tilde\zeta_\alpha(\supp\nu_{ac})$ and the corollary follows from
the well known result $\{\supp\nu\}' = \supp\nu_{ac} = \T\setminus\Gamma_{\lambda,a}$.

\epr

If $a=0$, Corollary \ref{LOPEZ}.1 is a direct consequence of Theorem \ref{LIMITS}.1,
while Corollary \ref{LOPEZ}.2 can be derived from Szeg\H o's Theorem for OP on $\T$:
Theorem \ref{GENERAL}.2 implies that
$\supp\mu_{ac}(\bsa,\bsalpha)=\supp\mu_{ac}(\bsa,\bsbeta)$ for $\beta_n=\alpha$ whenever
$\sum_n|\alpha_n-\alpha|<\infty$. $\mu(\bsa,\bsbeta)=\nu_\alpha$, where now
$\nu=\mu(\bsa,0)$, and the condition $\sum_n|a_n|<\infty$ gives $\supp\nu_{ac}=\T$
because $\nu$ is in the Szeg\H o class $\bsa=\cS_0(\nu)\in\ell^2$. Hence,
$\supp\mu_{ac}(\bsa,\bsalpha)=\supp(\nu_\alpha)_{ac}=\T$. In fact, this reasoning proves
that the equality $\supp\mu_{ac}(\bsa,\bsalpha)=\T$ holds under the more general
condition $\sum_n(|\alpha_n-\alpha|+|a_n|^2)<\infty$.

Corollary \ref{LOPEZ} of Theorem \ref{GENERAL2} can be understood also as an example of
the following general result. It says that, when $\bsalpha$ is convergent in $\D$, the
analysis of $\{\supp\mu(\bsa,\bsalpha)\}'$ and $\supp\mu_{ac}(\bsa,\bsalpha)$ can be
related to the much more known polynomial case, corresponding to $\bsalpha=0$.

\bt \label{LIMIT}

Let $\alpha\in\D$.
\be
\item
$\ds \lim_n\alpha_n=\alpha \;\Rightarrow\;
\{\supp\mu(\bsa,\bsalpha)\}'=\tilde\zeta_\alpha(\{\supp\mu(\bsa,0)\}')$.
\item
$\ds \sum_n|\alpha_n-\alpha|<\infty \;\Rightarrow\;
\supp\mu_{ac}(\bsa,\bsalpha)=\tilde\zeta_\alpha(\supp\mu_{ac}(\bsa,0))$.
\ee

\et

\bpr

Again $\bsalpha$ is compactly included in $\D$ because it is convergent in $\D$. So, if
$\beta_n=\alpha$, Theorem \ref{GENERAL} implies that $\{\supp\mu(\bsa,\bsalpha)\}' =
\{\supp\mu(\bsa,\bsbeta)\}'$ when $\lim_n\alpha_n=\alpha$, and
$\supp\mu_{ac}(\bsa,\bsalpha)=\supp\mu_{ac}(\bsa,\bsbeta)$ when
$\sum_n|\alpha_n-\alpha|<\infty$. On the other hand, $\mu(\bsa,\bsbeta)=\nu_\alpha$ with
$\nu=\mu(\bsa,0)$. As in the proof of Corollary \ref{LOPEZ}, the result follows from the
relation between $\nu$ and $\nu_\alpha$.

\epr

The importance of the above theorem is due to the numerous known results for the relation
between $\mu$ and $\bsa$ in the case of OP on $\T$. Theorem \ref{LIMIT} permits us to
translate some of these results to those ORF on $\T$ whose poles converge in $\E$. For
instance, Corollary \ref{LOPEZ}.1 can be understood as the translation to this kind of
ORF of a result for OP on $\T$ due to Barrios-L\'{o}pez (see \cite{BaLo99}). This result was
generalized later on in \cite{LaSi06} as an improvement of a partial extension appearing in
\cite{TRUNC}. The corresponding translation of this generalization to ORF states that
Corollary \ref{LOPEZ}.1 holds even if we substitute the condition $\lim_n|a_n|=a$ by the
more general one $\liminf_n|a_n|=a$.

All the above results provide only sufficient conditions on the sequences $\bsalpha$ and
$\bsa$ to ensure a certain property for the measure $\mu(\bsa,\bsalpha)$. On the
contrary, Krein's theorem permits us to characterize exactly those measures
$\mu(\bsa,\bsalpha)$ with a fixed finite set $\{\supp\mu(\bsa,\bsalpha)\}'$. The
characterization is in terms of the compactness of a matrix depending on $\bsa$ and
$\bsalpha$. The fact that, contrary to the polynomial case, this matrix is not banded
makes difficult to translate its compactness into equivalent conditions for the sequences
$\bsa$ and $\bsalpha$. Nevertheless, in the case of $\{\supp\mu(\bsa,\bsalpha)\}'$ with
at most two points we can find explicitly such equivalent conditions.

\bt \label{KREIN}

If $\bsalpha$ is compactly included in $\D$ and $\lambda,\lambda_1,\lambda_2\in\T$, then:
\be
\item
$\ds \{\supp\mu(\bsa,\bsalpha)\}'=\{\lambda\}$ iff
$$
\lim_n \tilde\zeta_n(-\overline a_n a_{n+1})=\lambda.
$$
\item
$\ds \{\supp\mu(\bsa,\bsalpha)\}'\subset\{\lambda_1,\lambda_2\}$ iff
$$
\ba{l}
\lim_n\rho_n\rho_{n+1}=0,
\medskip \cr
\ds \lim_n\rho_n\left(\frac{\varpi_n(\lambda_1)}{\varpi_n(\alpha_n)}k_n(\lambda_2)
-\frac{\varpi_{n-1}^*(\lambda_2)}{\varpi_{n-1}(\alpha_{n-1})}k_{n-1}(\lambda_1)\right)=0,
\medskip \cr
\ds \lim_n\left(\overline{k_n(\lambda_1)}k_n(\lambda_2)
+(\rho_n^-)^2\,\overline{\varpi_{n-1}^*(\lambda_1)}\,\varpi_{n-1}^*(\lambda_2)\,+\right.
\cr \kern92pt
\left.+\,(\rho_{n+1}^+)^2\,\overline{\varpi_{n+1}(\lambda_1)}\,\varpi_{n+1}(\lambda_2)\right)=0,
\ea
$$
where $k_n(z)=a_n\varpi_n^*(z)+a_{n+1}\varpi_n(z)$.
\ee

\et

\bpr

We are dealing only with measures $\mu$ with an infinite support on $\T$, thus,
$\supp\mu$ has at least one limit point in $\T$. Hence, from Krein's theorem,
$\{\supp\mu\}'=\{\lambda\}$ iff $\lambda-T_\mu$ is compact, i.e., iff $\lambda-\cU$ is
compact. (\ref{MOB}) yields \beq \label{MOB-C1} \lambda-\cU =
\lambda-\tilde\zeta_\cA(\cC) =
\eta_\cA^{-1}\varpi_\cA(\lambda)\,(\zeta_\cA(\lambda)-\cC)\,\tilde\varpi_\cA(\cC)^{-1}\eta_\cA.
\eeq Bearing in mind that $\eta_\cA$, $\varpi_\cA(\lambda)$ and $\tilde\varpi_\cA(\cC)$
are bounded with bounded inverse, the above expression shows that the compactness of
$\lambda-\cU$ is equivalent to the compactness of $\zeta_\cA(\lambda)-\cC$. On the other
hand, $\zeta_\cA(\lambda)-\cC$ is compact iff $\lim_n\rho_n=0$ and
$\lim_n(\zeta_n(\lambda)+\overline a_na_{n+1})=0$. However, the first of these conditions
is a consequence of the second one because $|\zeta_n(\lambda)+\overline a_na_{n+1}| \geq
1-|a_n|$ since $\lambda\in\T$. Also, taking into account (\ref{MOB}),
$$
\varpi_n(\lambda)\,(\zeta_n(\lambda)+\overline a_na_{n+1}) =
(\lambda-\tilde\zeta_n(-\overline a_na_{n+1}))\,\tilde\varpi_n(-\overline a_na_{n+1}).
$$
Therefore, $\lim_n(\zeta_n(\lambda)+\overline a_na_{n+1})=0$ iff
$\lim_n(\lambda-\tilde\zeta_n(-\overline a_na_{n+1}))=0$ because
$2>|\varpi_n(\lambda)|,|\tilde\varpi_n(-\overline a_na_{n+1})|\geq1-|\alpha_n|$ and
$\bsalpha$ is compactly supported in $\D$.

As for the case of two limit points, Krein's theorem implies that the inclusion
$\{\supp\mu\}'\subset\{\lambda_1,\lambda_2\}$ is equivalent to the compactness of the matrix
$(\lambda_1-\cU)(\lambda_2-\cU)$. To express this condition as the compactness of a band
matrix we use (\ref{MOB-C1}) for the factor $\lambda_2-\cU$, but for $\lambda_1-\cU$ we
use the equality
\beq \label{MOB-C2} \lambda-\cU = \lambda-\zeta_{-\cA}(\cC) =
\eta_\cA\varpi_{-\cA}(\cC)^{-1}(\zeta_\cA(\lambda)-\cC)\,\tilde\varpi_{-\cA}(\lambda)\,\eta_\cA^{-1},
\eeq
obtained from (\ref{MOB}) and the identity $\tilde\zeta_\cA=\zeta_{-\cA}$. Then,
similarly to the case of one limit point, we find that $(\lambda_1-\cU)(\lambda_2-\cU)$
is compact iff the 9-diagonal matrix
$(\zeta_\cA(\lambda_1)-\cC)\,\varpi_\cA(\lambda_1)\,\varpi_\cA(\cA)^{-1}
\varpi_\cA(\lambda_2)(\zeta_\cA(\lambda_2)-\cC)$ is compact. This compactness condition
can be equivalently formulated using a simpler band matrix obtained multiplying the above
one on the left and the right by the unitary matrices $\cC_o^\dag$ and $\cC_e^\dag$
respectively. Taking into account the identity $\varpi_\cA^*(z)=z\,\varpi_\cA(z)^\dag$,
$z\in\T$, we find in this way that $\{\supp\mu\}'\subset\{\lambda_1,\lambda_2\}$ iff the
five-diagonal matrix $K(\lambda_1)^\dag\varpi_\cA(\cA)^{-1}K(\lambda_2)$ is compact,
where $K(z)=\varpi_\cA^*(z)\,\cC_e^\dag-\varpi_\cA(z)\,\cC_o$. Now, it is just a matter
of calculating the diagonals of $K(\lambda_1)^\dag\varpi_\cA(\cA)^{-1}K(\lambda_2)$
to obtain the conditions given in the theorem.

\epr

The implication $\lim_n\tilde\zeta_n(-\overline a_na_{n+1})=\lambda\in\T \Rightarrow
\{\supp\mu\}'=\{\lambda\}$ was in fact a consequence of Theorem \ref{LIMITS}.2. Krein's
theorem adds the opposite implication. Concerning the case of two limit points notice
that, although the third condition is symmetric under the exchange of $\lambda_1$ and
$\lambda_2$, the second one does not show explicitly such a symmetry. However, a detailed
analysis of the second condition reveals that it is symmetric too.

It seems that there is no simple way to generalize the arguments given in the proof of
Theorem \ref{KREIN} to the case of more than two limit points. The reason is that, for
$n\geq3$, identities (\ref{MOB-C1}) and (\ref{MOB-C2}) are not enough to reduce the
compactness of $(\lambda_1-\cU)\cdots(\lambda_n-\cU)$ to the compactness of a band
matrix. So, contrary to the polynomial situation (see \cite{Go00} and \cite{MIN,LaSi06}),
the practical application of Krein's theorem to characterize in terms of the sequences
$\bsa$ and $\bsalpha$ those measures on $\T$ whose support has a finite set of more than
two limit points remains as an open problem in the rational case.

\section{Appendix: ORF on the real line} \label{ORF-R}

In what follows, a measure on the real line will be probability Borel measure $\mu$
supported on an infinite subset $\supp\mu$ of $\overline\R$. When $\infty$ is not a mass
point of $\mu$ we will refer to $\mu$ as a measure on $\R$. Notice that we are
considering all these measures as measures on $\overline\R$, no matter whether they have
a mass point at $\infty$ or not. This means that $\infty\in\supp\mu$ when $\infty$ is a
mass point of $\mu$ or when $\mu$ is a measure on $\R$ with unbounded standard support,
so that $\supp\mu$ is always closed in $\overline\R$.

Analogously to the case of the unit circle, for any measure $\mu$ on the real line it is
possible to consider ORF in $L^2_\mu$ with poles in the lower half plane $\LL=\{z\in\C :
\im(z)<0\}$. For this purpose we introduce for any $\alpha\in\U=\{z\in\C :
\im(z)>0\}$ the linear fractional transformation
$$
\zeta_\alpha(z) = \frac{\varpi_\alpha^*(z)}{\varpi_\alpha(z)}, \qquad
\cases{
\varpi_\alpha(z) =
z - \overline\alpha,
\medskip \cr
\varpi_\alpha^*(z) = z - \alpha,
}
$$
which maps $\overline\R$, $\U$ and $\LL$ onto $\T$, $\D$ and $\E$ respectively, and has
the inverse
$$
\tilde\zeta_\alpha(z) = \frac{\tilde\varpi_\alpha^*(z)}{\tilde\varpi_\alpha(z)}, \qquad
\cases{
\tilde\varpi_\alpha(z) = 1 - z,
\medskip \cr
\tilde\varpi_\alpha^*(z) = \alpha - \overline\alpha z.
}
$$
Notice that $\varpi_{\alpha}^*=\varpi_{\alpha*}$, where the $*$-involution is now defined
by $f_*(z)=\overline{f(\overline z)}$, but nothing similar holds for
$\tilde\varpi_\alpha^*$. Besides, for the distinguished value $\alpha_0=i$,
$\zeta=\zeta_{\alpha_0}$ is the Cayley transform and $\tilde\zeta=\tilde\zeta_{\alpha_0}$
its inverse.

Any sequence $\bsalpha=(\alpha_n)_{n\geq1}$ in $\U$ defines the products $(B_n)_{n\geq0}$
as in (\ref{BLAS}), but with the new meaning for $\zeta_{\alpha_n}$. The
orthonormalization in $L^2_\mu$ of $(B_n)_{n\geq0}$ leads to a sequence
$(\phi_n)_{n\geq0}$ of ORF with respect to $\mu$ with poles in
$(\overline\alpha_n)_{n\geq1}$, which will be called a sequence of ORF on the real line.
The study of ORF on the real line can be carried out in a completely analogous way to the
case of the unit circle, so most of the results described for the last ones translate
directly to the first ones with an obvious change of the meaning in the notations. In
particular, the sequence $(\phi_n)_{n\geq0}$ can be chosen such that it satisfies a
recurrence like (\ref{RR}) depending on a sequence $\bsa=(a_n)_{n\geq1}$ in $\D$, which
establishes a surjective application $\cS_\bsalpha\colon\frP\to\D^\infty$, where $\frP$
means now the set of probability measures on $\overline\R$. This application is a
bijection when $B_n(z)=\prod_{n=1}^\infty\zeta_n(z)$ diverges to zero for $z\in\U$, but
this is equivalent now to $\sum_{n=1}^\infty\im\,\alpha_n/(1+|\alpha_n|^2)=\infty$, which
means that the poles can not approach too quickly to $\overline\R$.

\medskip

Following the same strategy as in the case of the unit circle, we can develop a spectral
theory for ORF on the real line. The starting point is again recurrence (\ref{RR})
written in the form (\ref{RR1}), but now the positive factors $\eta_\alpha$,
$\alpha\in\U$, are defined by
$$
\eta_\alpha=\left(\frac{\varpi_\alpha(\alpha)}{2i}\right)^{1/2}=\sqrt{\im\,\alpha}.
$$
Both, the expressions for the unit circle and the real line can be combined in
$\eta_\alpha=(\varpi_\alpha(\alpha)/\varpi_{\alpha_0}(\alpha_0))^{1/2}$.

The form (\ref{RR1}) of the recurrence is the key tool to obtain the matrix
representations with respect to the ORF for the multiplication operator
$$
T_\mu \colon \mathop{L^2_\mu \to L^2_\mu} \limits_{f(z) \; \to \; zf(z)}
$$
where $\mu$ is the corresponding orthogonality measure on the real line. If $\supp\mu$ is
bounded, $T_\mu$ is an everywhere defined self-adjoint operator on $L^2_\mu$. In general,
$T_\mu$ is a densely defined self-adjoint operator on $L^2_\mu$ when the function $z$ is
finite $\mu$-a.e. (see \cite[page 259]{ReSi72}), that is, when $\infty$ is not a mass
point of $\mu$. In this case, $\sigma_p(T_\mu)=\{\hbox{mass points of } \mu\}$ and
$\sigma(T_\mu)=\supp\mu$ under the convention that $\infty\in\sigma(T_\mu)$ when $T_\mu$
has an unbounded standard spectrum. A way to deal with the case of measures with a mass
point at $\infty$ is to work with the operator multiplication by $\zeta$ in $L^2_\mu$,
i.e.,
$$
S_\mu \colon \mathop{L^2_\mu \longrightarrow L^2_\mu}
\limits_{f(z) \; \to \; \zeta(z)f(z)}
$$
This operator is unitary for any measure $\mu$ on $\overline\R$ and verifies the
identities $\sigma_p(S_\mu)=\zeta(\hbox{mass points of } \mu)$ and
$\sigma(S_\mu)=\zeta(\supp\mu)$. The matrix representations of $T_\mu$ and $S_\mu$ with
respect to the related ORF are related to the operator analogs of the new linear
fractional transformations $\zeta_\alpha$.

\medskip

To discuss such operator linear fractional transformations it is convenient to introduce
the notation
$$
\re\,T=\frac{1}{2}(T+T^\dag), \qquad \im\,T=\frac{1}{2i}(T-T^\dag),
$$
for any densely defined operator $T$ on $H$. The operator linear fractional
transformations of interest for ORF on the real line are
$$
\ba{l} \zeta_A(T) = \eta_A \,\varpi_A(T)^{-1} \varpi_A^*(T) \,\eta_A^{-1},
\kern10pt
\cases{ \varpi_A(T) = T - A^\dag,
\medskip \cr
\varpi_A^*(T) = T - A,
}
\medskip \cr
\tilde\zeta_A(T) = \eta_A^{-1} \tilde\varpi_A^*(T) \,\tilde\varpi_A(T)^{-1} \eta_A,
\kern12pt
\cases{ \tilde\varpi_A(T) = 1 - T,
\medskip \cr
\tilde\varpi_A^*(T) = \eta_A A \kern1pt \eta_A^{-1} \! - \eta_A A^\dag \eta_A^{-1} T,
}
\ea
$$
where
$$
\eta_A = \sqrt{\im\,A}
$$
and $A\in\B_H$ is such that $\im\,A\geq\varepsilon$ for some positive number
$\varepsilon$ (in short, $\im\,A>0$), so that $\eta_A$ is bounded with bounded inverse.
When $A$ is normal, as it is the case related to ORF on the real line,
$\tilde\varpi_A^*$ becomes
$$
\tilde\varpi_A^*(T) = A - A^\dag T.
$$

$\zeta_A$ is a bijection of $\U_H=\{T\in\B_H:\im\,T>0\}$ onto $\D_H$, and $\tilde\zeta_A$
is its inverse. To prove this assertion we start showing that $\zeta_A$ and
$\tilde\zeta_A$ map $\U_H$ and $\D_H$ respectively on $\B_H$. The statement for
$\tilde\zeta_A$ is a consequence of the fact that $\|T\|<1$ implies
$\tilde\varpi_A(T)^{-1}\in\B_H$. As for $\zeta_A$, the result follows from the fact that
the spectrum of any operator $T\in\B_H$ is included in the closure of its numerical range
$\{(x,Tx):\|x\|=1\}$. So, if $\im\,T>0$, then $\sigma(T)\subset\U$ and thus
$T^{-1}\in\B_H$. In consequence, $\varpi_A(T)^{-1}\in\B_H$ for any $T\in\U_H$ since
$\im\,(T-A^\dag)\geq\im\,A>0$.

On the other hand, using the equality
$A^\dag\eta_A^{-2}A=2iA+A\,\eta_A^{-2}A=A\,\eta_A^{-2}A^\dag$, we find the identities
\beq
\label{UB} \ba{l}
\varpi_A(T)\,\eta_A^{-1}(1-\zeta_A(T)\,\zeta_A(T)^\dag)\,\eta_A^{-1}\varpi_A(T)^\dag =
4\,\im\,T,
\smallskip \cr
\tilde\varpi_A(T)^\dag\eta_A^{-1}(\im\,\tilde\zeta_A(T))\,\eta_A^{-1}\tilde\varpi_A(T) =
1-T^\dag T, \ea
\eeq
which prove that $\zeta_A$ maps $\U_H$ on $\D_H$ and $\tilde\zeta_A$ does the opposite.
Moreover, a direct calculation shows that, for any $T\in\U_H$ and any $S\in\D_H$,
$S=\zeta_A(T)$ iff $T=\tilde\zeta_A(S)$. This completes the proof.

The above arguments can be easily generalized to see that
$\zeta_A$ extends to a transformation of $\overline\U_H=\{T\in\B_H:\im\,T\geq0\}$
onto $\{T\in\overline\D_H:1\notin\sigma(T)\}$, $\tilde\zeta_A$ being its inverse. In
consequence, $\zeta_A$ maps $\overline\U_H\setminus\U_H$ onto
$\{T\in\T_H:1\notin\sigma(T)\}$ and $\tilde\zeta_A$ does the converse.
Furthermore, (\ref{UB}) also implies that $\zeta_A$ maps the set of bounded self-adjoint
operators onto the set of unitary operators whose spectrum does not contain 1.

The above properties are verified in particular by the Cayley transform for operators,
since it is given by $\zeta=\zeta_{A=i}$. Indeed, $\zeta_A$ is nothing but the
composition of the Cayley transform with an operator transformation depending on $A$
which maps onto theirselves $\overline\U_H$, $\U_H$ and the set of self-adjoint
operators on $H$. More precisely, taking into account that
$$
\ba{l}
\eta_A^{-1}\varpi_A(T)\,\eta_A^{-1} = \eta_A^{-1}(T-\re\,A)\,\eta_A^{-1}+i,
\smallskip \cr
\eta_A^{-1}\varpi_A^*(T)\,\eta_A^{-1} = \eta_A^{-1}(T-\re\,A)\,\eta_A^{-1}-i,
\ea
$$
we obtain \beq \label{CAYLEY} \zeta_A(T)=\zeta(\eta_A^{-1}(T-\re\,A)\,\eta_A^{-1}). \eeq
It is known that the Cayley transform extends to a bijection between the set of (bounded
or unbounded) self-adjoint operators and the set of unitary operators whose point
spectrum does not contain 1, so the same holds for $\zeta_A$. The importance of this
property is that it permits us to formulate the spectral theory for ORF on the real line,
including the case of measures on $\R$ with unbounded support since they are associated
with unbounded self-adjoint multiplication operators.

Another advantage of relation (\ref{CAYLEY}) is that it expresses $\zeta_A$ as a product
of two commutative factors. This provides two equivalent representations of $\zeta_A$,
namely,
$$
\zeta_A(T) = \eta_A (T-A^\dag)^{-1} (T-A) \,\eta_A^{-1} =
\eta_A^{-1} (T-A) (T-A^\dag)^{-1} \eta_A,
$$
giving rise to two expressions for $\tilde\zeta_A$ too. From the above
result we find that $\zeta_A(T)^\dag=\zeta_{A^\dag}(T^\dag)$ and
$\tilde\zeta_A(T)^\dag=\tilde\zeta_{A^\dag}(T^\dag)$, as in the case of the unit circle.

Finally, if $\frI$ is a hermitian ideal of $\B_H$, similar arguments to those given in
the proof of Theorem \ref{IDEAL} prove that, for any normal operators $A,B\in\U_H$
such that $AB=BA$, the condition $A-B\in\frI$ implies the equivalences
$$
\ba{l}
T-S\in\frI \;\Leftrightarrow\; \zeta_A(T)-\zeta_B(S)\in\frI, \qquad
\forall\,T,S\in\overline\U_H,
\smallskip \cr
T-S\in\frI \;\Leftrightarrow\; \tilde\zeta_A(T)-\tilde\zeta_B(S)\in\frI, \qquad
\forall\,T,S\in\overline\D_H, \quad 1\notin\sigma(T)\cup\sigma(S).
\ea
$$
The right implication of each case is equivalent to the left implication of the other one
due to the fact that $\zeta_A$ and $\tilde\zeta_A$ are mutually inverse transformations.
As we pointed out, when $T$ is unitary the transformation $\zeta_A(T)$ is well defined
provided that 1 is not an eigenvalue of $T$. So, the right implication of the second
equivalence can be formulated in a more general context when the operators $T,S$ are
unitary. The analogous extension for the right implication of the first equivalence,
i.e., the case of $T,S$ unbounded self-adjoint operators, is not possible because we
suppose that $\frI$ is an ideal of $\B_H$ (as it is the case for the classes of
perturbations usually considered in operator theory). Therefore, we only can assure that
$$
\ba{l}
T-S\in\frI \;\Leftarrow\; \zeta_A(T)-\zeta_B(S)\in\frI, \qquad
\forall\,T,S \hbox{ self-adjoint},
\smallskip \cr
T-S\in\frI \;\Rightarrow\; \tilde\zeta_A(T)-\tilde\zeta_B(S)\in\frI, \qquad
\forall\,T,S \hbox{ unitary}, \quad 1\notin\sigma_p(T)\cup\sigma_p(S).
\ea
$$
These results, although weaker than the previous ones, are enough to apply perturbative
techniques to the spectral theory of ORF on the real line, even if the support of the
orthogonality measure is unbounded.

\medskip

With all these operator tools at hand we can develop the spectral theory for ORF on the
real line following the same steps as in the case of the unit circle. In fact, the
results for the unit circle are formulated throughout the paper in such a way that the
translation to the real line is just a matter of changing the meaning of the symbols
according to the previous discussion, together with some other obvious modifications. For
instance, if $\mu$ is a measure on $\overline\R$ and $(\phi_n)_{n\geq0}$ are the ORF
associated with an arbitrary sequence $\bsalpha=(\alpha_n)_{n\geq1}$ in $\U$, then the
sequence $(\chi_n)_{n\geq0}$ defined by (\ref{LORF}), with the new meaning for
$\zeta_n=\zeta_{\alpha_n}$, are the ORF associated with
$(\alpha_1,\hat\alpha_2,\alpha_3,\hat\alpha_4,\dots)$. $(\chi_n)_{n\geq0}$ is a basis of
$L^2_\mu$ when the odd and even products $B^o(z)$ and $B^e(z)$ converge to zero for
$z\in\U$, but this means now that
$\sum_{k=1}^\infty\im\,\alpha_{2k-1}/(1+|\alpha_{2k-1}|^2)=
\sum_{k=1}^\infty\im\,\alpha_{2k}/(1+|\alpha_{2k}|^2)=\infty$. Also, if
$\bsa=(a_n)_{n\geq1}$ are the parameters of the recurrence for $\phi_n$, the zeros of
$\phi_n$ are the eigenvalues of $\cU^{(n)}=\tilde\zeta_{\cA_n}(\cC_n)$, where
$\cC=\cC(\bsa)$ and $\cA=\cA(\bsalpha)$ as in the case of the unit circle, but
$\tilde\zeta_{\cA_n}$ is the new operator linear fractional transformation given in this
section. Notice that $\cA_n\in\U_{\C^n}$ and $\cC_n\in\T_{\C^n}$ has its eigenvalues in
$\D$ because they are the zeros of the $n$-th OP associated with the parameters $\bsa$.
Hence, $\cU^{(n)}$ is a well defined matrix of $\overline\U_{\C^n}\setminus\U_{\C^n}$,
which agrees with the fact that the zeros of $\phi_n$ lie on $\U$.

Other results for the unit circle can be translated to the real line in a similar way,
but two of the main results need a special discussion. The first one concerns the
representation of the self-adjoint multiplication operator $T_\mu$ for a measure $\mu$ on
$\R$, and the other one is related to the representation of the self-adjoint
multiplication operator $T_{\mu_n^v}$ corresponding to the finitely supported measure
$\mu_n^v$ associated with the PORF $Q_n^v$.

Following the same steps as in Theorem \ref{BANDORF}, we would find that, if $\mu$ is a
measure on $\R$, for any sequence $\bsalpha$ compactly included in $\U$, the matrix
representation of $T_\mu$ with respect to the ORF $(\chi_n)_{n\geq0}$ associated with
$(\alpha_1,\hat\alpha_2,\alpha_3,\hat\alpha_4,\dots)$ is $\cU=\tilde\zeta_\cA(\cC)$,
where $\cA=\cA(\bsalpha)$, $\cC=\cC(\bsa)$ and $\bsa=\cS_\bsalpha(\mu)$. However, since
the matrix $\cC$ is unitary, we can assure that $\tilde\zeta_\cA(\cC)$ provides a well
defined (self-adjoint) operator only when 1 is not an eigenvalue of $\cC$. That is, in
the case of the real line, the matrix representation $\cU=\tilde\zeta_A(\cC)$ is valid
provided that $1\notin\sigma_p(\cC)$. To understand the meaning of this condition we will
relate $\cC$ to the matrix representation with respect to $(\chi_n)_{n\geq0}$ of $S_\mu$.
When $1\notin\sigma_p(\cC)$ the matrix of $S_\mu=\zeta(T_\mu)$ is $\zeta(\cU)$, but, as
we will see, an expression for the matrix representation of $S_\mu$ can be obtained for
any measure $\mu$ on the real line, even if it has a mass point at $\infty$. This
discussion will lead also to a relation between the operator linear fractional
transformations in the real line and the unit circle.

\medskip

Since we are going to consider at the same time the linear fractional transformations
used on the real line and on the unit circle, in what follows we will distinguish between
both cases with a superscript ${}^\R$ or ${}^\T$ respectively. Let $A\in\U_H$. Due to the
properties of the Cayley transform, $B=\zeta(A)\in\D_H$. A direct computation gives
$$
\im\,A  = (1-B)^{-1}(1-BB^\dag)\,(1-B^\dag)^{-1}.
$$
Therefore, $\eta_A^\R=|\eta_B^\T(1-B^\dag)^{-1}|$ and, as a consequence of the polar
decomposition,
\beq \label{ETA-U}
\eta_B^\T(1-B^\dag)^{-1}=U\eta_A^\R, \qquad U \hbox{ unitary}.
\eeq
If we change $A$ by $-A^\dag$, then $B$ changes to $B^\dag$, thus,
\beq \label{ETA-V}
\eta_{B^\dag}^\T(1-B)^{-1}=V\eta_A^\R, \qquad V \hbox{ unitary}.
\eeq
When $A$ is normal, $B$ is normal too and $\eta_A^\R=|1-B|^{-1}\eta_B^\T$, so $U=V^\dag=\xi_B$,
where
$$
\xi_B=(1-B)|1-B|^{-1}.
$$
In the general case, using (\ref{ETA-U}) and (\ref{ETA-V}), we find that
$$
\zeta_B^\T(\zeta(T))=U\zeta_A^\R(T)V^\dag,
$$
hence
\beq \label{ZR-ZT}
\zeta(\tilde\zeta_A^\R(T))=\tilde\zeta_B^\T(UTV^\dag).
\eeq
Denoting $w=\zeta(z)$ and $S=UTV^\dag$, a straightforward calculation gives
\beq \label{R-T-1}
\varpi_A^{*\R}(z) - \varpi_A^\R(z) \, T_A^\R = \frac{2i}{1-w}
\bigg(\varpi_B^{*\T}(w) - \varpi_B^\T(w) \, S_B^\T\bigg) (1-B)^{-1},
\eeq
where $T_A^\R=(\eta_A^\R)^{-1} T \, \eta_A^\R$ and $S_B^\T=(\eta_B^\T)^{-1} S \,
\eta_B^\T$. Since equations (\ref{MOB2}) and (\ref{MOB1}) hold for the real line too, the
above equality can be written equivalently as
\beq \label{R-T-2}
z\,\tilde\varpi_A^\R(T_A^\R) - \tilde\varpi_A^{*\R}(T_A^\R) = \frac{2i}{1-w}
\bigg(w\,\tilde\varpi_B^\T(S_B^\T) - \tilde\varpi_B^{*\T}(S_B^\T)\bigg) (1-B)^{-1}.
\eeq
Using (\ref{ETA-U}) and (\ref{ETA-V}) we obtain $S_B^\T=(1-B^\dag)^{-1}T_A^\R\kern1pt(1-B)$.
Taking this relation into account, a direct computation yields
$$
1-\tilde\zeta_B^\T(S) = 1 - \tilde\varpi_B^{*\T}(S_B^\T)\,\varpi_B^\T(S_B^\T)^{-1} =
(\eta_A^\R)^{-1} (1-T) \kern1pt V^\dag \tilde\varpi_B^\T(S)^{-1} \eta_{B^\dag},
$$
which implies that
\beq \label{R-T-3}
1\in\sigma(\tilde\zeta_B^\T(S)) \;\Leftrightarrow\; 1\in\sigma(T),
\qquad
1\in\sigma_p(\tilde\zeta_B^\T(S)) \;\Leftrightarrow\; 1\in\sigma_p(T).
\eeq

\medskip

Assume now that $\bsalpha$ is compactly included in $\U$ and $\mu$ is a measure on $\R$
such that $1\notin\sigma_p(\cC)$. From (\ref{ZR-ZT}) we see that the matrix
representation $\zeta(\tilde\zeta_\cA^\R(\cC))$ of $S_\mu$ can be expressed alternatively
as $\tilde\zeta_\cB^\T(\xi_\cB\,\cC\,\xi_\cB)$, with $\cB=\zeta(\cA)$. Nevertheless,
contrary to $\zeta(\tilde\zeta_\cA^\R(\cC))$, $\tilde\zeta_\cB^\T(\xi_\cB\,\cC\,\xi_\cB)$
is always a well defined (unitary) matrix, no matter whether 1 is an eigenvalue of $\cC$
or not, because $\xi_\cB\,\cC\,\xi_\cB$ is unitary and $\tilde\zeta_\cB^\T$ maps unitary
operators into unitary operators. Actually, we are going to prove that, if $\bsalpha$
compactly included in $\U$, $\tilde\zeta_\cB^\T(\xi_\cB\,\cC\,\xi_\cB)$ is the matrix
representation of $S_\mu$ with respect to $(\chi_n)_{n\geq0}$ for any measure $\mu$ on
$\overline\R$. Following similar arguments to those given in the proof of Theorem
\ref{BANDORF} we find that, for any measure $\mu$ on $\overline\R$, the ORF
$(\chi_n)_{n\geq0}$ satisfy equation (\ref{RR-ORF}) too, but substituting
$\hat\cC=\cC_\cA^\T$ by $\hat\cC=\cC_\cA^\R$, and $\varpi_\cA^\T$, $\varpi_\cA^{*\T}$ by
$\varpi_\cA^\R$, $\varpi_\cA^{*\R}$ respectively. Applying (\ref{R-T-1}) and using
(\ref{MOB2}) and (\ref{MOB1}) we conclude that, for $\bsalpha$ compactly included in
$\U$,
$$
\pmatrix{\chi_0(z) & \chi_1(z) & \cdots}
\left( \zeta(z) - \tilde\zeta_\cB^\T(\xi_\cB\,\cC\,\xi_\cB) \right) = 0,
\qquad \cB=\zeta(\cA),
$$
which means that $\tilde\zeta_\cB^\T(\xi_\cB\,\cC\,\xi_\cB)$ is the matrix of $S_\mu$
with respect to $(\chi_n)_{n\geq0}$. As a consequence of this result and (\ref{R-T-3}),
we have the equivalences
$$
\kern-5pt
\ba{c}
1\in\sigma(\cC) \,\Leftrightarrow\,
1\in\sigma(\tilde\zeta_\cB^\T(\xi_\cB\,\cC\,\xi_\cB)) \,\Leftrightarrow\,
1\in\sigma(S_\mu) \,\Leftrightarrow\,
1\in\zeta(\supp\mu),
\smallskip \cr
1\in\sigma_p(\cC) \,\Leftrightarrow\,
1\in\sigma_p(\tilde\zeta_\cB^\T(\xi_\cB\,\cC\,\xi_\cB)) \,\Leftrightarrow\,
1\in\sigma_p(S_\mu) \,\Leftrightarrow\,
1\in\zeta(\hbox{mass points of } \mu).
\ea
$$
Thus, we have reached the following result.

\bt \label{1-INF}

Let $\bsalpha$ be a sequence compactly included in $\U$, $\mu$ a measure on
$\overline\R$ and $\cC=\cC(\bsa)$ with $\bsa=\cS_\bsalpha(\mu)$. Then,
$$
1\in\sigma(\cC) \;\Leftrightarrow\; \infty\in\supp\mu, \qquad
1\in\sigma_p(\cC) \;\Leftrightarrow\; \infty \hbox{ is a mass point of } \mu.
$$

\et

Therefore, $\mu$ is a measure on $\R$ iff its related sequence $\bsa$ satisfies
$1\notin\sigma_p(\cC)$. Thus, $\cU=\tilde\zeta_\cA^\R(\cC)$ provides a well defined
matrix representation of $T_\mu$ for any measure $\mu$ on $\R$. Moreover, the measures on
$\R$ with bounded support are characterized by the fact that $\bsa$ is such that
$1\notin\sigma(\cC)$.

In the case of an arbitrary measure $\mu$ on $\overline\R$, including the possibility of
a mass point at $\infty$, we can study the relation $\mu(\bsa,\bsalpha)$ throughout the
spectral analysis of the matrix representation
$\tilde\zeta_\cB^\T(\xi_\cB\,\cC\,\xi_\cB)$ of $S_\mu$ or, alternatively, we can deal
with a pair of operators. More precisely, relation (\ref{R-T-2}) implies that the spectra
of $\tilde\zeta_\cB^\T(\xi_\cB\,\cC\,\xi_\cB)$ and the pair
$(\tilde\varpi_\cA^{*\R}(\cC),\tilde\varpi_\cA^\R(\cC))$ are related by the Cayley
transform, so
$$
\supp\mu = \sigma(\tilde\varpi_\cA^{*\R}(\cC),\tilde\varpi_\cA^\R(\cC)) =
\sigma(\cC_o+\cA\cC_e^\dag,\cC_e^\dag+\cA^\dag\cC_o).
$$
Also, the eigenvalues of the pair are the mass points of $\mu$ and the eigenvectors of
the pair with eigenvalue $\lambda$ are spanned by
$\pmatrix{\chi_0(\lambda)&\cdots&\chi_{n-1}(\lambda)} \eta_{\cA_n}^{-1/2}$. That is,
while the spectral methods that use linear fractional transformations
$\tilde\zeta_\cA^\R$ of five-diagonal matrices only work for measures on $\R$, their
formulation in terms of pairs of band matrices are valid for any measure on
$\overline\R$.

\medskip

Similar results hold too for the finitely supported measures associated with the PORF.
Given an arbitrary measure $\mu$ on $\overline\R$, consider the measure $\mu_n^v$
supported on the zeros of the PORF $Q_n^v=\phi_n+v\phi_n^*$, $v\in\T$. As in the case of
the unit circle, $Q_n^v$ has $n$ different zeros, but now they lie on $\overline\R$.
Besides, if $u=\tilde\zeta_{a_n}$, the matrix representation
$\cU^{(n;u)}=\tilde\zeta_{\cA_n}(\cC_n^u)$ of $T_{\mu_n^v}$ with respect to
$(\chi_k)_{k=0}^{n-1}$ is well defined whenever $1\notin\sigma(\cC_n^u)$. Concerning this
condition, an analogous argument to that of the measure $\mu$ proves that
$$
1\in\sigma(\cC_n^u) \;\Leftrightarrow\; \infty\in\supp\mu_n^v,
$$
i.e., the matrix representation $\cU^{(n;u)}$ of $T_{\mu_n^v}$ is valid for any measure
$\mu_n^v$, except for the value $v=-\phi_n^*(\infty)/\phi_n(\infty)$ which locates a zero
of $Q_n^v$ at $\infty$. Nevertheless, analogously to the previous discussion, the
spectral interpretation of the PORF in terms of pairs of band matrices given for the unit
circle after Theorem \ref{BANDPORFzeros} holds for any PORF on the real line too.

\medskip

Concerning the applications of the spectral theory for ORF on the real line, from the
previous comments we know that, if $\frI$ is an ideal of $\B_{\ell^2}$, for any sequences
$\bsalpha,\bsbeta$ compactly included in $\U$ and any sequences $\bsa,\bsb$ in
$\overline\D$ such that $1\notin\sigma_p(\cC(\bsa))\cup\sigma_p(\cC(\bsb))$,
$$
\cA(\bsalpha)-\cA(\bsbeta),\kern2pt\cC(\bsa)-\cC(\bsb)\in\frI \;\Rightarrow\;
\cU(\bsa,\bsalpha)-\cU(\bsb,\bsbeta)\in\frI.
$$
This permits us to extend to ORF on $\R$ the applications for ORF on $\T$ discussed
in Section \ref{APPL}.

\medskip

Equation (\ref{ZR-ZT}) provides a connection between the real line and the unit circle
representations. Let $\bsalpha=(\alpha_n)_{n\geq1}$ be a sequence compactly included in
$\U$, and consider the sequence $\bsbeta=(\beta_n)_{n\geq1}$ in $\D$ given by
$\beta_n=\zeta(\alpha_n)$. Following the previous notation we also have $\alpha_0=i$, so
$\beta_0=1$. Consider two sequences $\bsa=(a_n)_{n\geq1}$ and $\bsb=(b_n)_{n\geq1}$ in
$\D$ related by
$$
b_n=\xi_0^2\xi_1^2\cdots\xi_{n-1}^2\,a_n, \qquad \xi_n=\frac{1-\beta_n}{|1-\beta_n|}.
$$
We have the identities $\cC_o(\bsb)=\Lambda^\dag\kern1pt\xi_\cB\,\cC_o(\bsa)\,\Gamma$ and
$\cC_e(\bsb)=\Gamma^\dag\kern1pt\cC_e(\bsa)\,\xi_\cB\,\Lambda$, where $\cB=\zeta(\cA)$,
$\cA=\cA(\bsalpha)$ and
\beq \label{GAMMA-LAMBDA}
\ba{l}
\Gamma=\pmatrix{\gamma_0 \cr & \gamma_1 \cr && \ddots}, \quad
\gamma_0=1, \quad
\gamma_n=
\cases{
\overline\xi_0^2\overline\xi_2^2\cdots\overline\xi_{n-1}^2 & odd $n$,
\smallskip \cr
\xi_1^2\xi_3^2\cdots\xi_{n-1}^2 & even $n$,
}
\medskip \cr
\Lambda=\pmatrix{\lambda_0 \cr & \lambda_1 \cr && \ddots}, \quad
\lambda_0=1, \quad
\lambda_n=
\cases{
\gamma_{n-1}\xi_n & odd $n$,
\smallskip \cr
\gamma_{n-1}\overline\xi_n & even $n$,
}
\ea
\eeq
Therefore, $\cC(\bsb)=\Lambda^\dag\kern1pt\xi_\cB\,\cC(\bsa)\,\xi_\cB\,\Lambda$ and,
thus, equation (\ref{ZR-ZT}) implies that
\beq \label{UR-UT}
\zeta(\cU^\R(\bsa,\bsalpha))=\Lambda\,\cU^\T(\bsb,\bsbeta)\,\Lambda^\dag.
\eeq

This relation can be understood taking into account that the ORF on the real line and the
unit circle are related by the Cayley transform. More precislely, $\phi_n(z)$ are ORF on
the real line iff $\phi_n(\tilde\zeta(z))$ are ORF on the unit circle. If $\mu$ is the
orthogonality measure on $\overline\R$, the corresponding measure $\nu$ on $\T$ is given
by $\nu(\Delta)=\mu(\tilde\zeta(\Delta))$ for any Borel subset $\Delta$ of $\T$. Also,
the parameters $\alpha_n$ and $\beta_n$ associated respectively with the poles of
$\phi_n(z)$ and $\phi_n(\tilde\zeta(z))$ are related by $\beta_n=\zeta(\alpha_n)$.
Moreover, $\phi_n$ satisfies the analogue of recurrence (\ref{RR}) on the real line with
coefficients $a_n$ iff $\widehat\phi_n=\xi_0^2\xi_1^2\cdots\xi_{n-1}^2\xi_n\phi_n$
satisfies such a recurrence on the unit circle with coefficients
$b_n=\xi_0^2\xi_1^2\cdots\xi_{n-1}^2a_n$. If $\chi_n$ and $\widehat\chi_n$ are the
associated ORF (given by the corresponding version of (\ref{LORF}) on $\overline\R$ and
$\T$ respectively), then $\widehat\chi_n=\lambda_n\chi_n$ with $\lambda_n$ as in
(\ref{GAMMA-LAMBDA}). Therefore, if $\bsalpha$ is compactly included in $\U$, the matrix
representation $\cU^\R(\bsa,\bsalpha)$ of $T_\mu$ with respect to $(\chi_n)_{n\geq0}$ and
the matrix representation $\cU^\T(\bsb,\bsbeta)$ of $T_\nu$ with respect to
$(\widehat\chi_n)_{n\geq0}$ are related by (\ref{UR-UT}).

\bigskip

\noindent {\bf Acknowledgements}

This work was partially realized during a stay of the author at the Norwegian University
of Science and Technology financed by Secretar\'{\i}a de Estado de Universidades e
Investigaci\'{o}n from the Ministry of Education and Science of Spain. The work of the author
was also partly supported by a research grant from the Ministry of Education and Science
of Spain, project code MTM2005-08648-C02-01, and by Project E-64 of Diputaci\'on General
de Arag\'on (Spain).

The author is very grateful to Professor Olav Nj\aa stad for his hospitality at the
Norwegian University of Science and Technology and for his constant interest and
invaluable help during the development of this work. The discussions with him were
indispensable to find the right way to generalize to the rational case the matrix
approach to orthogonal polynomials on the unit circle.

The author also thanks Professor Barry Simon for his useful comments on operator M\"{o}bius
transformations.


\begin{thebibliography}{99}

\bibitem{Ak56}
N.I. Akhiezer,
\emph{Theory of Approximation},
Frederic Ungar Publ. Co., New York, 1956.

\bibitem{AkKr62}
N.I. Akhiezer, M.G. Krein,
\emph{Some Questions in the Theory of Moments},
Transl. Math. Monographs, Vol.2, AMS, Providence, RI, 1962;
Russian original, Kharkov, 1938.

\bibitem{Al77}
M. Alfaro,
\emph{El operador multiplicaci\'on en la teor\'{\i}a de polinomios ortogonales
sobre la circunferencia unidad},
Proc. II Spanish-Portuguese Mathematical Conference (Madrid, 1973), pp. 13--21,
Consejo Sup. Inv. Cient., Madrid, 1977.

\bibitem{AzIo89}
T.Ya. Azizov, I.S. Iokhvidov,
\emph{Linear operators in spaces with an indefinite metric},
John Wiley $\&$ Sons, Ltd., Chichester, 1989.

\bibitem{BaLo99}
D. Barrios, G. L\'{o}pez,
\emph{Ratio asymptotics for polynomials orthogonal on arcs of the unit circle},
Constr. Approx. {\bf 15} (1999) 1--31.

\bibitem{Be71}
I.D. Berg,
\emph{An extension of the Weyl-von Neumann theorem to normal operators},
Trans. Amer. Math. Soc. {\bf 160} (1971) 365--371.

\bibitem{Bi62}
M.S. Birman,
\emph{On existence conditions for wave operators},
Dokl. Akad. Nauk SSSR {\bf 143} (1962) 506--509 (Russian).

\bibitem{BiKr62}
M.S. Birman, M.G. Krein,
\emph{On the theory of wave operators and scattering operators},
Dokl. Akad. Nauk SSSR {\bf 144} (1962) 475--478 (Russian).

\bibitem{BGHN99a}
A. Bultheel, P. Gonz\'{a}lez-Vera, E. Hendriksen, O. Nj{\aa}stad,
\emph{A density problem for orthogonal rational functions},
J. Comput. Appl. Math. {\bf 105} (1999) 199-212.

\bibitem{BGHN99}
A. Bultheel, P. Gonz\'{a}lez-Vera, E. Hendriksen, O. Nj{\aa}stad,
\emph {Orthogonal rational functions},
Cambridge Monographs on Applied and Computational Mathematics, 5,
Cambridge University Press, Cambridge, 1999.

\bibitem{MPOP}
M.J. Cantero, L. Moral, L. Vel\'azquez,
\emph{Measures and para-orthogonal polynomials on the unit circle},
East J. Approx. {\bf 8} (2002) 447--464.

\bibitem{FIVE}
M.J. Cantero, L. Moral, L. Vel\'azquez,
\emph{Five-diagonal matrices and zeros of orthogonal polynomials on the unit circle},
Linear Algebra Appl. {\bf 362} (2003) 29--56.

\bibitem{MIN}
M.J. Cantero, L. Moral, L. Vel\'azquez,
\emph{Minimal representations of unitary operators and orthogonal polynomials on the unit circle},
Linear Algebra Appl. {\bf 408} (2005) 40--65.

\bibitem{TRUNC}
M.J. Cantero, L. Moral, L. Vel\'azquez,
\emph{Measures on the unit circle and unitary truncations of unitary operators},
J. Approx. Theory {\bf 139} (2006) 430--468.

\bibitem{Ge44}
Ya.L. Geronimus,
\emph{On polynomials orthogonal on the circle, on trigonometric moment problem,
and on allied Carath\'eodory and Schur functions},
Mat. Sb. {\bf 15} (1944) 99--130.

\bibitem{Go00}
L. Golinskii,
\emph{Singular measures on the unit circle and their reflection coefficients},
J. Approx. Theory {\bf 103} (2000) 61--77.

\bibitem{Go00b}
L. Golinskii,
\emph{Operator theoretic approach to orthogonal polynomials on an arc of the unit circle},
Matematicheskaya fizika, analiz, geometriya {\bf 7} (2000) 3--34.

\bibitem{Go00c}
L. Golinskii,
\emph{On the spectra of infinite Hessenberg and Jacobi matrices},
Matematicheskaya fizika, analiz, geometriya 7 (2000) 284--298.

\bibitem{GoNeAs95}
L. Golinskii, P. Nevai, W. Van Assche,
\emph{Perturbation of orthogonal polynomials on an arc of the unit circle},
J. Approx. Theory {\bf 83} (1995) 392--422.

\bibitem{Gr93}
W.B. Gragg,
\emph{Positive definite Toeplitz matrices, the Arnoldi process for isometric operators,
and Gaussian quadrature on the unit circle},
J. Comput. Appl. Math. {\bf 46} (1993) 183--198;
Numerical Methods of Linear Algebra, pp. 16--32, Moskov. Gos. Univ., Moskow, 1982.

\bibitem{JoNjTh89}
W.B. Jones, O. Nj\aa stad, W.J. Thron,
\emph{Moment theory, orthogonal polynomials, quadrature, and continued fractions
associated with the unit circle},
Bull. London Math. Soc. {\bf 21} (1989) 113--152.

\bibitem{Ka57}
T. Kato,
\emph{Perturbation of continuous spectra by trace class operators},
Proc. Japan. Acad. {\bf 33} (1957) 260--264.

\bibitem{Kr50}
M.G. Krein,
\emph{On an application of the fixed point principle in the theory of linear transformations
of spaces with an indefinite metric},
Uspehi Matem. Nauk (N.S.) {\bf 5} (1950), no. 2(36), 180--190 (Russian).

\bibitem{Kr64}
M.G. Krein,
\emph{A new application of the fixed-point principle in the theory of operators in a space
with indefinite metric},
Dokl. Akad. Nauk SSSR {\bf 154} (1964) 1023--1026 (Russian).

\bibitem{KrSm67}
M.G. Krein, Yu.L. \v Smuljan,
\emph{On linear-fractional transformations with operator coefficients},
Mat. Issled {\bf 2} (1967), no. 3, 64--96 (Russian);
English transl. in Amer. Math. Soc. Transl., Ser. 2, {\bf 103} (1974) 125--152.

\bibitem{LaSi06}
Y. Last, B. Simon,
\emph{The essential spectrum of Schr\"{o}dinger, Jacobi and CMV operators},
J. Anal. Math. {\bf 98} (2006) 183--220.

\bibitem{MaGo91}
F. Marcell\'an, E. Godoy,
\emph{Orthogonal polynomials on the unit circle: distribution of zeros},
J. Comput. Appl. Math. {\bf 37} (1991) 195--208.

\bibitem{ReSi72}
M. Reed, B. Simon,
\emph{Methods of Modern Mathematical Physics, Vol. 1: Functional Analysis},
Academic Press, New York-London, 1972.

\bibitem{Si71}
W. Sikonia,
\emph{The von Neumann converse of Weyl's theorem},
Indiana Univ. Math. J. {\bf 21} (1971/72) 121--124.

\bibitem{Si105}
B. Simon,
\emph{Orthogonal Polynomials on the Unit Circle, Part 1: Classical Theory},
AMS Colloquium Series, American Mathematical Society, Providence, RI, 2005.

\bibitem{Si205}
B. Simon,
\emph{Orthogonal Polynomials on the Unit Circle, Part 2: Spectral Theory},
AMS Colloquium Series, American Mathematical Society, Providence, RI, 2005.

\bibitem{Si}
B. Simon,
\emph{CMV matrices: Five years after},
to appear in the Proceedings of the W.D. Evans 65th Birthday Conference,
arXiv:math.SP/0603093, 2006.

\bibitem{Te92}
A.V. Teplyaev,
\emph{The pure point spectrum of random polynomials orthogonal on the unit circle},
Soviet Math. Dokl. {\bf 44} (1992) 407--411;
Dokl. Akad. Nauk SSSR 320 (1991) 49--53.

\bibitem{Th88}
W.J. Thron,
\emph{$L$-polynomials orthogonal on the unit circle},
Nonlinear numerical methods and rational approximation (Wilrijk, 1987),
pp. 271--278, Math. Appl., vol. 43, Reidel, Dordrecht, 1988.

\bibitem{We09}
H. Weyl,
\emph{\"{U}ber beschra\"{a}nkte quadratische Formen, deren Differenz vollstetig ist},
Rend. Circ. Mat. Palermo {\bf 27} (1909) 373--392.

\bibitem{Wa93}
D.S. Watkins,
\emph{Some perspectives on the eigenvalue problem},
SIAM Rev. {\bf 35} (1993) 430--471.

\end{thebibliography}
\end{document}